\documentclass[11pt]{article}
\usepackage{CJK}
\usepackage{graphicx}
\usepackage{fancyhdr}
\usepackage{amsmath,amsthm,amssymb}
\usepackage{bm}
\usepackage{mathrsfs}
\usepackage{stmaryrd}
\usepackage{multirow}
\usepackage{indentfirst}
\usepackage{url}
\usepackage{color}
\usepackage{epstopdf}

\newtheorem{theorem}{Theorem}[section]
\numberwithin{equation}{section}
\makeatletter
\newcommand{\Rmnum}[1]{\expandafter\@slowromancap\romannumeral #1@}
\makeatother

\allowdisplaybreaks
\def\hDash{\bot\!\!\!\bot}

\begin{document}

\title{\bf Integrated conditional moment test and beyond: when the number of covariates is divergent
\footnote{Lixing Zhu is a Chair professor of Department of Mathematics
at Hong Kong Baptist University, Hong Kong, China. He was supported by a grant from the
University Grants Council of Hong Kong, Hong Kong, China. }}%The authors thank the editor, the associate editor and two referees for their constructive suggestions that led to the  improvement of an early manuscript.}}
\author{Falong Tan$^1$ and Lixing Zhu$^2$ \\~\\
{\small {\small {\it $^1$ College of Finance and Statistics, Hunan University, Changsha, China} }}\\
{\small {\small {\it $^2$ Department of Mathematics, Hong Kong Baptist University, Hong Kong} }}}
\date{}
\maketitle

\begin{abstract}
The classic integrated conditional moment test is a  promising method for testing regression model misspecification. However, it severely suffers from the curse of dimensionality. To extend it to handle the testing problem for parametric multi-index models with diverging number of covariates, we investigate three issues in inference in this paper. First, we study the consistency and asymptotically linear representation of the least squares estimator of the parameter matrix at faster rates of divergence than those in the literature for nonlinear models. Second,   we propose, via sufficient dimension reduction techniques, an adaptive-to-model version of the integrated conditional moment  test. We study the asymptotic properties of the new test under both the null and alternative hypothesis to examine its ability of significance level maintenance and its sensitivity to the global and local alternatives that are distinct from the null at the fastest possible rate in hypothesis testing. Third, we derive the consistency of the bootstrap approximation for the new test in the diverging dimension setting.  The numerical studies show that the new test can very much enhance the performance of the original ICM test in high-dimensional scenarios. We also apply the test to a real data set for illustrations.  \\

\end{abstract}

\newpage
%\baselineskip=16pt

%\newpage

\setcounter{equation}{0}
\section{Introduction}
To test for model misspecification in regression analysis, integrated conditional moment (ICM) test (Bierens 1982; Bierens and Ploberger 1997) is a classic approach  that attracts, when  the dimension $p$ of the covariate vector is fixed,  many follow-up studies and develops many variants in the literature. The research described in this paper is motivated by extending it to handle the goodness-of-fit testing for parametric multiple-index regression models when the dimension $p$ of the covariate vector is large regarding it as a diverging number as the sample size $n$ goes to infinity.
%This is to prevent possible wrong conclusions in any further statistical analysis.
%Note that when $p$ is fixed, the relevant theories of inference have been very mature. However, in high dimensional data analysis, for large $p$ that could be treated as a divergent number when the sample size $n$ tends to infinity, such a model checking problem has not yet been well studied in the literature. To this end, we need to estimate the parameter of interest in the models and then to construct a test.
%An early revelent reference is the seminal paper of Neyman and Scott (1948). Huper (1973) consider robust estimates for linear regressions when the dimension $p$ of predictors grows to infinity with the sample size $n$ at a convergence rate $p=o(n^{1/3})$. Most existing revelent literatures focus on estimation problems for divergent $p$.
As the classic ICM test suffers severely from the curse of dimensionality that will be described below, we then try to solve three problems in inference. First, we will study at which rate of divergence of $p$, the asymptotically linear presentation of the least squares estimation of the parameters in the regression model can hold. Second, we develop an adaptive-to-model version of the ICM test such that it can greatly avoid the curse of dimensionality. Third, as the limiting null distribution is in general intractable, we investigate the consistency of bootstrap approximation in this diverging setting.  As everything in this paper relies on the sample size $n$, we stress $n$ in the notations.
%For contemporary literature, one can refer to Donoho (2000) and Fan and Li (2006).

%Dimension-reduction models such as parametric single-index models are often used in high dimensional regression analysis. Thus, in the present paper, we restrict ourselves to testing a parametric single-index model,
We begin with describing the model checking problem in detail. Consider a parametric multiple-index model
\begin{equation}\label{1.1}
Y = g(\bm{\beta}_0^{\top}X, \vartheta_0) + \varepsilon
\end{equation}
where $Y$ represents the real-valued response variable with the $p$-dimensional explanatory variable $X$, $g$ is a given smooth function, $\vartheta_0 \in \mathbb{R}^l$ and $\bm{\beta}_0=(\beta_{10}, \cdots, \beta_{d0}) \in \mathbb{R}^{p \times d}$ is the unknown matrix of regression parameters, $\varepsilon=Y-E(Y|X)$ is the unpredictable part of $Y$ given $X$, and the notation $\top$ denotes the transpose. Without loss of generality, we assume that $\beta_{10}, \cdots, \beta_{d0}$ are orthogonal.
%Single index models~(\ref{1.1}) always be attractive to practitioners because these models constitute a compromise between a parametric and a completely nonparametric models. Since we do not know the underlying model in advance, it is necessary to check the mis-specification of the regression models for further statistical inference to avoid misleading conclusion. Thus, a saturated alternative model
%\begin{equation}\label{1.2}
%Y = G(X) + \varepsilon
%\end{equation}
%is also considered where $G(\cdot)=E(Y|X=\cdot)$ is an unknown smooth function.

When the dimension $p$ of covariates is fixed, there are a number of proposals in the literature for testing the parametric regression models consistently which can also be applied to test the null hypothesis of model~(\ref{1.1}). We just list a few such as H\"{a}rdle and Mammen (1993), Zheng (1996), Dette (1999), Fan and Huang (2001), Horowitz and Spokoiny (2001), Koul and Ni (2004) and van Keilegom et al. (2008) by using nonparametic estimation methods. They are called local smoothing tests. Some of tests in this class  have tractable limiting null distributions, but  the use of nonparametric regression estimations usually cause them to suffer severely from the curse of dimensionality in high-dimensional cases. Guo et al. (2016) gave some detailed comments on this issue. The empirical process-based tests are constructed in terms of converting conditional expectations $E(\varepsilon|X)=0$ to infinite and parametric unconditional orthogonality restrictions, that is
\begin{equation}\label{1.2}
  E(\varepsilon|X)=0 \quad \Leftrightarrow \quad E[\varepsilon w(X,t)]=0 \quad \forall \ t  \in \Gamma,
\end{equation}
where $\Gamma$ is some proper space. There exist several parametric families $w(X,t)$ such that the equivalence~(\ref{1.2}) holds; see Bierens and Ploberger (1997) and Escanciano (2006b) for more details on the primitive conditions of $w(\cdot,t)$ to satisfy this equivalence.
The indicator function $w(X,t)=I(X \leq t)$ is commonly used as a weight function in the literature; see, e.g., Stute (1997), Stute, Gonz$\rm \acute{a}$lez Manteiga and Presedo Quindimil (1998), Stute, Thies and Zhu (1998), Zhu (2003), Khmadladze and Koul (2004), Stute, Xu and Zhu (2008), among many others. Bierens (1982) used the characteristic function $w(X,t)=\exp(\mathrm{i} t^{\top} X)$ as the weight function and constructed an integrated conditional moment (ICM) test, where $\mathrm{i}=\sqrt{-1}$ denotes the imaginary unit. The test statistic is given by
\begin{equation}\label{1.3}
ICM_n= \int_{\Gamma} |\frac{1}{\sqrt{n}} \sum_{j=1}^{n} \hat{\eta}_j \exp(\mathrm{i} t^{\top} \Phi(X_j))|^2 d\mu(t),
\end{equation}
where $\hat{\eta}_j$ is the residual, $\Gamma$ is a compact subset in $\mathbb{R}^p$ with non-empty interior, $\Phi(\cdot)$ is a bounded smoothing function, and  $\mu$ is a probability measure on $\Gamma$. Note that $ICM_n$ integrates over a compact subset $\Gamma \subset \mathbb{R}^p$. It is well known that high-dimensional numerical integrations are extremely difficult to handle in practice. When $p$ is large, the computation of the integral in~(\ref{1.3}) or its approximation becomes difficult,  time consuming and expensive. Escanciano (2006a) suggested to choose the standard normal measure in~(\ref{1.3}) and the high-dimensional numerical integrations can be avoided. The test statistic $ICM_n$ then becomes
\begin{equation}\label{1.4}
ICM_n = \int_{\mathbb{R}^p}|\frac{1}{\sqrt{n}} \sum_{j=1}^{n} \hat{\eta}_j \exp(\mathrm{i} t^{\top} X_j))|^2 \phi(t) dt = \frac{1}{n} \sum_{j,k=1}^{n} \hat{\eta}_j \hat{\eta}_k \exp(-\frac{1}{2}\|X_j-X_k\|^2),
\end{equation}
where $\phi(t)$ is the standard normal density on $\mathbb{R}^p$. Escanciano (2006a, 2006b) commented that Bierens' test (1982) is less impacted by the dimension $p$ as it is based on one-dimensional projections.  Although it does perform well in many cases,  we found that Bierens' test (1982) still suffers from the dimensionality problem in some other cases. More specifically,  when $p$ is large it has difficulty to maintain the significance level and is less powerful. We give a simple simulation study to show this phenomenon by the following toy example:
$$H_{1}:  Y = \beta_0^{\top}X+ a \exp(-(\beta_0^{\top}X)^2) +\varepsilon; $$
where $\beta_0=(1, \cdots, 1)^{\top}/\sqrt{p}$ and $X$ is $N(0, I_p)$ independent with the standard normal distributed error $\varepsilon$. The significance level is set to be $\alpha=0.05$. The constant $a=0$ corresponds to the null and $a \neq 0$ to the alternatives. The empirical sizes and powers of Bierens' test (1982) $ICM_n$ are presented in  Figure 1. We can find that not only the empirical powers drop quickly as the dimension $p$ increases, but also the empirical sizes are almost zero for large $p$. This means that the $ICM_n$ test does not work even when the dimension is moderate.
\begin{figure}
  \centering
  \includegraphics[width=16cm,height=5cm]{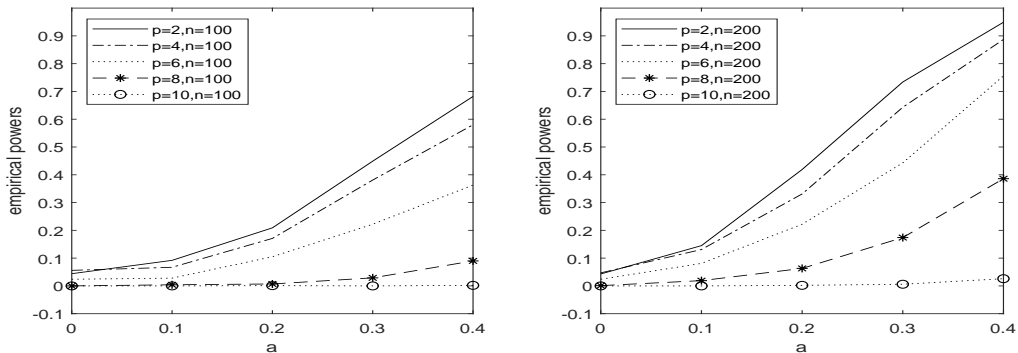}
  \caption{The empirical powers of $ICM_n$ for the model $H_1$.}\label{Figure1}
\end{figure}

To better understand this phenomenon in large dimension scenarios, we take a special case for illustration. Suppose that the components $\{X_{ik}: k=1, \cdots, p\}$ of $X_i$ are independent and identically distributed. Note that $\|X_i-X_j\|^2=\sum_{k=1}^p (X_{ik}-X_{jk})^2$. If the dimension $p$ tends to infinity, it follows that $(1/p)\|X_i-X_j\|^2$ tends to $2 Var(X_{i1})$. Thus $\exp\{-\|X_i-X_j\|^2/2\}$ tends to 0 with probability 1. This means that the Bierens' test statistic tends to be small under both the null and alternatives for large $p$.

%To alleviate the dimensionality problem, we use sufficient dimension reduction and projection technique to construct an adaptive-to-model ICM test.
%Motivated by such observations,
The purpose of this paper is to extend the Bierens' (1982) ICM test to the diverging dimension setting and simultaneously avoid the curse of dimensionality. As we need
the residuals $Y_j- g(\bm{\hat{\beta}}^{\top}X_j, \hat \vartheta)$ to construct the test statistic where $ \bm{\hat{\beta}}$ and $\hat \vartheta$ are estimators of $\bm{{\beta}}_0$ and $ \vartheta_0$ respectively, we first study the asymptotic properties of the estimator $(\bm{\hat{\beta}}, \hat \vartheta)$ under both the null and alternatives in the divergent dimension cases. This problem has been well investigated for linear regression models, see Huber (1973), Portnoy (1984, 1985) and Zou and Zhang (2009), etc. Nevertheless, for nonlinear models, it has rarely been studied in the literature.
Under some wild conditions, Tan and Zhu (2019a) obtained the estimation consistency and  the asymptotically linear representation of $(\bm{\hat{\beta}}, \hat \vartheta)$ at the rates $p=o(n^{1/4}) $  and $p=o(n^{1/5})$ of divergence  respectively. In this paper, with the help of some high dimensional empirical processes techniques developed in Tan and Zhu (2019b), we will greatly improve these rates of divergence to $p=o(n^{1/2}) $  and $p=o(n^{1/3}/\log n)$ and obtain the uniformly asymptotically linear representation of $(\bm{\hat{\beta}}, \hat \vartheta)$. These improvements enhance our ability to handle higher dimension in model checking problems. Under the conditions we design, the rates are the best in the relevant research in the literature to the best of our knowledge. Shi et al. (2019) obtained similar results for generalized linear models with fixed design in the diverging setting. But, as we actually deal with general nonlinear models with random design here, their approach is significantly different from ours.

To attack the curse of dimensionality, we notice that the main reason that the Bierens' (1982) test statistic $ICM_n$ does not work is that the weight function $\exp(-\|X_j-X_k\|^2 /2)$ in (\ref{1.4}) degenerates to zero at an exponential rate when the dimension $p$ tends to infinity.
To avoid this problem, we adopt projection pursuit and sufficient dimension reduction (SDR) techniques to reduce the original dimension $p$ to a much low dimension and then   develop an adaptive-to-model version of the classic $ICM_n$ test under the above convergence rates. % possibly optimal we guess, otherwise, the test would have no weak limit. %Next we introduce some definitions and notations used in SDR theory.
The asymptotic properties of the proposed test are investigated under both the null and alternative hypothesis when $p=o(n^{1/3}/\log{n})$. We also show that the test is consistent against all global alternatives and can detect the local alternatives converging to the null at the parametric rate $1/\sqrt{n}$ in this diverging dimensional scenario. This means that the test is sensitive to the local alternatives at a parametric rate even when the dimension $p$ is divergent. As the limiting null distribution of the proposed test admits a complicated structure, its critical values generally cannot be tabulated in practice. Thus we propose a bootstrap method to approximate the asymptotic null distribution of the test. We also show that the consistency of the bootstrap approximation remains valid when $p=o(n^{1/3}/\log{n})$. This extends the results of Stute, Gonz$\rm \acute{a}$lez Manteiga and Presedo Quindimil (1998) to the diverging dimension cases.

To implement the adaptive-to-model test construction, we note that  model (\ref{1.1}) under the null admits a dimension reduction structure. Thus
to make full use of the model structure under both the null and the alternative hypothesis,
consider the following general alternative model
\begin{equation}\label{1.5}
Y = G( \bm{B}^{\top}X) + \varepsilon,
\end{equation}
where $E(\varepsilon|X)=0$, $G(\cdot)$ is an unknown smooth function, and $\bm{B}$ is a $p \times q$ orthonormal matrix such that $\mathcal{S}_{E(Y|X)}={\rm span}(\bm{B})$. Here $\mathcal{S}_{E(Y|X)}$ is the central mean space of $Y$ with respect to $X$ which will be specified in Section 3.
Note that this model is even  more general  than the nonparametric model $Y=G(X)+\varepsilon$ that is a special case when $\bm{B}$ is an $p \times p$ identity matrix with $q=p$.

The rest of this article is organised as follows. In Section~2 we present the asymptotic results of the least squares estimation of the parameters under the null and alternative models. Section~3 describes the basic test construction. Since the sufficient dimension reduction techniques are crucial for the model adaptation property, we give a short review in this section. To alleviate the computational complexity of the test statistic, we also discuss the choice of the weight functions in this section. Section~4 is devoted to present the asymptotic properties of the test under both the null  and  alternatives. In Section~5, the simulation results  are reported and a real data analysis is used as an illustration of application. Appendix contains the technical proofs for the  theoretical results.

\section{Parametric estimation}
In this section, we consider parameter estimation for the model~(\ref{1.1}), and also study its asymptotic properties in the diverging dimension setting.
%when the underlying model is actually the alternative model~(\ref{1.5}) %\subsection{Parametric estimation}
%Recall that the purpose of this paper is to check whether the regression function $m(\cdot)=E(Y|X=\cdot)$ belongs to a given parametric family $\mathcal{G}=\{g(\bm{\beta}^{\top} \cdot, \vartheta): \bm{\beta} \in \mathbb{R}^{p\times d} \ {\rm and} \ \vartheta \in \mathbb{R}^l \}$ when the dimension $p$ diverges as the sample size $n$ tends to infinity and $d$ and $l$ are fixed. Thus the null hypothesis we want to test is that
Recall the null model~(\ref{1.1})
\begin{eqnarray*}
 Y= g(\bm{\beta}_0^{\top} X, \vartheta_0)+ \varepsilon.
\end{eqnarray*}
To estimate the unknown parameters $(\bm{\beta}_0, \vartheta_0)$, we here restrict ourselves to the ordinary least squares method.
For  notational simplicity, define $\theta=[vec(\bm{\beta})^{\top},\vartheta^{\top}]^{\top}$ where $vec(\bm{\beta})=[\beta_1^{\top}, \cdots, \beta_d^{\top}]^{\top}$. Then
we rewrite the null hypothesis $H_0$ in the following parametric regression format:
\begin{eqnarray*}
Y= g(\theta_0, X)+ \varepsilon  \quad {\rm for\ some } \ \theta_0  \in \mathbb{R}^m,
\end{eqnarray*}
where $m=pd+l$, $E(\varepsilon|X)=0$ and the number $p$ diverges as the sample size $n$ tends to infinity. Next we consider the ordinary least squares estimator of $\theta_0$ and its asymptotic properties.  Let %In model~(\ref{1.1}). Let
\begin{eqnarray}\label{2.1}
\hat{\theta} = \mathop{\rm arg\ min}\limits_{\theta \in \mathbb{R}^m} \sum_{i=1}^{n}[Y_i-g(\theta, X_i)]^2.
\end{eqnarray}
To analyse the asymptotic properties of $\hat{\theta}$, set
\begin{eqnarray}\label{2.2}
\tilde{\theta}_0 = \mathop{\rm arg\ min}\limits_{\theta \in \mathbb{R}^m} E[Y-g(\theta, X)]^2.
\end{eqnarray}
%Here $\tilde \theta_0=\theta_0$ under the null model~(\ref{1.1}) and $\tilde \theta_0 \not =\theta_0$ under the alternative model~(\ref{1.4}).
To obtain the asymptotic properties of $\hat{\theta}$ when the dimension of $\theta$ is divergent, we need some regularity conditions. In the following, $C$ always stands for a constant which may be different from place to place.

(A1). Suppose the function $g(\theta, x)$ admits second derivatives in $\theta$. Let
\begin{eqnarray*}
\psi(\theta,x,y)=[y-g(\theta, x)]g'(\theta, x)=[y-g(\theta, x)] \frac{\partial g(\theta, x)}{\partial \theta},\\
\psi_j(\theta,x,y)=[y-g(\theta, x)]g_j'(\theta, x)=[y-g(\theta, x)] \frac{\partial g(\theta, x)}{\partial \theta_j}.
\end{eqnarray*}
Assume that (i) $E|\psi_j(\tilde{\theta}_0,X,Y)|^2 < C$ for all $j=1, \cdots, m$;
(ii) $|\frac{\partial \psi_j}{\partial \theta_i}(\theta, x,y)| < F(x,y)$ with $E|F(X,Y)|^2 < \infty$ for all $i,j=1, \cdots, m$ and $\theta \in U(\tilde{\theta}_0)$, where $U(\tilde{\theta}_0)$ is some neighborhood of $\tilde{\theta}_0$.
%$\psi_j(\theta,x,y)$ satisfies the Lipschitz condition in some neighborhood $U(\tilde{\theta}_0)$ of $\tilde{\theta}_0$, that is
%$$|\psi_j(\theta_1,x,y)-\psi_j(\theta_2,x,y)| \leq \|\theta_1-\theta_2\| F(x,y), $$
%where $\theta_1, \theta_2 \in U(\tilde{\theta}_0)$.
%and the moment generating function of the random variable $F(X,Y)$ exists around zero.

(A2). Let $P \psi_{\theta}=E[\psi(\theta, X, Y)]$ and $P \psi_{j\theta}=E[\psi_j(\theta, X, Y)]$. Suppose that $P \psi_{\theta}$ is twice differentiable in the neighborhood $U(\tilde{\theta}_0)$ of $\tilde{\theta}_0$. Set
$$ \dot{P}\psi_{\theta}=\frac{\partial P\psi_{\theta}}{\partial \theta},  \ \ddot{P}\psi_{\theta}=\frac{\partial P\psi_{\theta}}{\partial \theta^{\top} \partial \theta}, \ {\rm and} \ \ddot{P}\psi_{j \theta}=\frac{\partial P\psi_{j \theta}}
{\partial \theta^{\top} \partial \theta}. $$
Define $\Sigma=-\dot{P}\psi_{\tilde{\theta}_0}$. Assume that the matrix $\Sigma$ is nonsingular and satisfies the following condition
$$
0<\underline{\lambda} \leq \lambda_{\min}(\Sigma) \leq \lambda_{\max}(\Sigma) <\overline{\lambda} < \infty,
$$
where $\underline{\lambda}$ and $\overline{\lambda}$ are two constants that do not depend on $n$ and $\lambda_{\min}(\Sigma)$ and $\lambda_{\max}(\Sigma)$ denote the smallest and largest eigenvalue of $\Sigma$ respectively. Furthermore, let $\lambda_1(\ddot{P}\psi_{j \theta}), \cdots, \lambda_m(\ddot{P}\psi_{j \theta})$ be the eigenvalues of the matrix $\ddot{P}\psi_{j \theta}$ and assume
\begin{eqnarray*}
\max_{1\leq i,j \leq m} |\lambda_i(\ddot{P}\psi_{j \theta})|\le C, \ \forall \ \theta  \in U(\tilde{\theta}_0)
\end{eqnarray*}
where $C$ is a constant that is free of $n$.

(A3). The vector $\tilde{\theta}_0$ is the unique minimizer of~(\ref{2.2}). Thus it is easy to see that $\tilde{\theta}_0=\theta_0$ under the model~(\ref{1.1}).

The regularity conditions (A1) and (A3) are standard for the nonlinear least squares estimation, see, e.g., Jennrich (1969) and White (1981). (A2) is similar to the regularity condition on the information matrix $I_n$ of the likelihood function proposed in Fan and Peng (2004) which facilitated the theoretical derivations. Although we only need the asymptotic property of the parameters for multiple index models in this paper, it also holds for general parametric regression models.

\begin{theorem}\label{theorem2.1}
Suppose the regularity conditions (A1)-(A3) hold. If $p^2/n \to 0$, then it follows that $\hat{\theta}$ is a norm consistent estimator of $\tilde{\theta}_0$ in the sense that $\|\hat{\theta}-\tilde{\theta}_0\|=O_p(\sqrt{p/n})$ where $\|\cdot\|$ denotes the Frobenius norm.
\end{theorem}

For the asymptotically linear representation of $\hat{\theta}-\tilde{\theta}_0$, we need some extra regularity conditions. We say a random variable $X \in \mathbb{R}$ is sub-Gaussian with variance proxy $\sigma^2$ if
$$ E\{\exp[u(X-EX)]\} \leq \exp(\frac{\sigma^2 u^2}{2}), \ \forall \ u \in \mathbb{R}. $$

(A4). The function $\psi_j(\theta,x)$ admits second derivatives in $\theta$ and $|\frac{\partial^2 \psi_j}{\partial \theta_k \partial \theta_l}(\theta, x,y)| < F(x,y)$ with $E|F(X,Y)|^2 < \infty$ for all $\theta \in U(\tilde{\theta}_0)$ and $i,j=1, \cdots, m$. Moreover, $v^{\top} \frac{\partial \psi}{\partial \theta}(\tilde{\theta}_0, X, Y)v$ is sub-Gaussian with variance proxy $\sigma^2$ for all $v \in \mathcal{S}^{m-1}$, where $\mathcal{S}^{m-1}$ is the unit sphere in $\mathbb{R}^m$.

The ``sub-Gaussian" condition is also used in Cai et al. (2010) to obtain the optimal rates of convergence for large covariance matrix estimation. It seems too restrictive in high dimensional settings. However, in sufficient dimension reduction, it can be rather reasonable in the following sense. When $p \to \infty$ as $n \to \infty$, Hall and Li (1993) proved that the linear combinations of the covariates are approximately normally distributed. Consequently, in the high dimension setting, the random variable $v^{\top} \frac{\partial \psi}{\partial \theta}(\tilde{\theta}_0, X, Y)v$ approximately has a Gaussian tail that is equivalent to the fact that  $v^{\top} \frac{\partial \psi}{\partial \theta}(\tilde{\theta}_0, X, Y)v$ is approximately sub-Gaussian. See Section 2 of Wainwright (2019) for more details.

\begin{theorem}\label{theorem2.2}
Suppose the regularity conditions (A1)-(A4) hold. If $(p\log{n})^3/n \to 0$, it follows that
\begin{equation}\label{2.3}
\sqrt{n} \alpha^{\top}(\hat{\theta}-\tilde{\theta}_0)=\frac{1}{\sqrt{n}} \sum_{i=1}^{n}[Y_i-g(\tilde{\theta}_0, X_i)] \alpha^{\top} \Sigma^{-1} g'(\tilde{\theta}_0, X_i)+o_p(1),
\end{equation}
where the remaining term $o_p(1)$ is uniformly in $\alpha \in \mathcal{S}^{m-1}$.
\end{theorem}

The convergence rate $O_p(\sqrt{p/n})$ in Theorem~\ref{theorem2.1} is in line with the results of the $M$-estimator studied by Huber (1973) and Portnoy (1984) when the number of parameters is divergent. But they only considered linear regression models $g(\theta, X)=\theta^{\top}X$.
In Theorem~\ref{theorem2.2} we obtain the uniformly asymptotically linear representation of $\hat{\theta}-\tilde{\theta}$. This result is much stronger than that in Portnoy (1985) who also only considered linear regression models. Furthermore, in the linear regression setting, we can derive, when $p=o(n)$ and $p=o(n^{1/2})$, the norm consistency and the uniformly asymptotically linear representation of $\hat{\theta}$ to $\tilde{\theta}_0$.

\begin{theorem}\label{theorem2.3}
Suppose the regularity conditions (A1)-(A3) hold for linear regression models and $v^{\top}X$ is sub-Gaussian with the variance proxy $\sigma^2$ for all $v \in \mathcal{S}^{p-1}$. If $p/n \to 0$, then we have
$\|\hat{\theta}-\tilde{\theta}_0\|=O_p(\sqrt{p/n})$ where $\|\cdot\|$ denotes the Frobenious norm. Moreover, if $p^2/n \to 0$, it follows that
$$\sqrt{n}\alpha^{\top}(\hat{\theta}-\tilde{\theta}_0)=\frac{1}{\sqrt{n}}
\sum_{i=1}^{n}(Y_i-\tilde{\theta}_0^{\top}X_i) \alpha^{\top} \Sigma^{-1}X_i+o_p(1),$$
where $\Sigma=E(XX^{\top})$ and the term $o_p(1)$ is uniformly in $\alpha \in \mathcal{S}^{p-1}$.
\end{theorem}

\section{Adaptive-to-model version of the ICM test}
%Recall that the purpose of this paper is to check whether the regression function $m(\cdot)=E(Y|X=\cdot)$ belongs to a given parametric family $\{ g(\beta_1^{\top} \cdot, \cdots, \beta_d^{\top} \cdot): \beta_1, \cdots, \beta_d \in \mathbb{R}^p \}$ when the dimension $p$ diverges as the sample size $n$ tends to infinity and $d$ is fixed.
Now we restate the null hypothesis as
\begin{eqnarray*}
H_0:\ E(Y|X)= g(\bm{\beta}_0^{\top} X, \vartheta_0) \quad {\rm for \ some} \ \bm{\beta}_0 \in \mathbb{R}^{p \times d} \ {\rm and} \ \vartheta_0 \in \mathbb{R}^l .
\end{eqnarray*}
While the alternative hypothesis is that for any $\bm{\beta} \in \mathbb{R}^{p \times d}$ and $\vartheta \in \mathbb{R}^l$,
\begin{eqnarray*}
H_1: \ E(Y|X)= G(\bm{B}^{\top} X) \neq g(\bm{\beta}^{\top} X, \vartheta),
\end{eqnarray*}
where $G$ is an unknown smoothing function, $\bm B$ is a $p\times q$ orthonormal matrix satisfying $\mathcal{S}_{E(Y|X)}= {\rm span}(\bm{B})$.
Here $\mathcal{S}_{E(Y|X)}$ is the central mean space of $Y$ with respect to $X$ which is defined as the intersection of all subspaces ${\rm span}(\bm{A})$ such that $Y \hDash E(Y|X)|\bm{A}^{\top}X$ where $\hDash$ means statistical independence and ${\rm span}(\bm{A})$ means the subspace spanned by the columns of $\bm{A}$. Under mild conditions, such a subspace $\mathcal{S}_{E(Y|X)}$ always exists (see Cook and Li (2002)). The dimension of $\mathcal{S}_{E(Y|X)}$ is called the structural dimension which is $q$ under the alternatives. Similarly, under the null we have $\mathcal{S}_{E(Y|X)}= {\rm span}(\bm{ B})={\rm span}(\bm{\beta_0})$ with a structural dimension $d$. Here we use the same notation $\bm B$ under both the null and alternatives which should be different in these two cases.

%For simplicity we define $\bm{ \beta_0}=(\beta_{01}, \cdots, \beta_{0q})$, $\bm{\beta}=(\beta_1, \cdots, \beta_d)$

Recall that $\varepsilon=Y-E(Y|X)$ and
$$(\bm{\tilde{\beta}_0}, \tilde{\vartheta}_0) = \mathop{\rm arg\ min}\limits_{\bm{\beta} \in \mathbb{R}^{p\times d}, \vartheta \in \mathbb{R}^l} E[Y-g(\bm{\beta}^{\top}X, \vartheta)]^2.$$
%as defined in Zheng (1996). We assume throughout that $\tilde{\beta}_0 \in \mathcal{S}_{E(Y|X)}$ under both the null and the alternatives  where $\mathcal{S}_{E(Y|X)}$ is the central mean subspace of $Y$ with respect to $X$. This is obvious for the null hypothesis. Under the alternative hypothesis, if $g(\beta^{\top}X, \theta)=\beta^{\top}X$ follows a linear regression model, it is easy to see $\tilde{\beta}_0=(EXX^{\top})^{-1}EXY \in \mathcal{S}_{E(Y|X)}.$ For other situations, it would typically depend on the form of $g(\cdot, \cdot)$.
Let $\eta=Y-g(\bm{\tilde{\beta}_0}^{\top}X, \tilde{\vartheta}_0)$ and $\bm{\tilde{\mathfrak{B}}}=(\bm{\tilde{\beta}_0}, \bm{B})$. Under the null hypothesis, we have $\bm{\tilde{\beta}_0}=\bm{\beta_0}$ and ${\rm span}(\bm{B})=\rm{span}(\bm{\beta_0})$.
Consequently,
$$E(\eta|\bm{\tilde{\mathfrak{B}}}^{\top}X)=E(Y|\bm{\tilde{\beta}_0}^{\top}X, \bm{B}^{\top}X)-g(\bm{\tilde{\beta}_0}^{\top}X, \tilde{\vartheta}_0)
=E(Y|\bm{\beta_0}^{\top}X)-g(\bm{\beta_0}^{\top}X, \vartheta_0)=0.$$
Under the alternative, we have $\eta=G(\bm{B}^{\top}X)-g(\bm{\tilde{\beta}_0}^{\top}X , \tilde{\vartheta}_0)+\varepsilon$. Then we obtain that
$$E(\eta|\bm{\tilde{\mathfrak{B}}}^{\top}X)=G(\bm{B}^{\top}X)-g(\bm{\tilde{\beta}_0}^{\top}X, \tilde{\vartheta}_0) \neq 0.$$
By Theorem 1 of Bierens (1982), under the null, we have $E[\eta \exp(\mathrm{i} t^{\top} \bm{\tilde{\mathfrak{B}}}^{\top}X)] = 0$ for all $t \in \mathbb{R}^{2d}$. Then it follows that
\begin{equation}\label{3.1}
\int_{\mathbb{R}^{2d}} |E[\eta \exp(\mathrm{i} t^{\top} \bm{\tilde{\mathfrak{B}}}^{\top}X)]|^2 \varphi(t) dt =0.
\end{equation}
While under the alternative, there exist $t_0 \in \mathbb{R}^{d+q}$ such that $E[\eta \exp(\mathrm{i} t_0^{\top} \bm{\tilde{\mathfrak{B}}}^{\top}X)] \neq 0$. Therefore,
\begin{equation}\label{3.2}
\int_{\mathbb{R}^{d+q}} |E[\eta \exp(\mathrm{i} t^{\top} \bm{\tilde{\mathfrak{B}}}^{\top}X)]|^2 \varphi(t) dt >0,
\end{equation}
where $\varphi(t)$ denotes a positive weight function which will be specified later.
Thereby we reject the null hypothesis for ``large value'' of the empirical version of the left-hand side of (\ref{2.1}). Let $\{ (X_1, Y_1), \cdots, (X_n, Y_n) \}$ be a random sample from the distribution of $(X, Y)$. Then we propose an adaptive-to-model integrated conditional moment test statistic as
\begin{equation}\label{3.3}
  \hat{V}_n(t) = \frac{1}{\sqrt{n}} \sum_{j=1}^{n} [Y_j-g(\bm{\hat{\beta}}^{\top}X_j, \hat{\vartheta}_0)] \exp(\mathrm{i} t^{\top} \bm{\hat{\mathfrak{B}}}^{\top}X_j),
\end{equation}
\begin{equation}\label{3.4}
% \hat{V}_n=\int_{\hat{\Gamma}} |\hat{V}_n(t)|^2 d \mu(t),
  AICM_n=\int_{\mathbb{R}^{d+\hat{q}}} |\hat{V}_n(t)|^2 \varphi(t) dt,
\end{equation}
%where $\hat{\Gamma}$ is some compact subset of $\mathbb{R}^{\hat{q}}$,
where $\bm{\hat{\mathfrak{B}}}=(\bm{\hat{\beta}}, \bm{\hat{B}})$, $\bm{\hat{B}}$ is a sufficient dimension reduction estimator of $\bm B$ with an estimated structural dimension $\hat{q}$ of $q$, and $\bm{\hat{\beta}}$ is a norm consistent estimators of $\bm{\tilde{\beta}_0}$. In this paper, we restrict ourselves to the ordinary least square estimator of $\bm{\tilde{\beta}_0}$.

%and $\varphi(t)$ is a positive function such that $\varphi(t)=\varphi(-t)$. %

%Our new test statistic is closely related to that of Bierens (1982) who proposed an Integrated Conditional Moment (ICM) test as follows:
%Recall that Bierens's (1982) original $ICM_n$ test statistic is given by
%$$ ICM_n= \int_{\Gamma} |\frac{1}{\sqrt{n}} \sum_{j=1}^{n} [Y_j-g(\hat{\beta}_n^{\top}X_j, \hat{\theta}_n)] \exp(i t^{\top} \Phi(X_j))|^2 d\mu(t). $$
%Since $\Gamma$ is a compact subset in $\mathbb{R}^p$, when $p$ is large, computation of the integral over $\Gamma$ become much more difficult. In contrast, we utilize the dimension reduction technique to reduce the original dimension $p$ to $d+s$. when $d+s$ is much smaller than the original dimension $p$, computational complexity of $AICM_n$ will be largely alleviated. We will show the advantages in the simulations.

It is worth mentioning that Guo et al. (2016) first used sufficient dimension reduction techniques to construct a goodness of fit test for parametric single-index models when the dimension $p$ of the covariates is fixed. But they only used the matrix $\bm B$ rather than $(\bm{\tilde{\beta}_0}, \bm{B})$ as we consider here. To make sure $E(\eta|{\bm B}^{\top}X) \neq 0$ under the alternative hypotheses, they need an extra condition that $\tilde{\beta}_0 \in {\rm span}(\bm B)$. Note that this extra condition does not always hold. Thus we use the matrix $\bm{\tilde{\mathfrak{B}}}$ in the conditional distribution to avoid this extra condition.

Another test related to ours is that of Stute and Zhu (2002), who developed a dimension reduction type test for generalized linear models in the fixed dimension cases. Their test statistic is based on the fact that $E(\eta|\beta_0^{\top}X)=0$ under the null hypothesis, where $\beta_0$ is the fixed direction involved in the generalized linear model. Stute and Zhu's (2002) test has been turned out to be very powerful for large $p$ in some cases. However, it is not omnibus as it only used one projection direction in the test statistic. Escanciano (2006a) gave a detail comment on this issue. While our test combines these two methods of Guo et al. (2016) and Stute and Zhu (2002) to avoid their shortcomings at the cost of dealing with one more dimension than that in Guo et al. (2016)'s test.

\subsection{The choice of $\varphi(t)$}
%Note that our test statistics $AICM_n$ involve a numerical integral in the space $\mathbb{R}^{d+\hat{q}}$.
The choice of weight functions $\varphi(t)$ is flexible.
%Bierens (1982) first suggested $w(t)=I(t \in \Gamma)$ where $\Gamma=[-\tau_1, \tau_1] \times \cdots \times [-\tau_p, \tau_p]$ (see, Bierens (1982), 109). However, as pointed by Bierens (1982), there are cases where $E[\eta \exp(\mathrm{i} t_0^{\top}B^{\top}X)] \neq 0$ only holds for $t_0$ far apart from $0$ and then $t_0$ may lie outside of the subset $\Gamma$. The second part of Lemma 1 avoid this question. So Bierens (1982) suggested a bounded one-to-one measurable mapping from $\mathbb{R}^p$ to $\mathbb{R}^p$, for instance, $ \Phi(x)=(\arctan(x_1), \cdots, \arctan(x_p))^{\top}$.
%Then Bierens (1982)'s original $ICM_n$ test statistic is given by
%$$ ICM_n= \int_{\Gamma} |\frac{1}{\sqrt{n}} \sum_{j=1}^{n} [Y_j-g(\hat{\beta}_n^{\top}X_j, \hat{\theta}_n)] \exp(\mathrm{i} t^{\top} \Phi(X_j))|^2 dt. $$
%Here we take $\tau_i=1/2$ for $i=1, \cdots, p$. Similarly, the corresponding adaptive integrated conditional moment test statistic is defined as
%$$ \hat{V}_{U,n}= \int_{\hat{\Gamma}} |\frac{1}{\sqrt{n}} \sum_{j=1}^{n} [Y_j-g(\hat{\beta}_n^{\top}X_j, \hat{\theta}_n)] \exp(\mathrm{i} t^{\top} \Phi(\hat{B}_n^{\top} X_j))|^2 dt. $$
%Since $\Gamma$ is a compact subset in $\mathbb{R}^p$, the computation of the integral over $\Gamma$ become much more difficult for large $p$. Even when we utilize the dimension reduction technique to reduce the original dimension $p$ to the structure dimension $q$, the calculation of the integral is still a big issue for large $q$. To reduce the computational complexity, we consider other choice of the weight functions $w(t)$.
Recall that our test statistic $AICM_n=\int_{\mathbb{R}^{d+\hat{q}}} |\hat{V}_n(t)|^2 \varphi(t) dt$.  Suppose that $\varphi(t)=\varphi(-t)$ and then some elementary calculations lead to
\begin{eqnarray*}
AICM_n
&=&\frac{1}{n} \int_{\mathbb{R}^{d+\hat{q}}} |\sum_{j=1}^{n} \hat{\eta}_j \exp(\mathrm{i} t^{\top} \bm{\hat{\mathfrak{B}}}^{\top}X_j)|^2 \varphi(t) dt \\
&=& \frac{1}{n} \sum_{j,k=1}^{n} \hat{\eta}_j \hat{\eta}_k \int_{\mathbb{R}^{d+\hat{q}}} \exp\{\mathrm{i} t^{\top} (\bm{\hat{\mathfrak{B}}}^{\top} X_j-\bm{\hat{\mathfrak{B}}}^{\top}X_k)\}\varphi(t)dt\\
&=& \frac{1}{n} \sum_{j,k=1}^{n} \hat{\eta}_j \hat{\eta}_k \int_{\mathbb{R}^{d+\hat{q}}} \cos\{ t^{\top} (\bm{\hat{\mathfrak{B}}}^{\top} X_j- \bm{\hat{\mathfrak{B}}}^{\top}X_k)\} \varphi(t) dt
\end{eqnarray*}
where $\hat{\eta}_j=Y_j-g(\bm{\hat{\beta}}^{\top}X_j, \hat{\vartheta}_j)$.
%The following idea is similar to Henze et al. (2005) and Zden$\rm \check{e}$k Hl$\rm \acute{a}$vka et al. (2017).
Let $K_{\varphi}(x)=\int \cos(t^{\top}x)\varphi(t)dt$, then it follows that
$$ AICM_n=\frac{1}{n} \sum_{j,k=1}^{n} \hat{\eta}_j \hat{\eta}_k K_{\varphi} (\bm{\tilde{\mathfrak{B}}}^{\top} X_j- \bm{\tilde{\mathfrak{B}}}^{\top}X_k). $$
%Let $w_a(x)=w(ax)$ be the weight function with $a>0$. Then the test statistic becomes
%$$ \hat{V}_{n,a}=\frac{1}{na^{\hat{q}}} \sum_{j,k=1}^{n} \hat{\eta}_j \hat{\eta}_k K_w (\frac{\hat{B}_n^{\top} X_j- \hat{B}_n^{\top}X_k}{a} ). $$
Yet the calculation of $K_{\varphi}(x)$ is still complex in high dimension settings, even when we use sufficient dimension reduction techniques here.
To facilitate the calculation of the test statistic,  a close form of the function $K_{\varphi}(x)$ is preferred.
There is a large class of weight functions $\varphi(t)$ available for this purpose. For instance, if we choose the density function $\phi(t)$ of a standard Gaussian distribution as a weight function, then
%Note that the characteristic function of a  spherically symmetric stable distribution (see, Vladimir V. Uchaikin and Vladimir M. Zolotarev (1999)) is given by
$$ K_{\phi}(x)=\int_{\mathbb{R}^{d+\hat{q}}} \cos(t^{\top}x)\phi(t)dt=\exp(- \frac{1}{2} \|x \|^2), $$
where $\|\cdot\|$ denotes the Frobenius norm.
%The multivariate standard Gaussian and Cauchy distributions are the spacial cases of spherical symmetric stable distributions with $s=2$ and $s=1$ respectively.
Consequently, the test statistic can be stated as
%$$ \hat{V}_{n,a}=\frac{1}{na^{\hat{q}}} \sum_{j,k=1}^{n} \hat{e}_j \hat{e}_k \exp(- \| \frac{ \hat{B}_n^{\top} X_j- \hat{B}_n^{\top}X_k}{a} \|^s ). $$
\begin{equation}\label{3.5}
AICM_n=\frac{1}{n} \sum_{j,k=1}^{n} \hat{\eta}_j \hat{\eta}_k \exp\{- \frac{1}{2} (\| \bm{\hat{B}}^{\top} X_j- \bm{\hat{B}}^{\top}X_k \|^2+\| \bm{\hat{\beta}}^{\top} X_j- \bm{\hat{\beta}}^{\top}X_k \|^2  ) \}.
\end{equation}
Note that the characteristic function of a spherically symmetric distribution has a form of $\psi(\|x\|)$ (see, Chapter 2 of Fang, Kotz and Ng (1990)). Thus the density function of any spherically symmetric distribution is also suitable as a weight function. In the simulation studies, we will use the density function of a multivariate standard Gaussian distribution.

Unlike Bierens's (1982) original test statistic $ICM_n$ in~(\ref{1.3}), we here utilize the sufficient dimension reduction and projection pursuit techniques to reduce the original dimension $p$ to $d+\hat{q}$. Then the exponential weight in (\ref{3.5}) would not deteriorate to zero as $p$ tends to infinity. In practice, when $d+\hat{q}$ is much smaller than $p$, the dimensionality difficulty will be largely alleviated.  We will show the advantages in the simulations.

%As pointed by Zden$\rm \check{e}$k Hl$\rm \acute{a}$vka et al. (2017), the choice of the weight function $w(x)$ is much less important than the choice of the parameter $a$. This is similar to the case of the nonparametric estimation where the kernel function is much less important than the choice of bandwidth. Therefore, in the simulation studies, we will use the density function of multivariate standard Gaussian distributions as the weight functions.

\subsection{Model adaptation property}
To achieve the model adaptation, we need the sufficient dimension reduction (SDR) techniques to identify the structural dimension $q$ and the central mean subspace $\mathcal{S}_{E(Y|X)}$. When $p$ is fixed, there are several methods in the literature to identify the central mean subspace $\mathcal{S}_{E(Y|X)}$, such as principal Hessian directions (pHd, Li (1992)). However, when $p$ is divergent, we have no corresponding asymptotic results about the estimated structural dimension $\hat q$ and orthonormal matrix $\hat{\bm{B}}$ for these methods. To overcome this difficulty, we consider the central subspace $\mathcal{S}_{Y|X}$ instead of the central mean subspace $\mathcal{S}_{E(Y|X)}$. The central subspace $\mathcal{S}_{Y|X}$ defined as the intersection of all subspaces ${\rm span}(\bm{A})$ such that $Y \hDash X|\bm{A}^{\top}X$. It is easy to see that $\mathcal{S}_{E(Y|X)} \subset \mathcal{S}_{Y|X}$. Thus we further assume that $\mathcal{S}_{E(Y|X)} = \mathcal{S}_{Y|X}$. This can be achieved when the error terms $\varepsilon$ under the null and alternatives have dimension reduction structures: $\varepsilon=\sigma_1(\bm{\beta_0}^{\top}X) \tilde{\varepsilon}$ and $\varepsilon=\sigma_2(\bm{B}^{\top}X) \tilde{\varepsilon}$ with $\tilde{\varepsilon} \hDash X$ respectively. More details about this issue can be found in Tan and Zhu (2019a).

%Recall that the central mean subspace $\mathcal{S}_{E(Y|X)}$ is defined as the intersection of all subspaces ${\rm span}(\bm{A})$ spanned by the columns of $\bm{A}$ such that $Y \hDash E(Y|X)|\bm{A}^{\top}X$.
%where $\hDash$ denotes statistical independence.
%The dimension of $\mathcal{S}_{E(Y|X)}$ is called the structural dimension, denoted by $d_{E(Y|X)}$.
%Under mild conditions, the subspace $\mathcal{S}_{E(Y|X)}$ always exists (see Cook and Li, (2002)).
%If $\mathcal{S}_{E(Y|X)}={\rm span}(\bm B)$, then  $E(Y|X)=E(Y|\bm{B}^{\top}X)$.
%Under the null hypothesis (\ref{1.1}),  $d_{E(Y|X)}=d$ and $\mathcal{S}_{E(Y|X)}=\rm{span}(\bm{\beta_0})$. Under the alternative (\ref{1.2}),  $d_{E(Y|X)}=q$ and $\mathcal{S}_{E(Y|X)}={\rm span}(\bm{B})$.
%For simplicity, we assume throughout this paper that  $\mathcal{S}_{E(Y|X)}=\mathcal{S}_{Y|X}$. Here $\mathcal{S}_{Y|X}$ is the central subspace of $Y$ with respect to $X$ (see, Cook (1998)).

When $p$ is fixed, there exist many estimation approaches available in the literature to identify the central subspace $\mathcal{S}_{Y|X}$, such as sliced inverse regression (SIR,Li (1991)), sliced average variance estimation (SAVE, Cook and Weisberg (1991)), minimum average variance estimation (MAVE, Xia et.al. (2002)), directional regression (DR, Li and Wang, (2007)), discretization-expectation estimation (DEE, Zhu, et al. (2010a)) etc. Zhu, Miao, and Peng (2006) first discussed the asymptotic properties of SIR in the divergent dimension setting. In this paper we adapt cumulative slicing estimation (CSE, Zhu, Zhu, and Feng (2010b)) to identify the central subspace $\mathcal{S}_{Y|X}$, as it allows the dimension $p$ to diverge to infinity at a rate $p^2/n \to 0$ and is very easy to implement in practice. The basic procedure of CSE is given below.

For simplicity, we assume $E(X)=0$ and $Var(X)=\bm{I_p}$. Under the linearity condition (see Li (1991)), it can be shown that $E[Xh(Y)] \in \mathcal{S}_{Y|X}$ for any function $h(\cdot)$. Therefore, we obtain infinity amount of vectors in $\mathcal{S}_{Y|X}$. Zhu et al. (2010b) suggested the determining class of indicator functions $\mathcal{H}=\{h_t(Y)=I(Y \leq t): t\in \mathbb{R}\}$ to substitute $h(\cdot)$.
%If $h_t(Y)=I(Y \leq t)$, it follows that
%$$ Y \hDash X|\bm{B}^{\top}X  \Longleftrightarrow h_t(Y) \hDash X|\bm{B}^{\top}X,  \quad \forall \ t \in \mathbb{R}. $$
Then define a target matrix as follows
\begin{equation}\label{3.6}
  \bm{M}=\int E[X h_t(Y)] E[X^{T} h_t(Y)] dF_{Y}(t),
\end{equation}
where $F_{Y}(t)$ is the cumulative distribution function of $Y$. If the rank of $\bm{M}$ is $q$, Zhu et al. (2010b) showed that ${\rm span}(\bm M)=\mathcal{S}_{Y|X}$. Let $Z_i$ be the standardized $X_i$ and $\hat{\alpha}_t=\frac{1}{n} \sum_{i=1}^{n} Z_i I(Y_i \leq t)$.
Then the sample version of $\bm M$ is given by
\begin{equation}\label{3.7}
 \bm{\hat{M}}=\frac{1}{n} \sum_{j=1}^{n} \hat{\alpha}_{Y_j} \hat{\alpha}_{Y_j}^{T}.
\end{equation}
If the structural dimension $q$ is known, the estimator $\bm{\hat{B}}(q)$ of $\bm B$ consists of the eigenvectors corresponding to the largest $q$ eigenvalues of $\bm{\hat{M}}$.

Yet we need to estimate the structural dimension $q$. Inspired by the idea of Xia, Xu  and Zhu (2015), we suggest a minimum ridge-type eigenvalue ratio estimator (MRER) to identify $q$. Let $\{\hat{\lambda}_j, 1\leq j \leq p \}$ and $\{\lambda_j, 1\leq j \leq p\}$ be the eigenvalues of the matrix $\bm {\hat{M}}$ and $\bm M$ respectively. Suppose $\hat{\lambda}_{j+1} \leq \hat{\lambda}_{i}$ and $\lambda_{j+1} \leq \lambda_{j}$. Since $rank(\bm M)=q$, we have
$$ \lambda_p = \cdots =\lambda_{q+1}=0 < \lambda_q \leq \cdots \leq \lambda_1 .$$
Hence we estimate the structural dimension $q$ by
\begin{equation}\label{2.7}
\hat{q}=\arg\min_{1\leq j \leq p-1}\left\{j: \frac{\hat{\lambda}^2_{j+1}+c_n}{\hat{\lambda}^2_j+c_n}\right\}.
\end{equation}
Here $c_n$ is a positive constant relying on the sample size $n$. The following result is a slight extension of Tan and Zhu (2019a) that shows the consistency of MRER and model adaptation to the underlying models, if $c_n$ is selected appropriately.

\begin{theorem}\label{theorem3.1}
Suppose the regularity conditions of Theorem 3 in Zhu et al. (2010b) hold and let $\bm{\hat{B}}(q)$ be a matrix whose columns are the eigenvectors associating with the largest $q$ eigenvalues of $\bm{\hat{M}}$. If $0 < c \leq \lambda_q \leq \lambda_1 \leq C < \infty $ and $c_n=\log{n}/n$, then it follows that \\
(1) under $H_0$, we have $\mathbb{P}(\hat{q}=d) \to 1$ and $\| \bm{\hat{B}}(d) - \bm{B} \|=O_p(\sqrt{p/n})$, \\
(2) under $H_1$, we have $\mathbb{P}(\hat{q}=q) \to 1$ and $\| \bm{\hat{B}}(q)- \bm{B} \|=O_p(\sqrt{pq/n})$,\\
where $c$ and $C$ are two positive constants free of $n$ and $\bm B$ under the null $H_0$ satisfied that ${\rm span}(\bm{B})={\rm span}(\bm{\beta_0})$.
\end{theorem}

\section{Asymptotic properties of the test statistic}
\subsection{Limiting null distribution}
%First we consider the asymptotic property of the test statistic $AICM_n$ under the null hypothesis.
Recall that
$$AICM_n= \int_{\mathbb{R}^{d+\hat{q}}} |\frac{1}{\sqrt{n}} \sum_{j=1}^{n} \hat{\eta}_j \exp(\mathrm{i} t^{\top} \bm{\hat{\mathfrak{B}}}^{\top}X_j)|^2 \varphi(t) dt,$$
where $\hat{\eta}_j=Y_j-g(\bm{\hat{\beta}}^{\top}X_j, \hat{\vartheta}_j)$.
To facilitate the derivation of the asymptotical properties, we define the following empirical process
$$\hat{V}_n^1(t)=\frac{1}{\sqrt{n}} \sum_{j=1}^{n} \hat{\eta}_j [\cos( t^{\top}
\bm{\hat{\mathfrak{B}}}^{\top}X_j)+
\sin( t^{\top}\bm{\hat{\mathfrak{B}}}^{\top}X_j)].$$
If $\varphi(t)=\varphi(-t)$, then it follows that
$$AICM_n=\int_{\mathbb{R}^{d+\hat{q}}} |\hat{V}_n^1(t)|^2 \varphi(t) dt. $$

To obtain the large-sample properties of $AICM_n$, we need some regularity conditions. Also recall that $g(\bm{\beta}^{\top}X, \vartheta)=g(\beta_1^{\top}X, \cdots, \beta_d^{\top}X, \vartheta)$. Put
$$g'_i(t_1, \cdots, t_d, \vartheta)=\frac{\partial g}{\partial t_i}(t_1, \cdots, t_d, \vartheta) \quad {\rm for} \ i=1, \cdots, d.    $$

(B1) Assume that $|g'_i(\beta_1^{\top}x, \cdots, \beta_d^{\top}x, \vartheta)| \leq F(x)$ for all $\beta_1, \cdots, \beta_d, \vartheta$ and $i=1, \cdots, d$. Further, assume that
$$ \lambda_{\max}\{E[F(X)XX^{\top}]\} \leq C \quad {\rm and} \quad \lambda_{\max}\{E [\varepsilon^2 F(X)XX^{\top}]\} \leq C,   $$
where $\lambda_{\max}(\cdot)$ is the largest eigenvalue of the matrix and $C$ is a positive constant.

(B2) The weight function $\varphi(t)$ is positive and satisfies that $\varphi(t)=\varphi(-t)$, $\int_{\mathbb{R}^{d+\hat{q}}} \varphi(t)dt < \infty$,  $\int_{\mathbb{R}^{d+\hat{q}}} \|t\|^6 \varphi(t)dt < \infty$, and $ \int_{\mathbb{R}^{2d}} \|M(t)\|^2 \varphi(t)dt = O(1)$. Here $M(t)$ is a shift term which is given by $M(t)=E\{g'(\theta_0,X)[\cos(t^{\top}\bm{\mathfrak{B}_0}^{\top}X)+
\sin(t^{\top}\bm{\mathfrak{B}_0}^{\top}X)]\}$ and $\bm{\mathfrak{B}_0}=(\bm{\beta}_0, \bm{B})$.

(B3) Assume that $E(\varepsilon^8) < \infty$ and $E(X_{jk}^8) \leq C$ where $X_{jk}$ is the $k$-component of $X_j$.

The regularity condition (B1) is similar as that in condition (A2) which is usually used in the diverging dimensional statistical inference, see Fan and Peng (2004) and Zhang and Huang (2008) for instance. Condition (B2) is necessary for the convergence of the remainder in the decomposition of the test statistic in the diverging dimension scenarios. Condition (B3) is standard in model checking literature, see, e.g., Stute (1997) and Escanciano (2006a).
%\begin{eqnarray*}
%\bm{\hat{\beta}_n} &=& \mathop{\rm arg \min}\limits_{\bm{\beta}} \sum_{i=1}^{n}[Y_i-g(\bm{\beta}^{\top}X_i)]^2;\\
%\bm{\tilde{\beta}_0} &=& \mathop{\rm arg \min}\limits_{\bm{\beta}} E[Y-g(\bm{\beta}^{\top}X)]^2.
%\end{eqnarray*}
%If the regression function $m(\cdot) \in \mathcal{G}$, it is easy to see that $\bm{\tilde{\beta}_0}=\bm{\beta_0}$. If $m(\cdot) \notin \mathcal{G}$, $\bm{\tilde{\beta}_0}$ would typically rely on the distribution of $(X, Y)$. To obtain the asymptotical properties of $\bm{\tilde{\beta}_0}$, we first give some notations used in this paper. Suppose that $g(\beta_1^{\top}x, \cdots,  \beta_d^{\top}x)$ is third differentiable with respective to $(\beta_1, \cdots ,\beta_d)$. Putting
%\begin{eqnarray*}
%g'(\bm{\beta}, x)=\frac{\partial g(\bm{\beta}^{\top}x)}{\partial \bm{\beta}} \quad {\rm and} \quad
%g''(\bm{\beta}, x)=\frac{\partial g'(\bm{\beta}, x)}{\partial \bm{\beta}}.
%\quad g'''(\beta, \theta, x)=\frac{\partial g''(\beta, \theta, x)}{\partial (\beta, \theta)}.
%\end{eqnarray*}
%The $pd \times pd$ matrix $g''(\bm{\beta}, x)$ is used in following matrix $\Sigma$ which will play a crucial role in deriving the asymptotical properties of $\bm{\hat{\beta}_n}$:
%\begin{eqnarray*}
%\Sigma = E[g'(\bm{\tilde{\beta}_{0}}, X) g'(\bm{\tilde{\beta}_{0}}, X)^{\top}]-E[ \eta g''(\bm{\tilde{\beta}_0}, X)] =: \Sigma_1- \Sigma_2.
%\end{eqnarray*}

Now we can obtain the asymptotic distribution of $AICM_n$ under the null hypothesis.
By Theorem~\ref{theorem3.1}, $\mathbb{P}(\hat{q}=d) \to 1$ under the null hypothesis. Thus we only work on the event $\{\hat{q}=d \}$. Consequently, $\bm{\hat{\mathfrak{B}}}=[\bm{\hat{\beta}_0}, \bm{\hat{B}}(d)]$ and $AICM_n$ becomes
$$AICM_n=\int_{\mathbb{R}^{2d}} |\hat{V}_n^1(t)|^2 \varphi(t) dt. $$
Under the regularity conditions (A1)-(A4) and (B1)-(B3) and on the event $\{\hat{q}=d\}$, we can show that under the null hypothesis,
\begin{equation}\label{4.1}
\hat{V}_n^1(t) = V_n^1(t) + R_n(t),
\end{equation}
where
$$V_n^1(t)=\frac{1}{\sqrt{n}}\sum_{j=1}^{n}\varepsilon_j
[\cos(t^{\top}\bm{\mathfrak{B}_0}^{\top}X_j)+\sin(t^{\top} \bm{\mathfrak{B}_0}^{\top}X_j)-M(t)^{\top} \Sigma^{-1}g'(\bm{\beta_0}, X_j)],$$
and $R_n(t)$ is a remainder satisfying
$$ \int_{\mathbb{R}^{2d}} |R_n(t)|^2 \varphi(t)dt=o_p(1). $$
The proof of (\ref{4.1}) will be given in the Appendix. Then we can obtain the following result.

%To facilitate the study, we define the following process
%\begin{equation*}
%V_{n}^0(t)= \frac{1}{\sqrt{n}}\sum_{j=1}^n[Y_i- g(\bm{\tilde{\beta}_0}^{\top}X_j)] \exp(\mathrm{i} t^{\top} \bm{\tilde{B}}^{\top}X_j)
%\end{equation*}
%Under the null hypothesis, we have $q=d$ and $\bm{\tilde{\beta}_0}=\bm{\beta_0}$ . %Thus
%$$V^0_{n}(t)= \frac{1}{\sqrt{n}}\sum_{j=1}^n \varepsilon_j \exp(\mathrm{i} t^{\top} \bm{\tilde{B}}^{\top}X_j).$$
%Let $\psi_n(s,t)$ be the covariance function of $V_n^0(t)$. Then
%$$ \psi_n(s, t)=cov(V^0_{n}(s), \overline{V^0_{n}(t)}) = E[\varepsilon^2 \exp(\mathrm{i} (s-t)^{\top} \bm{\tilde{B}}^{\top}X)]. $$
%If $\psi_n(s,t) \to \psi(s,t)$ pointwisely in $(s, t)$, by Lindeberg-Feller Central limit theorem, it is easy to see that finite dimensional distributions of $V^0_{n}(t)$ converge to multivariate normal distributions with zero-mean. An argument in Appendix shows that $V^0_{n}(t)$ is asymptotically tight. Altogether these yield the following result.

%\begin{theorem}\label{Theorem3.1}
%If $\psi_n(s,t) \to \psi(s,t)$ pointwisely in $(s, t)$, under the null hypothesis, we %have
%\begin{eqnarray}\label{3.1}
%  V^0_{n}(t)  \longrightarrow  V_{\infty}(t)  \quad {\rm in\ distribution \ in \ the \ space} \ \mathbb{C}(\mathbb{R}^{2d}),
%\end{eqnarray}
%where $\mathbb{C}(\mathbb{R}^{2d})$ is the space of complex-valued continuous functions on $\mathbb{R}^{2d}$, and $V_{\infty}$ is a zero-mean complex Gaussian process with covariance function $\psi(s,t)$.
%\end{theorem}

\begin{theorem} \label{theorem4.1}
Suppose that the regularity conditions (A1)-(A4) and (B1)-(B3) hold. If $(p\log(n))^3/n \to 0$, under the null $H_0$, we have in distribution
\begin{eqnarray}\label{3.2}
AICM_n  \longrightarrow  \int_{R^{2d}} |V^1_{\infty}(t)|^2 \varphi(t) dt,
\end{eqnarray}
where $V^1_{\infty}(t)$ is a zero-mean Gaussian process with a covariance function $K(s,t)$ which is the point-wise limit of $K_n(s,t)$. Here $K_n(s,t)$ is the covariance function of $V_n^1(t)$, that is,
\begin{eqnarray*}
K_n(s,t)&=& Cov(V_n^1(s), V_n^2(t)) \\
&=& E\{ \varepsilon^2 [\cos((s-t)^{\top} \bm{\mathfrak{B}}_0^{\top}X)+\sin((s+t)^{\top} \bm{\mathfrak{B}}_0^{\top}X)]\} \\
&& -M(t)^{\top} \Sigma^{-1} E\{ \varepsilon^2 [\cos(s^{\top} \bm{\mathfrak{B}}_0^{\top}X)+\sin(s^{\top} \bm{\mathfrak{B}}_0^{\top}X)] g'(\beta_0^{\top}X)\} \\
&& -M(s)^{\top} \Sigma^{-1} E\{ \varepsilon^2 [\cos(t^{\top} \bm{\mathfrak{B}}_0^{\top}X)+\sin(t^{\top} \bm{\mathfrak{B}}_0^{\top}X)] g'(\beta_0^{\top}X)\} \\
&&+M(s)^{\top}\Sigma^{-1} E[\varepsilon^2 g'(\theta_0, X) XX^{\top}] \Sigma^{-1} M(t).
\end{eqnarray*}
\end{theorem}

\subsection{Limiting distribution under the alternative hypotheses}
Now we discuss the asymptotic property of $AICM_n$ under the alternative hypotheses. Consider the following sequence of alternative hypotheses
\begin{eqnarray}\label{4.3}
H_{1n}:\ Y_n=g(\bm{\beta_0}^{\top}X, \vartheta_0)+ r_n G(\bm{B}^{\top}X)+\varepsilon,
\end{eqnarray}
where $E(\varepsilon|X)=0$, $G(\bm{B}^{\top}X)$ is a random variable satisfying $E[G(\bm{B}^{\top}X)]=0$ and $\mathbb{P}(G(\bm{B}^{\top}X)=0)<1$. The convergence rate $r_n$ satisfies $r_n=1/\sqrt{n}$ or $\sqrt{n}r_n \to \infty$.
%The quantity $a=1$ or can be relied on $n$ tending to zero.
To obtain the asymptotical distribution of $AICM_n$ under the alternatives~(\ref{4.3}), we first derive the asymptotic properties of the estimators $\hat{q}$ and $\bm{\hat{\beta}}$, when the dimension $p$ diverges to infinity.

\begin{theorem}\label{theorem4.2}
Suppose that the regularity conditions of Theorem 3 in Zhu et al. (2010b) hold. Let $\bm{\hat{B}}(d)$ be the eigenvectors associating with the largest $d$ eigenvalues of $\hat{M}_n$. If $p r_n \to 0$ and $c_n=r_n^2 \log{r_n^{-2}}$, then under $H_{1n}$, we have $\mathbb{P}(\hat{q}=d) \to 1$ and $\|\bm{\hat{B}}(d)-\bm{B}_{L}\|=O_p(\sqrt{p} r_n)$. Here $\bm{B}_L$ is a $p \times d$ orthonormal matrix satisfying ${\rm span}(\bm{B}_L)={\rm span}(\bm{\beta_0})$.
\end{theorem}

It is worth to mention that under the local alternatives $H_{1n}$ with $r_n = o_p(1/p)$, the estimated structural dimension $\hat q$ is not equal to the true structure dimension, but to $d$ asymptotically. This means $\hat q$ does not an consistent estimator of the true structural dimension in this case. A special case is that if $r_n=1/\sqrt{n}$, it follows that $\mathbb{P}(\hat{q}=d) \to 1$ and $\|\bm{\hat{B}}(d)-\bm{B}_{L}\|=O_p(\sqrt{p/n})$.
Yet we need to derive the asymptotic properties of $(\bm{\hat{\beta}}, \hat{\vartheta})$ with respective to $(\bm{\beta_{0}}, \vartheta_0)$. Recall that $\hat{\theta}=[vec(\bm{\hat{\beta}})^{\top},\hat{\vartheta}^{\top}]^{\top}$ and $\theta_0=[vec(\bm{\beta_0})^{\top},\vartheta_0^{\top}]^{\top}$ where $vec(\bm{\beta})=[\beta_1^{\top}, \cdots, \beta_d^{\top}]^{\top}$.
%and the asymptotical decomposition of $vec(\bm{\hat{\beta}_n})-vec(\bm{\beta_{0}})$ under $H_{1n}$. %Here $\hat{\gamma}_n=(\hat{\beta}_n^{\top}, \hat{\theta}_n^{\top})^{\top}$ and $\gamma_{0}=(\beta_{0}^{\top}, \theta_{0}^{\top})^{\top}$ as mentioned before.

\begin{theorem}\label{theorem4.3}
Suppose the regularity conditions (A1)-(A3) and~(\ref{4.3}) hold. If $pr_n \to 0$, then $\hat{\theta}$ is a norm consistent estimator for $\theta_0$ with $\|\hat{\theta}-\theta_0\|=O_p(\sqrt{p} r_n)$. Moreover, if the regularity conditions (A1)-(A4) hold and $n (p\log{n})^3 r_n^4 \to 0$, then under the alternatives~(\ref{4.3}), we have
\begin{equation}\label{4.4}
\sqrt{n} \alpha^{\top}(\hat{\theta}-\theta_0)= \frac{1}{\sqrt{n}} \sum_{j=1}^{n} \varepsilon_j \alpha^{\top}\Sigma^{-1}g'(\theta_0, X_j)
+\sqrt{n} r_n \alpha^{\top}\Sigma^{-1} E[g'(\theta_0, X_j)G(\bm{B}^{\top}X_j)] + o_p(1),
\end{equation}
where the term $o_p(1)$ is uniformly in $\alpha \in \mathcal{S}^{p-1}$.
\end{theorem}

Under the alternatives $H_{1n}$ with $r_n=1/\sqrt{n}$, if $p^2/n \to 0$, then we have $\hat{\theta}$ is a norm consistent estimator for $\theta_0$ with $\|\hat{\theta}-\theta_0\|=O_p(\sqrt{p/n})$. Moreover, if $(p\log{n})^3 /n \to 0$, it follows that
$$\sqrt{n} \alpha^{\top}(\hat{\theta}-\theta_0)= \frac{1}{\sqrt{n}} \sum_{j=1}^{n} \varepsilon_j \alpha^{\top}\Sigma^{-1}g'(\theta_0, X_j)
+\alpha^{\top}\Sigma^{-1} E[g'(\theta_0, X)G(\bm{B}^{\top}X)] + o_p(1). $$

\begin{theorem}\label{theorem4.4}
Suppose the regularity conditions (A1)-(A4) and (B1)-(B3) hold. \\
(1) If $(p\log{n})^3/n \to 0$, under the global alternative $H_1$, we have in probability
\begin{eqnarray*}
\frac{1}{n} \int_{\mathbb{R}^{d+\hat{q}}} |\hat{V}_{n}(t)|^2 \varphi(t) dt \longrightarrow  C_1 ,
\end{eqnarray*}
where $C_1$ is a positive constant. This means $AICM_{n}$ diverges to infinity at the rate of $n $;  \\
(2) If $n (p\log{n})^3 r_n^4 \to 0$, under the local alternatives $H_{1n}$ with $\sqrt{n} r_n \to \infty$, we have in probability
\begin{eqnarray*}
\frac{1}{n r_n^2} \int_{\mathbb{R}^{d+\hat{q}}} |\hat{V}_{n}(t)|^2 \varphi(t) dt \longrightarrow  C_2,
\end{eqnarray*}
where $C_2$ is a positive constant. This means $AICM_{n}$ diverges to infinity at the rate of $n r_n^2 $; \\
(3) If $(p\log{n})^3/n \to 0$, under the local alternatives $H_{1n}$ with $r_n=1/\sqrt{n}$, we have in distribution
\begin{eqnarray*}
\int_{\mathbb{R}^{d+\hat{q}}} |\hat{V}_{n}(t)|^2 \varphi(t) dt \longrightarrow \int_{\mathbb{R}^{2d}} | V^1_{\infty}(t)+L_1(t)-L_2(t)|^2 \varphi(t) dt,
\end{eqnarray*}
where $V^1_{\infty}$ is a zero-mean Gaussian process given by (\ref{3.2}) and $L_1(t)$ and $L_2(t)$ are the uniformly limits of $L_{n1}(t)$ and $L_{n2}(t)$, respectively. Two functions $L_{n1}(t)$ and $L_{n2}(t)$ are given by
\begin{eqnarray*}
L_{n1}(t)&=& M^{\top}(t) \Sigma^{-1} E[g'(\theta_0, X)G(\bm{B}^{\top}X)] \\
L_{n2}(t)&=& E\{G(\bm{B}^{\top}X)[\cos(t^{\top} \bm{\mathfrak{B}_0}^{\top}X) +\sin(t^{\top} \bm{\mathfrak{B}_0}^{\top}X)]\}.
\end{eqnarray*}
This means $AICM_{n}$ is still sensitive to the local alternatives distinct from the null at the rate of $1/\sqrt n$.
\end{theorem}

\section{Bootstrap approximation}
Note that limiting null distribution of our test statistic $AICM_n$ depends on the parameters $\bm{\beta_0}$ and the matrix $\bm{B}$, and thus is not tractable for the critical value determination.  A typical method used in the literature is the wild bootstrap proposed by Wu (1986). Stute, Gonz$\acute{a}$lez Manteiga and Presedo Quindimil (1998) proved that the wild bootstrap yields a valid approximation of residual marked empirical processes based on the indictor functions, when the dimension of covariates is fixed. We now show that the wild bootstrap still works in the diverging dimension setting. More specially, set
$$X_j^*=X_j \quad {\and} \quad Y_j^*=g(\hat{\theta}, X_j)+\hat{\eta}_j V_j^* , j=1, \cdots, n,$$
where $\hat{\eta}_j$ is the residual and $\{V_j^*\}_{j=1}^n$ are i.i.d. bounded random variables with zero mean and unit variance, independent of the original sample $\{(X_i, Y_i) \}_{j=1}^n$.
%Here $E^*$ is the expectation under the probability measure $P^*$ induced by the wild bootstrap resampling conditional on the original sample $\{(X_i, Y_i): i=1, \cdots, n\}$.
An often used example of $\{V_j^*\}_{j=1}^n$ is the i.i.d. Bernoulli variates with
$$ \mathbb{P}(V_j^*=\frac{1-\sqrt{5}}{2})=\frac{1+2\sqrt{5}}{2\sqrt{5}}, \quad \mathbb{P}(V_j^*=\frac{1+\sqrt{5}}{2})=1-\frac{1+2\sqrt{5}}{2\sqrt{5}}. $$
For other examples of $\{V_j^*\}_{j=1}^n$, one can refer to Mammen (1993).

Let $\hat{\theta}^*$ be the bootstrap estimator obtained by the ordinary least square based on the bootstrap sample $\{(X_j, Y_j^*)\}_{j=1}^n$. Then we approximate the limiting null distribution of $AICM_n$ by that of
$$AICM_n^*= \int_{\mathbb{R}^{d+\hat{q}}} |\hat{V}_n^*(t)|^2 \varphi(t)dt,  $$
where
$$ \hat{V}_n^*(t)=\frac{1}{\sqrt{n}} \sum_{j=1}^{n} [Y_j^*-g(\hat{\theta}^*, X_j)]\exp(\mathrm{i}t^{\top}\hat{\bm{\mathfrak{B}}}^{\top}X_j). $$
To determine the critical value in practice, repeat the above process a large number times, say $B$ times. For a given a nominal level $\tau$, the critical value is determined by the upper $\tau$ quantile of the bootstrap distribution $\{\hat{V}_{nk}^*: k=1, \cdots, B \} $.

In the next theorem we establish the validity of the wild bootstrap in the diverging dimension setting.
\begin{theorem}\label{theorem5.1}
Suppose that the bootstrap sample is generated from the wild bootstrap and the regularity conditions (A1)-(A4) holds.\\
(i) If $(p\log{n})^3/n \to 0$, under the null hypothesis $H_0$ or under the alternative hypothesis with $r_n=1/\sqrt{n}$,  we have with probability 1,
$$ AICM_n^* \longrightarrow \int_{\mathbb{R}^{2d}} |V_{\infty}^{1*}|^2 \varphi(t)dt \quad {\rm in \ distribution}, $$
where $V_{\infty}^{1*}$ have the same distribution as the Gaussian process  $V_{\infty}^1$ given in Theorem~{\ref{theorem4.1}}. \\
(ii) If $n (p\log{n})^3 r_n^4 \to 0$, under the local alternatives $H_{1n}$ with $\sqrt{n} r_n \to \infty$, the result in (i) continues to hold. \\
(iii) If $(p\log{n})^3/n \to 0$, under the global alternative $H_1$, the distribution of $AICM_n^*$ converge to a finite limiting distribution which may different from the limiting null distribution.
\end{theorem}

\section{Numerical studies}

\subsection{Simulations}
In this subsection we conduct some numerical studies to show the performance of the proposed test in finite samples when the dimension of covariates is relatively large. From the theoretical view in this paper, we set $p=[3n^{1/3}]-5$ with $n=100, 200, 400 \ {\rm and}\ 600$. We also compare our test with some existing competitors proposed in the literature, although most of them dealt with fixed dimension.
%We have shown before that the new tests are not distribution-free, that is their asymptotic null distributions are case dependent. Thus we can not obtain the exact critical values for general cases. Here we adopt the wild bootstrap to determine the critical values. For more details, one can refer to Stute et al. (1998a) or Escanciano (2006b).

%Note that our test can be considered as an extension of Bierens' (1982) original test $ICM_n$.
Based on the standard normal density function, the Bierens' (1982) test statistic becomes
$$ ICM_n = \int_{\mathbb{R}^p}|\frac{1}{\sqrt{n}} \sum_{j=1}^{n} \hat{\eta}_j \exp(\mathrm{i} t^{\top} X_j))|^2 \phi(t) dt = \frac{1}{n} \sum_{j,k=1}^{n} \hat{\eta}_j \hat{\eta}_k \exp(-\frac{1}{2}\|X_j-X_k\|^2), $$
where $\hat{\eta}_j$ is the residual and $\phi(t)$ is the standard normal density on $\mathbb{R}^p$. Escanciano (2006a) also use this test statistic for comparison.

Zheng (1996) proposed a local smoothing test for parametric regression models as
$$ T_n^{ZH}=\frac{\sum_{i \neq j} K((X_i-X_j)/h) \hat{\eta}_i \hat{\eta}_j} {\{\sum_{i \neq j} 2K^2((X_i-X_j)/h)\hat{\eta}_i^2\hat{\eta}_j^2\}^{1/2}}. $$
Here we use $K(u)=(15/16)(1-u^2)^2I(|u| \leq 1)$ as the kernel function and the bandwidth $h=1.5n^{-1/(4+p)}$.

Escanciano (2006a) developed a global smoothing test for parametric regression models based on a projected residual marked empirical process. The test statistic is defined as follows,
$$ PCvM_n=\frac{1}{n^2}\sum_{i,j,r=1}^{n} \hat{\eta}_i \hat{\eta}_j \int_{\mathbb{S}^{p-1}}I(\alpha^{\top}X_i \leq \alpha^{\top}X_r)I(\alpha^{\top}X_j \leq \alpha^{\top}X_r) d\alpha, $$
where $\mathbb{S}^{p-1}=\{\alpha \in \mathbb{R}^p: \|\alpha\|=1\}$ and $d\alpha$ denotes the uniform density on the unit sphere $\mathbb{S}^{p-1}$. The critical value of Escanciano's (2006a) test is determined by the wild bootstrap.

Stute and Zhu (2002) proposed a dimension reduction-based test for generalized linear models based on residual marked empirical processes. A martingale transformation leads the test to be asymptotically distribution-free. Their test statistic is given by
$$ T_n^{SZ}=\frac{1}{\hat{\psi}_n(x_0)} \int_{-\infty}^{x_0} |\hat{T}_nR_n^1|^2 \hat{\sigma}_n^2 dF_n, $$
where
\begin{eqnarray*}
R_n^1(u)&=& \frac{1}{\sqrt{n}} \sum_{i=1}^{n} \hat{\eta}_i I(\hat{\beta}^{\top} X_i \leq u);\\
\hat{T}_n R_n^1(u) &=& R_n^1(u)- \int_{-\infty}^{u} \hat{a}_n(z)^{\top} \hat{A}_n^{-1}(z) \left(\int_{z}^{\infty}\hat{a}_n(v) R_n^1(dv)\right) \hat{\sigma}_n^2(z) F_n(dz).
\end{eqnarray*}
More details of the definitions of $\hat{T}_n R_n^1$ can be found in their paper. Under the right specification of the generalized linear model and some regularity conditions,
$$ T_n^{SZ}\longrightarrow  \int_{0}^{1} B^2(u)du  \quad {\rm in \ distribution}, $$
where $B(\cdot)$ is the standard Brownian motion.

Guo, Wang and Zhu (2016) introduced a model-adaptive local smoothing test for parametric single index models that largely alleviate the dimensionality problem, although they also considered in the fixed dimension cases. Their test statistic is given by
\begin{equation*}
T_n^{GWZ}=\frac{h^{1/2}\sum_{i \neq j} \hat{\eta}_i \hat{\eta}_j \frac{1}{h^{\hat{q}}} K(\bm{\hat{B}}^{\top}(X_i-X_j)/h) } {\{2\sum_{i \neq j} \hat{\eta}_i^2 \hat{\eta}_j^2 \frac{1}{h^{\hat{q}}} K^2(\bm{\hat{B}}^{\top}(X_i-X_j)/h)\}^{1/2}},
\end{equation*}
where the kernel function $K(u)=(15/16)(1-u^2)^2I(|u| \leq 1)$ and the bandwidth $h=1.5n^{1/(4+\hat{q})}$ as suggested in Guo, Wang and Zhu (2016) and $\bm{\hat{B}}$ is a sufficient dimension estimator of $\bm{B}$ with an estimated structural dimension $\hat{q}$ of $q$.

Recently, Tan and Zhu (2019a) proposed a projected adaptive-to-model test for parametric single index models where they also allow the dimension $p$ to diverge with the sample size $n$. Their test statistic is given by
$$ ACM_n^2 = \frac{1}{\hat{\psi}_n(u_0)^2} \int_{-\infty}^{u_0} \sup_{\hat{\alpha} \in \mathcal{S}_{\hat{q}}^{+}}  |\hat{T}_n V_{n}(\hat{\alpha}, u)|^2 \hat{\sigma}_n^2(u) F_n(du), $$
where
\begin{eqnarray*}
\hat{T}_n V_n(\hat{\alpha}, u) &=& \frac{1}{n^{1/2}}\sum_{i=1}^n \hat{\eta}_i I(\hat{\alpha}^{\top}\bm{\hat{B}}^{\top} X_i \leq u)-  \frac{1}{n^{3/2}} \sum_{i,j=1}^n I(\hat{\alpha}^{\top}\bm{\hat{B}}^{\top} X_i \leq u) \hat{a}_n(\hat{\alpha}^{\top}\bm{\hat{B}}^{\top} X_i)^{\top}   \\
&& \times \hat{A}_n^{-1}(\hat{\alpha}^{\top}\bm{\hat{B}}^{\top} X_i) \hat{a}_n(\hat{\alpha}^{\top}\bm{\hat{B}}^{\top} X_j)
I(\hat{\alpha}^{\top}\bm{\hat{B}}^{\top} X_j \geq \hat{\alpha}^{\top} \bm{\hat{B}}^{\top} X_i) \hat{\eta}_j\hat{\sigma}_n^2(\hat{\alpha}^{\top} \bm{\hat{B}}^{\top} X_i).
\end{eqnarray*}
For the quantities $\hat{\psi}(u_0), \hat{\sigma}_n^2, \hat{a}_n$ and $\hat{A}_n$, one can refer their paper for details. Under the null hypothesis and some regularity conditions, we have
$$ ACM_n^2  \longrightarrow  \int_{0}^{1} B^2(u)du  \quad {\rm in \ distribution}, $$
Thus this test is asymptotically distribution-free and its critical values can be tabulated.

%First we compare our tests $AICM_n$ with Bierens's (1982) original test $ICM_n$, . As we mentioned before, the integral in Bierens's (1982) test becomes extremely difficulty to handle when the dimension $p$ of predictors is large. Thus we suggest to use the standard normal probability measure in Bierens's (1982) test statistic and then its test statistic becomes
%$$ ICM_n = \int_{\mathbb{R}^p}|\frac{1}{\sqrt{n}} \sum_{j=1}^{n} \hat{e}_j \exp(\mathrm{i} t^{\top} X_j))|^2 \phi(t) dt
%= \frac{1}{n} \sum_{j,k=1}^{n} \hat{e}_j \hat{e}_k \exp(-\frac{1}{2}\|X_j-X_k\|^2), $$
%where $\phi(t)$ is the standard normal density on $\mathbb{R}^p$.

In the simulations that follows, $a=0$ corresponds to the null while $a\neq 0$ to the alternatives. The significance level is $\alpha= 0.05$. The simulation results are based on the average of $1000$ replications and the bootstrap approximation of  $B=500$ replications.

$Study$ 1. Consider the following regression models
\begin{eqnarray*}
H_{11}:  Y &=& \beta_0^{\top}X+ a \exp(-(\beta_0^{\top}X)^2) +\varepsilon;  \\
H_{12}:  Y &=& \beta_0^{\top}X+ a \cos(0.6\pi \beta_0^{\top}X) +\varepsilon;\\
H_{13}:  Y &=& \beta_1^{\top}X+ a (\beta_2^{\top}X)^2 +\varepsilon; \\
H_{14}:  Y &=& \beta_1^{\top}X+ a \exp(\beta_2^{\top}X)+\varepsilon;
%H_{5}:  Y &=& \exp(\beta_1^{\top}X)+ a(\beta_2^{\top}X)+\varepsilon.
\end{eqnarray*}
where $\beta_0=(1, \cdots, 1)^{\top}/\sqrt{p}$, $\beta_1=(\underbrace{1,\dots,1}_{p_1},0,\dots,0)^{\top}/\sqrt{p_1}$ and $\beta_2=(0,\dots,0,\underbrace{1,\dots,1}_{p_1})/\sqrt{p_1}$ with $p_1=[p/2]$.  The covariate $X$ is $N(0, I_p)$ independent of the standard Gaussian error term $\varepsilon$.
%and  $N(0,\Sigma)$ with $\Sigma=(1/2^{|i-j|})_{p\times p}$ in order to check the influence of correlation between the covariates.
Note that $H_{12}$ is a high-frequency model and the others are low-frequency models. The structure dimension $q=1$ under both the null and alternative hypotheses in the models $H_{11}$ and $H_{12}$, while the structure dimension $q=2$ under the alternative hypothesis in models $H_{13}$ and $H_{14}$.

The empirical sizes and powers are presented in Tables 1-4. First we can see that the Bierens' (1982) original test $ICM_n$ performs the worst among these tests. It almost has no empirical powers at all and its empirical sizes are far away from the nominal sizes in the  large dimension cases. In our experience, when the dimension $p$ of  $X$ is smaller than 5, $ICM_n$ can maintain the significance level very well and the empirical power grows very quickly. But for higher dimension,  the Bierens's (1982) test does not work well.
For the other tests, we observe that $AICM_n$, $ACM_n$, $T_n^{SZ}$, and $T_n^{GWZ}$ can control the empirical sizes very well. The empirical sizes of $PCvM_n$ are also close to the significance level, but slightly unstable in some cases. While $T_n^{SZ}$ can not maintain the significance level in most cases and is generally conservative with smaller empirical sizes. For the empirical power, we can see that $AICM_n$, $ACM_n$, $T_n^{SZ}$, and $PCvM_n$ all perform very well for low frequency models $H_{11}, H_{13}$ and $H_{14}$, whereas $T_n^{GWZ}$ behaves slightly worse for the three low frequency models. For the high frequency model $H_{12}$, the test $T_n^{GWZ}$ beats all other competitors except for the new test $AICM_n$. This is somewhat surprised as local smoothing tests such as $T_n^{GWZ}$ usually performs better for high frequency models and global smoothing tests work better for low frequency models. While  the new test which can be viewed as a global smoothing test seems also to work well for high frequency models. In contrast, Zheng's (1996) test $T_n^{ZH}$ which is a typical local smoothing test has no empirical powers in all cases. This validates the well-known results that the traditional local smoothing tests suffer severely from the curse of dimensionality.
$$ \rm{Tables \ 1-4 \ about \ here} $$
%Next we only consider the other three tests that all allow the dimension $p$ to diverge to infinity with the sample size $n$. We can observe that both $\hat{V}_{U,n} $ and $\hat{V}_{G,n} $ can control the nominal size at $\alpha=0.1$ in linear null models. When $\alpha=0.05, 0.01$, the empirical sizes $\hat{V}_{U,n} $ and $\hat{V}_{G,n} $ seem slightly larger than their corresponding nominal sizes. In contract, $ACM_n^2$ can maintain the significance levels in all linear null models. In nonlinear null model $H_5$, empirical sizes of all tests are a little larger than the nominal sizes. For the empirical powers, we find that $\hat{V}_{G,n} $ usually has higher power than the other two in models $H_1$, $H_2$ and $H_5$. Moreover, empirical powers of $\hat{V}_{U,n}$ in models $H_2$ and $H_5$ seems a little low, even when $n=800$. Note that $H_2$ is a high frequency model. While, in models $H_3$ and $H_4$, $\hat{V}_{U,n}$ performs slightly better than the other two. %It seems no one can always be the winner.

The hypothetical models in study 1 are all single-index models. Next we consider multiple-index models in the second simulation study. As the tests $ACM_n, T_n^{SZ}$ and $T_n^{GWZ}$ only dealt with parametric single index models, we only compare our new test with $PCvM_n, ICM_n$ and $T_n^{ZH}$.

$Study$ 2. Generate data from the following models:
\begin{eqnarray*}
H_{21}: Y &=& \beta_1^{\top}X+\exp(\beta_2^{\top}X)+a (\beta_2^{\top}X)^2 +\varepsilon;  \\
H_{22}: Y &=& \beta_1^{\top}X+\exp(\beta_2^{\top}X)+a \cos(0.6 \pi \beta_2^{\top}X) +\varepsilon;  \\
H_{23}: Y &=& \beta_1^{\top}X+\exp(\beta_2^{\top}X)+a (\beta_1^{\top}X) (\beta_2^{\top}X) + \varepsilon;
\end{eqnarray*}
where $\beta_1, \beta_2, \varepsilon$ and $X$ are the same as in Study 1.

The simulation results are presented in Tables 5-7. We can observe that our test $AICM_n$ and $PCvM_n$ perform much better than the other two. While Bierens' (1982) test $ICM_n$ again does not work at all and Zheng's (1996) test $T_n^{ZH}$ can not maintain the nominal level and has no empirical powers in both cases. For the empirical size, both our test $AICM_n$ and $PCvM_n$ are slightly conservative with larger empirical sizes. This may be due to the inaccurate estimation of the related parameters when $p$ is large. The empirical powers of $AICM_n$ and $PCvM_n$ both grow fast under both the low frequency model $H_{21}$ and the high frequency model $H_{22}$. In model $H_{23}$, our test has much better power performance than  $PCvM_n$.
$$ \rm{Tables \ 5-7 \ about \ here} $$

In summary, the simulation results show that the proposed test  performs well and can detect different alternative hypotheses for large $p$ paradigms. In low frequency alternatives, the new test has the best power performance in high-dimensional cases among the global smoothing tests proposed by Bierens' (1982), Stute and Zhu (2002), Escanciano (2006a), and Tan and Zhu (2019a). While surprisingly, the new test, which can be viewed as a global smoothing test, also performs best or  is comparable to the best one among the local smoothing tests proposed by Zheng (1996) and  Guo, Wang and Zhu (2016).

\subsection{A real data example}
In this subsection we apply the test to the Boston house-price data set that is first analysed by Harrison and Rubinfeld (1978). This data set can be obtained through the website \url{http://lib.stat.cmu.edu/datasets/boston}. It contains 506 cases on one response variable: Median value of owner-occupied homes in 1000's (MEDV) $Y$ and 13 predictors: per capita crime rate by town (CRIM) $X_1$,
proportion of residential land zoned for lots over 25,000 sq.ft. (ZN) $X_2$,
proportion of non-retail business acres per town (INDUS) $X_3$,
Charles River dummy variable (CHAS, 1 if tract bounds river; 0 otherwise) $X_4$,
nitric oxides concentration (NOX parts per 10 million) $X_5$,
average number of rooms per dwelling (RM) $X_6$,
proportion of owner-occupied units built prior to 1940 (AGE) $X_7$,
weighted distances to five Boston employment centres (DIS) $X_8$,
index of accessibility to radial highways (RAD) $X_9$,
full-value property-tax rate per \$ 10,000 (TAX) $X_{10}$,
pupil-teacher ratio by town (PTRATIO) $X_{11}$,
$1000(B_k - 0.63)^2$ where $B_k$ is the proportion of blacks by town (B) $X_{12}$,
\% lower status of the population (LSTAT) $X_{13}$.
For easy explanation, we standardize all variables separately. Since the dimension of the predictor $13 \approx 506^{0.4119}$, it seems reasonable to apply our method.
To establish the relationship between the response $Y$ and the covariates $X=(X_1, \cdots, X_{13})^{\top}$, we first apply the dimension reduction techniques to the data set. When the cumulative slicing estimation (CSE) is used, we find that the estimated structural dimension of this data set is $\hat{q}=1$, which indicates that $Y$ may be conditional independent of $X$ given the projected covariates $\hat{\beta}_1^{\top}X$ where
\begin{eqnarray*}
\hat{\beta}_1 &=& (-0.0804, 0.0611, 0.0211, 0.0676, -0.3120, 0.1513, -0.1820, -0.3701, 0.3999, -0.3404, \\
              &&   -0.2876, 0.1233, -0.5663)^{\top}.
\end{eqnarray*}
The scatter plot of $Y$ against $\hat{\beta}_1^{\top}X$ is presented in Figure 2.
$$ \rm Figure \ 2 \ about \ here $$
From this Figure, a linear regression model seem to be plausible to fit the data set. We then apply our test to see whether these exists a model misspecification. The value of the test statistic $AICM_n$ is about $0.9639$ and the $p$-value is about $0.01$. Thus a linear regression model is not adequate to predict the response. Figure 3 presents the scatter plot of residuals from the linear regression model for $(Y,X)$ against $\hat{\beta}_1^{\top}X$. It also suggests that a linear relationship between $Y$ and $X$ may be not reasonable.
$$ \rm Figure \ 3 \ about \ here $$

To find the relationship between $Y$ and $X$, we give a more thorough search of the projected covariates. Consider the second projected covariates $\hat{\beta}_2^{\top}X$ and the scatter plot of $Y$ against $(\hat{\beta}_1^{\top}X, \hat{\beta}_2^{\top}X)$ is presented in Figure 4.
$$ \rm Figure \ 4 \ about \ here $$
From this figure, we can see that the second projected direction $\hat{\beta}_2^{\top}X$ is not necessary as the plot along this direction is nearly identical. This means the projection of the data onto the space $\hat{\beta}_1^{\top}X$ already contain almost all the regression information of the model structure. Further more, from Figure 3, it seems that there exists an exponential relationship between $Y$ and $\hat{\beta}_1^{\top}X$. Thus we consider $\log(Y)$ as the response variable that is also considered in Harrison and Rubinfeld (1978). Thus we use the following model to fit this data set:
\begin{equation}\label{6.1}
  \log(Y)=\beta^{\top}X+\varepsilon.
\end{equation}
When applying our test for the above model~(\ref{6.1}), we obtain that the value of the test statistic $AICM_n$ is about $0.1536$ and the $p$-value is about 0.368.
A scatter plot of residuals from model~(\ref{6.1}) against the fitted values also does not reveal any trend given in Figure 5.  Thus this model is plausible.
$$ \rm Figure \ 5 \ about \ here $$

\section{Discussions}
In this article we discuss three issues in statistical inference. First we extend Bierens' (1982) classic ICM test, via a model adaptation strategy,  to the diverging dimension scenarios for parametric multiple index models. In fixed dimension scenarios, it is generally believed that the ICM test is less sensitive to the dimension $p$ of covariates due to its projection-based nature. However, we note that even when $p$ is moderate, Bierens' original ICM test does not work at all. Thus, to attack the dimensionality problem, we construct an adaptive-to-model version of ICM test in terms of sufficient dimension reduction and projection techniques. The numerical studies suggest that the proposed test largely alleviates the adverse impact of the dimensionality in the sense that it can well maintain the significance level and enhance power performance. We also show that a large class of weight functions $\varphi(t)$ give the corresponding test statistics with the nice property of computational simplicity. Thus the test can be easily implemented in practice. Second, we obtain the uniformly asymptotically linear presentation of the least squares estimation of the parameters in the nonlinear regression model at a divergent rate $p=o(n^{1/3}/\log{n})$. This result can be useful for further studies in inference.  Third, the study on the consistency of the wild bootstrap approximation  in this diverging setting can also be used in other hypothesis testing problems.
%This paper also leaves some problems. One important issue is about how to relax the condition on the divergent rate of the dimension $p$ to the sample size $n$. In this paper, we can at most have a divergent rate $p=o(n^{1/3}/\log{n})$.
%Since this is not the main purpose of this paper, we leave it for further research.
Note that the test we proposed here is for the hypothetical models with dimension reduction structures. It is of great interest to investigate goodness of fit testing for hypothetical models without such a structure in the divergent dimension setting. The relevant studies are ongoing.

\section{Appendix}
\textbf{Proof of Theorem~\ref{theorem2.1}.} Let $\theta=\tilde{\theta}_0+\zeta$ and $F(\zeta)=\sum_{i=1}^{n}[Y_i-g(\tilde{\theta}_0+\zeta, X_i)]g'(\tilde{\theta}_0+\zeta,
X_i)$. Then it suffices to show that there is a root $\hat{\zeta}$ of $F(\zeta)$ such that $\|\hat{\zeta}\|^2=O_p(m/n)$. According to the results in (6.3.4) of Ortega and Rheinboldt (1970), it in turn needs to prove that $\zeta^{\top}F(\zeta)<0$ for $\|\zeta\|^2=C m/n$ where $C$ is some large enough constant.

Recall $\psi_{\theta}(x,y)=[y-g(\theta, x)]g'(\theta, x)$. Then $F(\zeta)=\sum_{i=1}^{n}\psi_{\tilde{\theta}_0+\zeta}(X_i, Y_i)$. Set
$$ P_n \psi_{\theta}=\frac{1}{n} \sum_{i=1}^n \psi_{\theta}(X_i, Y_i) \quad {\rm and} \quad P\psi_{\theta}=E[\psi_{\theta}(X, Y)].$$
Since $P\psi_{\tilde{\theta}_0}=0$, it follows that
\begin{eqnarray*}
\zeta^{\top}F(\zeta)&=&n \zeta^{\top} P_n\psi_{\tilde{\theta}_0+\zeta} \\
&=&n\zeta^{\top}[P_n(\psi_{\tilde{\theta}_0+\zeta}-\psi_{\tilde{\theta}_0})-
P(\psi_{\tilde{\theta}_0+\zeta}-\psi_{\tilde{\theta}_0})]+ \\
&& n \zeta^{\top} P_n\psi_{\tilde{\theta}_0} + n \zeta^{\top} P\psi_{\tilde{\theta}_0+\zeta}
\end{eqnarray*}
For the first term in $\zeta^{\top}F(\zeta)$, we have
\begin{eqnarray*}
&&E\{\zeta^{\top}[P_n(\psi_{\tilde{\theta}_0+\zeta}-\psi_{\tilde{\theta}_0})-P(\psi_{\tilde{\theta}_0+\zeta}-
\psi_{\tilde{\theta}_0})]\}^2 \\
&=&\frac{1}{n} Var\{\zeta^{\top} [\psi(\tilde{\theta}_0 +\zeta, X, Y)-\psi(\tilde{\theta}_0, X, Y)] \} \\
&\leq& \frac{1}{n} E|\zeta^{\top} [\psi(\tilde{\theta}_0 +\zeta, X, Y)-\psi(\tilde{\theta}_0, X, Y)]|^2\\
&=& \frac{1}{n} E|\zeta^{\top}\frac{\partial \psi}{\partial \theta}(\tilde{\theta}_1, X,Y)\zeta|^2 \leq C \frac{m^2}{n} \|\zeta\|^4,
\end{eqnarray*}
where $\tilde{\theta}_1$ lies between $\tilde{\theta}_0+\zeta$ and $\tilde{\theta}_0$.
Therefore,
$$n \zeta^{\top}[P_n(\psi_{\tilde{\theta}_0+\zeta}-\psi_{\tilde{\theta}_0})-
P(\psi_{\tilde{\theta}_0+\zeta}-\psi_{\tilde{\theta}_0})]=\sqrt{nm^2}\|\zeta\|^2O_p(1).
$$
Note that $P\psi_{\tilde{\theta}_0}=0$. Then it follows that
\begin{eqnarray*}
E\|P_n\psi_{\tilde{\theta}_0}\|^2 =\sum_{j=1}^m E|P_n\psi_{j,\tilde{\theta}_0}|^2 = \frac{1}{n} \sum_{j=1}^m E\psi_j(\tilde{\theta}_0, X, Y)^2 \leq \frac{m}{n} C.
\end{eqnarray*}
Thus we obtain that $n \zeta^{\top} P_n\psi_{\tilde{\theta}_0}=\sqrt{m n}\|\zeta\| O_p(1)$.

Next we consider the term $n \zeta^{\top} P\psi_{\tilde{\theta}_0+\zeta}$. Applying Taylor expansion,
\begin{eqnarray*}
n\zeta^{\top}P\psi_{\tilde{\theta}_0+\zeta}
&=& n\zeta^{\top}\dot{P}\psi_{\tilde{\theta}_0} \zeta + n \zeta^{\top} [\zeta^{\top}(\ddot{P}\psi_{\theta}|_{\theta=\tilde{\theta}_2})\zeta] \\
&=& -n \zeta^{\top}\Sigma\zeta + n \zeta^{\top} [\zeta^{\top} (\ddot{P} \psi_{\theta}|_{\theta=\tilde{\theta}_2}) \zeta],
\end{eqnarray*}
where $\tilde{\theta}_2$ lies between $\tilde{\theta}_0 +\zeta$ and $\tilde{\theta}_0$. By Assumption 2, we obtain
\begin{eqnarray*}
\|\zeta^{\top} (\ddot{P} \psi_{\theta}|_{\theta=\tilde{\theta}_2}) \zeta\|^2 =\sum_{j=1}^m |\zeta^{\top}(\ddot{P} \psi_{j\theta}|_{\theta=\tilde{\theta}_2})
\zeta|^2 \leq \sum_{j=1}^m \max_{1\leq i \leq m} |\lambda_i(\ddot{P} \psi_{j\theta}|_{\theta=\tilde{\theta}_2})|^2 \|\zeta\|^4 \leq C m \|\zeta\|^4 .
\end{eqnarray*}
Then it follows that
$$ n\zeta^{\top}P\psi_{\tilde{\theta}_0+\zeta}=-n\zeta^{\top} \Sigma\zeta + n\sqrt{m} \|\zeta\|^3 O_p(1). $$
Let $\zeta=\sqrt{m/n}U$. Altogether we obtain
\begin{eqnarray*}
\zeta^{\top}F(\zeta)&=&\sqrt{nm^2}\|\zeta\|^2 O_p(1)+ \sqrt{mn} \|\zeta\| O_p(1)-n\zeta^{\top} \Sigma \zeta +\sqrt{m}n\|\zeta\|^3 O_p(1)\\
&=& \frac{m^2}{\sqrt{n}} \|U\|^2 O_p(1)+m \|U\|O_p(1)-mU^{\top}\Sigma U + \frac{m^2}{\sqrt{n}} \|U\|^3 O_p(1) \\
&\leq& m \|U\|O_p(1)-m \lambda_{\min}(\Sigma) \|U\|^2+\frac{m^2}{\sqrt{n}} \|U\|^2 O_p(1)+\frac{m^2}{\sqrt{n}} \|U\|^3 O_p(1)\\
&=& m \|U\| \{ O_p(1)-\lambda_{\min}(\Sigma) \|U\|+ \frac{m}{\sqrt{n}} \|U\| O_p(1)+\frac{m}{\sqrt{n}} \|U\|^2 O_p(1)\}
\end{eqnarray*}
If $m^2/n \to 0$ and $\|U\|=C$ is large enough, for any $\epsilon > 0$, we have
$$\mathbb{P}(\zeta^{\top}F(\zeta)<0) \geq \mathbb{P}\{O_p(1)-\lambda_{\min}(\Sigma)\|U\|+\|U\| o_p(1)<0\} \geq 1-\epsilon .$$
Thus our result follows from (6.3.4) of Ortega and Rheinboldt (1970). Hence we complete the proof. \hfill$\Box$

\textbf{Proof of Theorem~\ref{theorem2.2}.} We use the same notations as those in the proof for Theorem~\ref{theorem2.1}. Let $G_n\psi_{\theta}=\sqrt{n}(P_n\psi_{\theta}-P\psi_{\theta})$ and $H=\{\theta: \theta \in \mathbb{R}^m, \|\theta-\tilde{\theta}_0\| \leq C \sqrt{(m \log n)/n} \}$. First we show that $G_n\alpha^{\top}(\psi_{\hat{\theta}}-\psi_{\tilde{\theta}_0})=o_p(1)$ uniformly in $\alpha \in \mathcal{S}^{q-1}$. Fix $\epsilon >0$. Then we have
\begin{eqnarray*}
&&\mathbb{P}(\sup_{\alpha \in \mathcal{S}^{m-1}}
|G_n\alpha^{\top}(\psi_{\hat{\theta}}-\psi_{\tilde{\theta}_0})|> 8\epsilon) \\
&\leq&\mathbb{P}(\sup_{\alpha \in \mathcal{S}^{m-1}}
|G_n\alpha^{\top}(\psi_{\hat{\theta}}-\psi_{\tilde{\theta}_0})| > 8\epsilon, \hat{\theta} \in H) + \mathbb{P}(\hat{\theta} \notin H ) \\
&\leq&\mathbb{P}(\sup_{\alpha \in \mathcal{S}^{m-1},\theta \in H}
|G_n\alpha^{\top}(\psi_{\theta}-\psi_{\tilde{\theta}_0})| > 8\epsilon) + \mathbb{P}(\|\hat{\theta}-\tilde{\theta}_0\| > C \sqrt{(m \log n)/n} )
\end{eqnarray*}
Since $\|\hat{\theta}-\tilde{\theta}_0\|=O_p(\sqrt{m/n})$, it suffices to show that
$$ \mathbb{P}(\sup_{\alpha \in \mathcal{S}^{m-1}, \theta \in H}
|G_n\alpha^{\top}(\psi_{\theta}-\psi_{\tilde{\theta}_0})| > 8\epsilon) \rightarrow 0. $$
Set $\mathcal{F}=\{\alpha^{\top}[\psi_{\theta}(x, y)-\psi_{\tilde{\theta}_0}(x, y)]: \alpha \in \mathcal{S}^{m-1}, \theta \in H \}$. By Assumption (A4), we obtain that
\begin{eqnarray*}
&&|\alpha^{\top}[\psi_{\theta}(x, y)-\psi_{\tilde{\theta}_0}(x, y)]|\\
&=&|\alpha^{\top}[\psi(\theta, x, y)-\psi(\tilde{\theta}_0, x, y)]| \\
&=&|\alpha^{\top}\frac{\partial \psi}{\partial \theta}(\tilde{\theta}_0, x,y)(\theta-\tilde{\theta}_0)+ \frac{1}{2} (\theta-\tilde{\theta}_0)^{\top}\{\alpha^{\top}\frac{\partial \psi}{\partial \theta^{\top} \partial \theta}(\tilde{\theta}_1, x,y) (\theta-\tilde{\theta}_0)\}|\\
&\leq& \|\theta-\tilde{\theta}_0\| \interleave \frac{\partial \psi}{\partial \theta}(\tilde{\theta}_0, x,y)\interleave_2 + m^{3/2} \|\theta-\tilde{\theta}_0\|^2 F_2(x,y)\\
&\leq&\|\theta-\tilde{\theta}_0\| \{\interleave\frac{\partial \psi}{\partial \theta}(\tilde{\theta}_0, x,y)\interleave_2 + m^2 \sqrt{(\log n)/n} F_2(x,y) \}
\end{eqnarray*}
where $\tilde{\theta}_1$ lies between $\theta$ and $\tilde{\theta}_0$ and $\interleave \cdot \interleave_2$ denotes the operator norm of a matrix in the sense that
$$\interleave A \interleave_2 =\sup_{\|v\|=1}\|Av\|.$$
Let $F(x,y)=\interleave \frac{\partial \psi}{\partial \theta}(\tilde{\theta}_0, x,y)\interleave_2 + m^2 \sqrt{(\log n)/n} F_2(x,y)$. Next we show that $E|F(X,Y)|^2 < Cm$.

Since $(m \log{n})^3/n \to 0$, it follows that
$ E[m^2 \sqrt{(\log n)/n} F_2(x,y)]^2 \leq C m $.
It remains to show that
$E[\interleave \frac{\partial \psi}{\partial \theta}(\tilde{\theta}_0, x,y)\interleave_2]^2 \leq C m$.
According to the arguments of Lemma 3 of Cai et al. (2010), there exists vectors $\{v_j \in \mathbb{S}^{m-1}, 1\leq j \leq 7^m\}$ such that
$$ \interleave A \interleave_2 \leq 4 \sup_{1\leq j \leq 7^m} |v_j^{\top} A v_j| $$
for any symmetric matrix $A$. Set $Z_j=v_j^{\top} \frac{\partial \psi}{\partial \theta}(\tilde{\theta}_0, X,Y)v_j$. Then it follows that
\begin{eqnarray*}
E \interleave \frac{\partial \psi}{\partial \theta}(\tilde{\theta}_0, X,Y) \interleave^2_2
\leq 4 E[\sup_{1\leq j \leq 7^m} |Z_j|^2]
\leq 8 E[\sup_{1\leq j \leq 7^m} |Z_j-EZ_j|^2]+8E[\sup_{1\leq j \leq 7^m} |E Z_j|^2]
\end{eqnarray*}
Note that $EZ_j=-v_j^{\top} \Sigma v_j$. Thus
$$8E[\sup_{1\leq j \leq 7^m} |E Z_j|^2]\leq C. $$
For the term $E[\sup_{1\leq j \leq 7^m} |Z_j-E Z_j|^2]$, we have
\begin{eqnarray*}
E[\sup_{1\leq j \leq 7^m} |Z_j-EZ_j|^2]
&=& \frac{1}{w} E \log \{\exp(\sup_{1\leq j \leq 7^m} w|Z_j-EZ_j|^2 )\} \\
&\leq& \frac{1}{w} \log E\{\sup_{1\leq j \leq 7^m} \exp(w|Z_j-EZ_j|^2)\}\\
&\leq& \frac{1}{w} \log\{\sum_{j=1}^{7^m} E[\exp(w|Z_j-EZ_j|^2)] \}.
\end{eqnarray*}
where $w$ is any positive number. Since $Z_j-EZ_j$ is sub-Gaussian with proxy variance $\sigma^2$, by Lemma 1.3 of Rigollet and H${\rm \ddot{u}}$tter (2017), we obtain that
\begin{eqnarray*}
E[\exp(w|Z_j-EZ_j|^2)]
&=& \int_{0}^{\infty} \mathbb{P}\{\exp(w|Z_j-EZ_j|^2) >t \}dt \\
&=& 1 + \int_{1}^{\infty} \mathbb{P}\{|Z_j-EZ_j| >\sqrt{\frac{\log{t}}{w}} \}dt \\
&\leq& 1+2 \int_{1}^{\infty} \exp(-\frac{\log t}{2w \sigma^2})dt \\
&=& 1+2 \int_{1}^{\infty} t^{-\frac{1}{2w \sigma^2}}dt.
\end{eqnarray*}
For any $w$ with $2w \sigma^2 <1$, we have $\int_{1}^{\infty} t^{-\frac{1}{2w \sigma^2}}dt < \infty$. Then it follows that
$$ E[\sup_{1\leq j \leq 7^m} |Z_j-EZ_j|^2] \leq C m .$$
Thus we obtain $E[\interleave \frac{\partial \psi}{\partial \theta}(\tilde{\theta}_0, x,y)\interleave_2]^2 \leq C m$ and then $E|F(X,Y)|^2 < C m$.

Now we can show that
$$ \mathbb{P}(\sup_{\alpha \in \mathcal{S}^{m-1}, \theta \in H}
|G_n[\alpha^{\top}(\psi_{\theta}-\psi_{\tilde{\theta}_0})]|> 8\epsilon) \rightarrow 0.
$$
Let $N(\epsilon, \mathcal{F}, L_1(Q))$ be the covering number of $\mathcal{F}$ with respect to the seminorm $L_1(Q)$ and let $N(\epsilon, H, \|\cdot\|)$ and $N(\epsilon, \mathcal{S}^{m-1}, \|\cdot\|)$ be the covering number of $H$ and $\mathcal{S}^{m-1}$ respectively, where $\|\cdot\|$ denotes the Frobenious norm. Since
$$|\alpha^{\top}[\psi_{\theta}(x, y)-\psi_{\tilde{\theta}_0}(x, y)]| \leq \|\theta-\tilde{\theta}_0\| F(x,y),$$
it follows that
$$ N(2\epsilon QF, \mathcal{F}_n, L_1(Q)) \leq N(\epsilon, H, \|\cdot\|) \cdot N(\epsilon, \mathcal{S}^{m-1}, \|\cdot\|). $$
Similar to the arguments of Lemma 1.18 of Rigollet and H${\rm \ddot{u}}$tter (2017), we have
$$
N(\epsilon,H,\|\cdot\|) \leq \left(\frac{3}{\epsilon} \sqrt{\frac{m\log n}{n}}\right)^m \ {\rm and} \ N(\epsilon, \mathcal{S}^{m-1}, \|\cdot\|) \leq \left(\frac{3}{\epsilon} \right)^m ,
$$
whence
$$N(2\epsilon QF, \mathcal{F}, L_1(Q)) \leq \left(\frac{9}{\epsilon^2} \sqrt{\frac{m \log n}{n}}\right)^m.$$
Let $\epsilon_n=\epsilon/\sqrt{n}$ and $\delta_n=\sqrt{(m\log n)/n}$. Since $\|\theta-\tilde{\theta}_0\| \leq \delta_n$, it follows that
\begin{eqnarray*}
\frac{Var\{P_n[\alpha^{\top}(\psi_{\theta}-\psi_{\tilde{\theta}_0})]\}}{(4\epsilon_n)^2}
\leq \frac{E[\|\theta-\tilde{\theta}_0\|^2 F(x,y)^2]}{16\epsilon^2} \leq \frac{C(m^2 \log n)/n}{16\epsilon^2} \rightarrow 0.
\end{eqnarray*}
Similar to the arguments of Lemma 2 of Tan and Zhu (2019b), we have
\begin{eqnarray*}
&&\mathbb{P}(\sup_{\alpha \in \mathcal{S}^{m-1}, \theta \in H}|G_n[\alpha^{\top}(\psi_{\theta}- \psi_{\tilde{\theta}_0})]|> 8\epsilon)\\
&=&\mathbb{P}(\sup_{\alpha \in \mathcal{S}^{m-1},\theta \in H}
|P_n[\alpha^{\top}(\psi_{\theta}-\psi_{\tilde{\theta}_0})]-P[\alpha^{\top}(\psi_{\theta}- \psi_{\tilde{\theta}_0})]|> 8 \epsilon_n)\\
&\leq& 4 E\{2N(\epsilon_n, \mathcal{F}_n, L_1(P_n)) \exp[-\frac{1}{2}n \epsilon_n^2/(\delta_n^2 P_n F^2)] \wedge 1 \} \\
&\leq& 8 E\{\left( \frac{36 \delta_n (P_nF)^2}{\epsilon_n^2} \right)^m \exp[-\frac{1}{2}n \epsilon_n^2/(\delta_n^2 P_n F^2)] \wedge 1\} \\
&\leq& 8 \left( \frac{36 \delta_n \tau_n^2}{\epsilon_n^2}\right)^m \exp\{-\frac{1}{2}n \epsilon_n^2/(\delta_n^2 \tau_n^2)\} + 4\mathbb{P}(P_nF^2 > \tau_n^2).
\end{eqnarray*}
Let $\tau_n=\sqrt{m \log{n}}$. Then it follows that
$$\mathbb{P}(P_nF^2 > \tau_n^2) \leq \frac{E(P_n F^2)}{\tau_n^2}=
\frac{E[F(X,Y)^2]}{\tau_n^2} \rightarrow 0. $$
Since $(m \log{n})^3/n \to 0$, it follows that
$$ \left( \frac{36 \delta_n \tau_n^2}{\epsilon_n^2} \right)^q \exp\{-\frac{1}{2}n \epsilon_n^2/(\delta_n^2 \tau_n^2)\} \rightarrow 0. $$
Thus we obtain that $\mathbb{P}(\sup_{\alpha \in \mathcal{S}^{m-1}, \theta \in H}|G_n[\alpha^{\top}(\psi_{\theta}- \psi_{\tilde{\theta}_0})]|> 8\epsilon)=o(1)$ and then $ G_n[\alpha^{\top}(\psi_{\hat{\theta}}-\psi_{\tilde{\theta}_0})]=o_p(1)$ uniformly in $\alpha \in \mathcal{S}^{m-1}$.

Note $P_n\psi_{\hat{\theta}}=0$ and $P\psi_{\tilde{\theta}_0}=0$. Then it follows that
\begin{eqnarray*}
G_n(\alpha^{\top} \psi_{\hat{\theta}})
&=&\sqrt{n}(P_n \alpha^{\top} \psi_{\hat{\theta}}-P\alpha^{\top} \psi_{\hat{\theta}})
=-\sqrt{n}(P\alpha^{\top}\psi_{\hat{\theta}}-P\alpha^{\top} \psi_{\tilde{\theta}_0})\\
&=& G_n(\alpha^{\top} \psi_{\tilde{\theta}_0})+G_n[\alpha^{\top} (\psi_{\hat{\theta}}- \psi_{\tilde{\theta}_0})]=G_n(\alpha^{\top} \psi_{\tilde{\theta}_0})+o_p(1)
\end{eqnarray*}
Thus we obtain $G_n(\alpha^{\top} \psi_{\tilde{\theta}_0})=
-\sqrt{n}(P\alpha^{\top} \psi_{\hat{\theta}}-P\alpha^{\top} \psi_{\tilde{\theta}_0})+o_p(1)$.
Replacing $\alpha$ with $\Sigma^{-1}\alpha$, it follows that
$$ \alpha^{\top} \Sigma^{-1} G_n(\psi_{\tilde{\theta}_0})=-\sqrt{n}\alpha^{\top} \Sigma^{-1}(P\psi_{\hat{\theta}}-P\psi_{\tilde{\theta}_0})+o_p(1). $$
Applying  Taylor expansion, we have
\begin{eqnarray*}
P\psi_{\hat{\theta}}-P\psi_{\tilde{\theta}_0}&=&\dot{P}\psi_{\theta}|_{\theta=\tilde{\theta}_0}
(\hat{\theta}-\tilde{\theta}_0)+(\hat{\theta}-\tilde{\theta}_0)^{\top} (\ddot{P}\psi_{\theta}|_{\theta=\tilde{\theta}_2})(\hat{\theta}-\tilde{\theta}_0)\\
&=& -\Sigma (\hat{\theta}-\tilde{\theta}_0)+(\hat{\theta}-\tilde{\theta}_0)^{\top} (\ddot{P}\psi_{\theta}|_{\theta=\tilde{\theta}_2})(\hat{\theta}-\tilde{\theta}_0)
\end{eqnarray*}
where $\tilde{\theta}_2$ lies between $\hat{\theta}$ and $\tilde{\theta}_0$. According to Assumption 2, we have
$$(\hat{\theta}-\tilde{\theta}_0)^{\top}(\ddot{P}\psi_{\theta}|_{\theta=\tilde{\theta}_2})
(\hat{\theta}-\tilde{\theta}_0) = \frac{m^{3/2}}{n} O_p(1), $$
whence,
$$ \sqrt{n} \alpha^{\top}(\hat{\theta}-\tilde{\theta}_0)=\alpha^{\top} \Sigma^{-1} G_n(\psi_{\tilde{\theta}_0}) + \sqrt{\frac{m^3}{n}} O_p(1)+o_p(1). $$
Since $ (m^3\log n)/n \to 0 $, it follows that
$$ \alpha^{\top}(\hat{\theta}-\tilde{\theta}_0)=\frac{1}{n} \sum_{i=1}^{n} [Y_i-g(\theta, X_i)]\alpha^{\top} \Sigma^{-1}g'(\theta, X_i)+o_p(\frac{1}{\sqrt{n}}). $$
Hence we complete the proof.  \hfill$\Box$

\textbf{Proof of Theorem~\ref{theorem2.3}.} We use the same arguments for Theorem~\ref{theorem2.1}. For the linear model, we have $\psi_{\theta}(X,Y)=X(Y-\theta^{\top}X)$. Thus
$$ \psi_{\tilde{\theta}_0+\zeta}(X,Y)-\psi_{\tilde{\theta}_0}(X,Y)=-XX^{\top}\zeta, $$
whence
$$P_n(\psi_{\tilde{\theta}_0+\zeta}-\psi_{\tilde{\theta}_0})-
P(\psi_{\tilde{\theta}_0+\zeta}-\psi_{\tilde{\theta}_0})=-(\frac{1}{n}\sum_{i=1}^n X_iX_i^{\top}-EXX^{\top})\zeta.$$
Next we show that
$$\interleave \frac{1}{n}\sum_{i=1}^n X_iX_i^{\top}-EXX^{\top}\interleave_2 =O_p(\sqrt{\frac{p}{n}}), $$
where $\interleave \cdot \interleave_2$ is the operator norm of a matrix $A$ that is define in the proof of Theorem~\ref{theorem2.2}.
Set $EX=\mu$. Then it follows that
\begin{eqnarray*}
\frac{1}{n}\sum_{i=1}^n X_iX_i^{\top}-EXX^{\top}
&=&\frac{1}{n}\sum_{i=1}^n (X_i-\mu)(X_i-\mu)^{\top}-E(X-\mu)(X-\mu)^{\top}+\\
&& \frac{1}{n}\sum_{i=1}^n \mu X_i^{\top}- \mu \mu^{\top}+\frac{1}{n}\sum_{i=1}^n X_i \mu ^{\top}- \mu \mu^{\top}\\
&=:& A_1+A_2+A_3.
\end{eqnarray*}
According to the proof of Lemma 3 of Cai et al. (2010), there exist vectors $\{v_j \in \mathbb{S}^{p-1}, 1\leq j \leq 7^p\}$ such that
$$ \interleave A \interleave_2 \leq 4 \sup_{1\leq j \leq 7^p} |v_j^{\top} A v_j| $$
for any symmetric matrix $A$, where $\mathbb{S}^{p-1}$ is the unit sphere with respect to the Frobenious norm.
Then it follows that
\begin{eqnarray*}
&& \mathbb{P}(\interleave \frac{1}{n}\sum_{i=1}^n (X_i-\mu)(X_i-\mu)^{\top}- E(X-\mu)(X-\mu)^{\top}\interleave_2>t) \\
&\leq&  \mathbb{P}(\sup_{1\leq j \leq 7^p} |\frac{1}{n}\sum_{i=1}^n v_j^{\top}(X_i-\mu)(X_i-\mu)^{\top}v_j-
Ev_j^{\top}(X-\mu)(X-\mu)^{\top}v_j|>t)      \\
&=&  \mathbb{P}(\sup_{1\leq j \leq 7^p} |\frac{1}{n}\sum_{i=1}^n |v_j^{\top}(X_i-\mu)|^2- E|v_j^{\top}(X-\mu)|^2|>t).
\end{eqnarray*}
Note that $v^{\top}(X-\mu)$ is sub-Gaussian with proxy variance $\sigma^2$. By Lemma 1.12 of Rigollet and H${\rm \ddot{u}}$tter (2017), we obtain that $|v^{\top}(X-\mu)|^2-E|v^{\top}(X-\mu)|^2 $ is sub-exponential with parameter $\lambda=16\sigma^2$. By Berstein's inequality, see e.g., Theorem 1.13 of Rigollet and H${\rm \ddot{u}}$tter (2017),
$$ \mathbb{P}(|\frac{1}{n}\sum_{i=1}^n [v^{\top}(X_i-\mu)]^2- E[v^{\top}(X-\mu)]^2 | >t) < 2\exp\{-\frac{1}{2}n(\frac{t^2}{\lambda^2} \wedge \frac{t}{\lambda})  \},$$
where $a \wedge b=\min(a,b)$.
Let $t=\sqrt{p/n}M$ and then $t/\lambda < 1$ for $n$ large enough. Thus $t^2/\lambda^2 < t/\lambda$. Therefore,
\begin{eqnarray*}
&& \mathbb{P}(\interleave \frac{1}{n}\sum_{i=1}^n (X_i-\mu)(X_i-\mu)^{\top}- E(X-\mu)(X-\mu)^{\top}\interleave_2>\sqrt{\frac{p}{n}}M) \\
&\leq&  \sum_{j=1}^{7^p}\mathbb{P}(|\frac{1}{n}\sum_{i=1}^n [v_j^{\top}(X_i-\mu)]^2- E[v_j^{\top}(X-\mu)]^2|>\sqrt{\frac{p}{n}}M) \\
&\leq&  7^p \cdot 2\exp(-\frac{pM^2}{2\lambda^2})=
2\exp\{p(\log7-\frac{M^2}{2\lambda^2})\}.
\end{eqnarray*}
Then  $\interleave A_1 \interleave_2=O_p(\sqrt{p/n})$.

For the second term $A_2$, we have
\begin{eqnarray*}
&& \mathbb{P}(\interleave \frac{1}{n}\sum_{i=1}^n \mu (X_i-\mu)^{\top} \interleave_2 >\sqrt{\frac{p}{n}}M) \\
&\leq& \sum_{j=1}^{7^p}\mathbb{P}(| \frac{1}{n}\sum_{i=1}^n v_j^{\top} \mu (X_i-\mu)^{\top}v_j | >\sqrt{\frac{p}{n}}M) \\
&=&    \sum_{j=1}^{7^p}\mathbb{P}(| \frac{1}{n}\sum_{i=1}^n v_j^{\top}(X_i-\mu) | >\sqrt{\frac{p}{n}} \frac{M}{|v_j^{\top} \mu|}).
\end{eqnarray*}
Since $v_j^{\top}(X_i-\mu)$ is sub-Gaussian with proxy variance $\sigma^2$, by Corollary 1.7 of Rigollet and H${\rm \ddot{u}}$tter (2017), we obtain
$$\mathbb{P}(| \frac{1}{n}\sum_{i=1}^n v_j^{\top}(X_i-\mu) | >t) \leq 2\exp(-\frac{nt^2}{2\sigma^2}),  $$
whence,
\begin{eqnarray*}
\mathbb{P}(\interleave \frac{1}{n}\sum_{i=1}^n \mu (X_i-\mu)^{\top} \interleave_2 >\sqrt{\frac{p}{n}}M) \leq 2\exp\{p(\log7-\frac{M^2}{2\sigma^2|v_j^{\top} \mu|^2})\}.
\end{eqnarray*}
Since $v_j^{\top} E(X-\mu)(X-\mu)^{\top}v_j \geq 0$, it follows that $v_j^{\top}\mu\mu^{\top}v_j \leq v_j^{\top}E(XX^{\top})v_j$. Note that $\Sigma = E(XX^{\top})$. Thus we obtain that
$$ |v_j^{\top}\mu|^2 \leq v_j^{\top} \Sigma v_j \leq \lambda_{\max}(\Sigma) \leq C. $$
Therefore
$$\mathbb{P}(\interleave \frac{1}{n}\sum_{i=1}^n \mu (X_i-\mu)^{\top} \interleave_2 >\sqrt{\frac{p}{n}}M) \leq 2\exp\{p(\log7-\frac{M^2}{2 C \sigma^2})\}.  $$
Hence  $\interleave A_2 \interleave_2=O_p(\sqrt{p/n})$. Similarly, we have $\interleave A_3 \interleave_2=O_p(\sqrt{p/n})$.
Consequently,
$$n \zeta^{\top}[P_n(\psi_{\tilde{\theta}_0+\zeta}-\psi_{\tilde{\theta}_0})-
P(\psi_{\tilde{\theta}_0+\zeta}-\psi_{\tilde{\theta}_0})]=\sqrt{pn}\|\zeta\|^2O_p(1).
$$
Similar to the arguments for Theorem~\ref{theorem2.1}, we derive
$$ n \zeta^{\top} P_n\psi_{\tilde{\theta}_0} = \sqrt{pn} \|\zeta\|O_p(1) \quad {\rm and} \quad n\zeta^{\top} P\psi_{\tilde{\theta}_0+\zeta}=-n\zeta^{\top}\Sigma\zeta. $$
Let $\zeta=\sqrt{p/n}U $. If $p/n \to 0$, then we have
\begin{eqnarray*}
\zeta^{\top}F(\zeta)
&=&\sqrt{pn} \|\zeta\|O_p(1)-n\zeta^{\top}\Sigma\zeta+\sqrt{pn}\|\zeta\|^2O_p(1) \\
&\leq& p\|U\|O_p(1)-p\lambda_{\min}(\Sigma)\|U\|^2+
\sqrt{\frac{p^{3/2}}{n}}\|U\|^2O_p(1) \\
&=& p\|U\| \{ O_p(1)-\lambda_{\min}(\Sigma)\|U\| + \|U\|o_p(1)\}.
\end{eqnarray*}
The rest of the proof follows the same arguments for Theorem~\ref{theorem2.1}. Hence we obtain the norm consistency of $\hat{\theta}$ to $\tilde{\theta}_0$.

Now we derive the uniformly asymptotically linear representation of
$\hat{\theta}-\tilde{\theta}_0$. Recall that $\psi_{\theta}(x,y)=x(y-\theta^{\top}x)$ and $\Sigma=E(XX^{\top})$. Since $$\psi_{\theta}(x,y)-\psi_{\tilde{\theta}_0}(x,y)=-xx^{\top}(\theta-\tilde{\theta}_0), $$
it follows that
\begin{eqnarray*}
G_n[\alpha^{\top}\Sigma^{-1}(\psi_{\hat{\theta}}-\psi_{\tilde{\theta}_0})]
&=& \sqrt{n}\alpha^{\top}\Sigma^{-1}[P_n(\psi_{\hat{\theta}}-\psi_{\tilde{\theta}_0})-
P(\psi_{\hat{\theta}}-\psi_{\tilde{\theta}_0})]\\
&=& -\sqrt{n}\alpha^{\top}\Sigma^{-1}(\frac{1}{n} \sum_{i=1}^{n}X_iX_i^{\top}-\Sigma)
(\hat{\theta}-\tilde{\theta}_0).
\end{eqnarray*}
Consequently,
\begin{eqnarray*}
|G_n[\alpha^{\top}\Sigma^{-1}(\psi_{\hat{\theta}}-\psi_{\tilde{\theta}_0})]|
&\leq& \sqrt{n} \|\alpha^{\top}\Sigma^{-1}\| \cdot \interleave \frac{1}{n} \sum_{i=1}^{n} X_iX_i^{\top}-\Sigma \interleave_2 \cdot \|\hat{\theta}-\tilde{\theta}_0\|  \\
&=& O_p(\sqrt{\frac{p^2}{n}})=o_p(1).
\end{eqnarray*}
Thus $G_n[\alpha^{\top}\Sigma^{-1}(\psi_{\hat{\theta}}-\psi_{\tilde{\theta}_0})]=o_p(1)$
uniformly in $\alpha \in \mathcal{S}^{q-1}$.
Note $P_n\psi_{\hat{\theta}}=0$ and $P\psi_{\tilde{\theta}_0}=0$. Then we have
\begin{eqnarray*}
\alpha^{\top}\Sigma^{-1}G_n\psi_{\hat{\theta}}
&=& \alpha^{\top}\Sigma^{-1}G_n\psi_{\tilde{\theta}_0}+o_p(1)
=\frac{1}{\sqrt{n}}\sum_{i=1}^{n}(Y_i-\tilde{\theta}_0^{\top}X_i) \alpha^{\top}\Sigma^{-1}X_i+o_p(1) \\
&=&-\sqrt{n}\alpha^{\top}\Sigma^{-1} P(\psi_{\hat{\theta}}-\psi_{\tilde{\theta}_0})
=\sqrt{n} \alpha^{\top} (\hat{\theta}-\tilde{\theta}_0),
\end{eqnarray*}
whence,
$$\sqrt{n}\alpha^{\top} (\hat{\theta}-\tilde{\theta}_0)=\alpha^{\top}\Sigma^{-1}
\frac{1}{\sqrt{n}}\sum_{i=1}^{n}(Y_i-\tilde{\theta}_0^{\top}X_i)X_i+o_p(1). $$
Hence we complete the proof.  \hfill$\Box$

\textbf{Proof of Theorem~\ref{theorem3.1}.} This proof is similar to that for proving Proposition~3 in Tan and Zhu (2019a). Thus we omit the details here.

\textbf{Proof of Theorem~\ref{theorem4.1}.}
%For the simplicity of notations, we only prove the theorem for parametric single index models. The arguments for multiple index models are similar.
Under the null hypothesis, Theorem~~\ref{theorem3.1} shows that $\mathbb{P}(\hat{q}=d) \to 1$. Thus we only need to work on the event $\{ \hat{q}=d \}$. Recall that
$$ \hat{V}_n^1(t) = \frac{1}{\sqrt{n}} \sum_{j=1}^{n} [Y_j-g(\hat{\bm{\beta}}^{\top}X_j, \hat{\vartheta})] [\cos(t^{\top} \hat{\bm{\mathfrak B}}^{\top}X_j)+\sin(t^{\top} \hat{\bm{\mathfrak B}}^{\top}X_j)].
$$
Here $t \in \mathbb{R}^{2d}$ and $\hat{\bm{\mathfrak B}}=(\bm{\hat{\beta}}, \bm{\hat{B}}) \in \mathbb{R}^{p\times 2d} $ is a random matrix. Let $\hat{\theta}=[vec(\hat{\bm{\beta}})^{\top},\hat{\vartheta}^{\top}]^{\top}$ as we have defined before. Then we rewrite $\hat{V}_n^1(t)$ as
$$ \hat{V}_n^1(t) = \frac{1}{\sqrt{n}} \sum_{j=1}^{n} [Y_j-g(\hat{\theta}, X_j)] [\cos(t^{\top} \hat{\bm{\mathfrak B}}^{\top}X_j)+\sin(t^{\top} \hat{\bm{\mathfrak B}}^{\top}X_j)].
$$
In the following we decompose the proof into three steps.

{\it Step 1}. To show that
$$ \hat{V}_n^1(t) = V_n^1(t)+R_n(t) \quad {\rm and} \quad \int_{\mathbb{R}^2} |R_n(t)|^2 \varphi(t)dt=o_p(1).
$$
Here $R_n(t)$ is a remainder and
$$ V_n^1(t)=  \frac{1}{\sqrt{n}} \sum_{j=1}^{n} \varepsilon_j [\cos(t^{\top} \bm{\mathfrak{B}_0}^{\top}X_j)+\sin(t^{\top} \bm{\mathfrak{B}_0}^{\top}X_j)]-M(t)^{\top} \Sigma^{-1} \frac{1}{\sqrt{n}}\sum_{j=1}^n \varepsilon_j g'(\theta_0, X_j)X_j
$$
with $\bm{\mathfrak{B}_0}=(\bm{\beta_0}, \bm{B})$ and $M(t)=E\{g'(\theta_0,X) [\cos(t^{\top}\bm{\mathfrak{B}_0}^{\top}X)+\sin(t^{\top}\bm{\mathfrak{B}_0}^{\top}X)]\}$.

Decompose $\hat{V}_n^1(t)$ as
\begin{eqnarray*}
\hat{V}_n^1(t) &=& \frac{1}{\sqrt{n}} \sum_{j=1}^{n} \varepsilon_j [\cos(t^{\top} \bm{\hat{\mathfrak B}}^{\top}X_j)+\sin(t^{\top} \bm{\hat{\mathfrak B}}^{\top}X_j)] \\
&& -\frac{1}{\sqrt{n}} \sum_{j=1}^{n} [g(\hat{\theta}, X_j)-g(\theta_0, X_j)] [\cos(t^{\top} \bm{\hat{\mathfrak B}}^{\top}X_j)+\sin(t^{\top} \bm{\bm{\hat{\mathfrak B}}}^{\top}X_j)] \\
&=& \frac{1}{\sqrt{n}} \sum_{j=1}^{n} \varepsilon_j [\cos(t^{\top} \bm{\mathfrak{B}_0}^{\top}X_j)+\sin(t^{\top} \bm{\mathfrak{B}_0}^{\top}X_j)] \\
&&  +\frac{1}{\sqrt{n}} \sum_{j=1}^{n} \varepsilon_j [\cos(t^{\top} \bm{\hat{\mathfrak B}}^{\top}X_j)-\cos(t^{\top} \bm{\mathfrak{B}_0}^{\top}X_j) + \sin(t^{\top} \bm{\hat{\mathfrak B}}^{\top}X_j)-\sin(t^{\top} \bm{\mathfrak{B}_0}^{\top}X_j)]  \\
&& -\frac{1}{\sqrt{n}} \sum_{j=1}^{n} [g(\hat{\theta}, X_j)-g(\theta_0, X_j)] [\cos(t^{\top} \bm{\mathfrak{B}_0}^{\top}X_j)+\sin(t^{\top} \bm{\mathfrak{B}_0}^{\top}X_j)]  \\
&& -\frac{1}{\sqrt{n}} \sum_{j=1}^{n} [g(\hat{\theta}, X_j)-g(\theta_0, X_j)] [\cos(t^{\top} \bm{\hat{\mathfrak B}}^{\top}X_j)-\cos(t^{\top} \bm{\mathfrak{B}_0}^{\top}X_j) + \sin(t^{\top} \bm{\hat{\mathfrak B}}^{\top}X_j)-\sin(t^{\top} \bm{\mathfrak{B}_0}^{\top}X_j)] \\
&=:& V_{n1}(t)+V_{n2}(t)-V_{n3}(t)-V_{n4}(t),
\end{eqnarray*}

First we show that
$$\int_{\mathbb{R}^2} V_{n2}(t)^2 \varphi(t)dt =o_p(1) \quad {\rm and } \quad \int_{\mathbb{R}^2} V_{n4}(t)^2 \varphi(t)dt =o_p(1).$$
Decompose further $V_{n2}(t)$ as
\begin{eqnarray*}
V_{n2}(t) &=& \frac{1}{\sqrt{n}} \sum_{j=1}^{n} \varepsilon_j [\cos(t^{\top} \bm{\hat{\mathfrak B}}^{\top}X_j)-\cos(t^{\top}
              \bm{\mathfrak{B}_0}^{\top}X_j)]+ \\
&&            \frac{1}{\sqrt{n}} \sum_{j=1}^{n} \varepsilon_j [\sin(t^{\top} \bm{\hat{\mathfrak B}}^{\top}X_j)-\sin(t^{\top}
              \bm{\mathfrak{B}_0}^{\top}X_j)]\\
&=:& V_{n21}(t)+V_{n22}(t)
\end{eqnarray*}
For the term $V_{n21}(t)$, by Taylor expansion, we have
\begin{eqnarray*}
V_{n21}(t) &=& -\frac{1}{\sqrt{n}} \sum_{j=1}^{n} \varepsilon_j (t^{\top} \bm{\hat{\mathfrak B}}^{\top}X_j-t^{\top}
               \bm{\mathfrak{B}_0}^{\top}X_j)\sin(t^{\top} \bm{\mathfrak{B}_0}^{\top}X_j)\\
&&             -\frac{1}{2 \sqrt{n}} \sum_{j=1}^{n} \varepsilon_j (t^{\top} \bm{\hat{\mathfrak B}}^{\top}X_j-t^{\top} \bm{\mathfrak{B}_0}^{\top}X_j)^2 \cos(t^{\top} \bm{\mathfrak{B}_0}^{\top}X_j)\\
&&             + \frac{1}{6 \sqrt{n}} \sum_{j=1}^{n} \varepsilon_j (t^{\top} \bm{\hat{\mathfrak B}}^{\top}X_j-t^{\top} \bm{\mathfrak{B}_0}^{\top}X_j)^3
              \sin(t^{\top} \bm{\mathfrak{B}}_1^{\top}X_j) \\
&=:&  -V_{n211}(t)-V_{n212}(t)+V_{n213}(t)
\end{eqnarray*}
where $\bm{\mathfrak{B}}_1$ lies between $\bm{\hat{\mathfrak B}}$ and $\bm{\mathfrak{B}_0}$. Then it follows that
\begin{eqnarray*}
\int_{\mathbb{R}^{2d}} |V_{n211}(t)|^2 \varphi(t)dt
\leq \|\bm{\hat{\mathfrak B}}-\bm{\mathfrak{B}_0}\|^2 \int_{\mathbb{R}^{2d}} \|t\|^2 \cdot \|\frac{1}{\sqrt{n}} \sum_{j=1}^{n} \varepsilon_j \sin(t^{\top} \bm{\mathfrak{B}_0}^{\top}X_j) X_j \|^2 \varphi(t)dt.
\end{eqnarray*}
Note that
\begin{eqnarray*}
&& E\{ \int_{\mathbb{R}^{2d}} \|t\|^2\|\frac{1}{\sqrt{n}} \sum_{j=1}^{n} \varepsilon_j \sin(t^{\top} \bm{\mathfrak{B}_0}^{\top}X_j) X_j \|^2 \varphi(t)dt  \} \\
&=& \int_{\mathbb{R}^{2d}} \|t\|^2 E\{ \|\frac{1}{\sqrt{n}} \sum_{j=1}^{n} \varepsilon_j \sin(t^{\top} \bm{\mathfrak{B}_0}^{\top}X_j) X_j \|^2 \} \varphi(t)dt \\
&=& \int_{\mathbb{R}^{2d}} \|t\|^2 E\{\|X\|^2\varepsilon^2 \sin(t^{\top} \bm{\mathfrak{B}_0}^{\top}X)^2\} \varphi(t)dt \leq C p.
\end{eqnarray*}
Consequently, we obtain that
$$\int_{\mathbb{R}^{2d}} |V_{n211}|^2 \varphi(t)dt = \frac{p^2}{n} O_p(1)=o_p(1).$$
Similarly, we can obtain
$$ \int_{\mathbb{R}^{2d}} |V_{n212}(t)|^2 \varphi(t)dt = \frac{p^4}{n^2} O_p(1)=o_p(1). $$
For the third term $V_{n213}(t)$ in $V_{n21}(t)$, we have
\begin{eqnarray*}
\int_{\mathbb{R}^{2d}} |V_{n213}(t)|^2 \varphi(t)dt &=& \int_{\mathbb{R}^{2d}} |\frac{1}{\sqrt{n}} \sum_{j=1}^{n} \varepsilon_j (t^{\top} \bm{\hat{\mathfrak B}}^{\top}X_j-t^{\top} \bm{\mathfrak{B}_0}^{\top}X_j)^3 \sin(t^{\top} \bm{\mathfrak{B}}_1^{\top}X_j)|^2 \varphi(t)dt \\
&\leq& \int_{\mathbb{R}^{2d}} \{ \frac{1}{\sqrt{n}} \sum_{j=1}^{n} |\varepsilon_j| \cdot |t^{\top} \bm{\hat{\mathfrak B}}^{\top}X_j-t^{\top} \bm{\mathfrak{B}_0}^{\top}X_j|^3 \}^2 \varphi(t)dt \\
&\leq& \|\bm{\hat{\mathfrak B}}-\bm{\mathfrak{B}_0}\|^6 \cdot  \int_{\mathbb{R}^{2d}} \|t\|^6 \varphi(t)dt \cdot \{ \frac{1}{\sqrt{n}} \sum_{j=1}^{n} |\varepsilon_j| \| X_j\|^3 \}^2.
\end{eqnarray*}
Since $E|\varepsilon_j|^8 < C $ and $E|X_{ji}|^8 < C$ for $1 \leq i \leq p$, it follows that
$$ E\{ \frac{1}{\sqrt{n}} \sum_{j=1}^{n} |\varepsilon_j| \| X_j\|^3 \}^2 < C np^3. $$
Combining this with $\|\bm{\hat{\mathfrak B}}-\bm{\mathfrak{B}_0}\|=O_p(\sqrt{p/n})$, we obtain
$$ \int_{\mathbb{R}^{2d}} |V_{n213}(t)|^2 \varphi(t)dt = \frac{p^6}{n^2} O_p(1) =o_p(1) .$$
Consequently, we have $ \int_{\mathbb{R}^{2d}} |V_{n21}(t)|^2 \varphi(t)dt = o_p(1) $. By the same arguments for the term $V_{n22}$, we have
$\int_{\mathbb{R}^{2d}} |V_{n22}(t)|^2 \varphi(t)dt = o_p(1) $.
Then it follows that
$$\int_{\mathbb{R}^{2d}} |V_{n2}(t)|^2 \varphi(t)dt =o_p(1).$$

Now we consider the term $\int_{\mathbb{R}^{2d}} V_{n4}(t)^2 \varphi(t)dt$.  Decompose $V_{n4}(t)$ as follows,
\begin{eqnarray*}
V_{n4}(t)&=& \frac{1}{\sqrt{n}} \sum_{j=1}^{n} [g(\hat{\theta}, X_j)-g(\theta_0, X_j)] [\cos(t^{\top} \bm{\hat{\mathfrak B}}^{\top}X_j)-\cos(t^{\top} \bm{\mathfrak{B}_0}^{\top}X_j)]+ \\
&& \frac{1}{\sqrt{n}} \sum_{j=1}^{n} [g(\hat{\theta}, X_j)-g(\theta_0, X_j)] [\sin(t^{\top} \bm{\hat{\mathfrak B}}^{\top}X_j)-\sin(t^{\top} \bm{\mathfrak{B}_0}^{\top}X_j)] \\
&=:& V_{n41}(t)+V_{n42}(t)
\end{eqnarray*}
By Taylor expansion, we have
\begin{eqnarray*}
V_{n41}(t)&=& \frac{1}{\sqrt{n}} \sum_{j=1}^{n}[(\hat{\theta}-\theta_0)^{\top}g'(\theta_0, X_j)+ \frac{1}{2}(\hat{\theta}-\theta_0)^{\top}g''(\theta_1, X_j)(\hat{\theta}-\theta_0)] \times \\
&&[-t^{\top} (\bm{\hat{\mathfrak B}} - \bm{\mathfrak{B}_0})^{\top} X_j \sin(t^{\top} \bm{{\mathfrak B}_0}^{\top}X_j)-\frac{1}{2} (t^{\top} \bm{\hat{\mathfrak B}}^{\top}X_j - t^{\top} \bm{\mathfrak{B}_0}^{\top}X_j)^2 \sin(t^{\top} \bm{{\mathfrak B}_1}^{\top}X_j)] \\
&=& -\frac{1}{\sqrt{n}} \sum_{j=1}^{n} (\hat{\theta}-\theta_0)^{\top} g'(\theta_0, X_j) X_j^{\top} \sin(t^{\top} {\mathfrak B}_0^{\top}X_j)
(\bm{\hat{\mathfrak B}} -\bm{\mathfrak{B}_0})t \\
&& -\frac{1}{2\sqrt{n}} \sum_{j=1}^{n} (\hat{\theta}-\theta_0)^{\top}g''(\theta_1, X_j)(\hat{\theta}-\theta_0) t^{\top} (\bm{\hat{\mathfrak B}} -  \bm{\mathfrak{B}_0})^{\top} X_j \sin(t^{\top} \bm{{\mathfrak B}_0}^{\top}X_j) \\
&& -\frac{1}{2\sqrt{n}} \sum_{j=1}^{n} (\hat{\theta}-\theta_0)^{\top}g'(\theta_0, X_j) (t^{\top} \bm{\hat{\mathfrak B}}^{\top}X_j - t^{\top} \bm{\mathfrak{B}_0}^{\top}X_j)^2 \sin(t^{\top} \bm{{\mathfrak B}_1}^{\top}X_j) \\
&& -\frac{1}{4\sqrt{n}} \sum_{j=1}^{n} (\hat{\theta}-\theta_0)^{\top}g''(\theta_1, X_j)(\hat{\theta}-\theta_0) (t^{\top} \bm{\hat{\mathfrak B}}^{\top}X_j - t^{\top} \bm{\mathfrak{B}_0}^{\top}X_j)^2 \sin(t^{\top} \bm{{\mathfrak B}_1}^{\top}X_j)\\
&=:& -V_{n411}(t)-V_{n412}(t)-V_{n413}(t)-V_{n414}(t),
\end{eqnarray*}
where $\theta_1$ lies between $\hat{\theta}$ and $\theta_0$ and $\bm{{\mathfrak B}_1}$ lies between $\hat{\bm{\mathfrak{B}}}$ and $\bm{{\mathfrak B}_0}$. By the same arguments as that for the term $V_{n213}$, we have
$$ \int_{\mathbb{R}^{2d}} |V_{n41k}|^2 \varphi(t)dt=o_p(1), \quad \ {\rm for } \ j=2,3,4. $$
For the term $V_{n411}$, set
$$ M_1(t)=E[g'(\theta_0, X) X^{\top} \sin(t^{\top} \bm{{\mathfrak B}_0}^{\top}X)]. $$
Then we have
\begin{eqnarray*}
&& \int_{\mathbb{R}^{2d}} |V_{n411}(t)|^2 \varphi(t)dt \\
&\leq& 2 n \int_{\mathbb{R}^{2d}} |(\hat{\theta}-\theta_0)^{\top} \{\frac{1}{n} \sum_{j=1}^{n} g'(\theta_0, X_j)X_j^{\top} \sin(t^{\top} \bm{{\mathfrak B}_0}^{\top}X_j)- M_1(t)\} \times \\
&& (\bm{\hat{\mathfrak B}} -\bm{\mathfrak{B}_0})t|^2 \varphi(t)dt + 2 n \int_{\mathbb{R}^{2d}} |(\hat{\theta}-\theta_0)^{\top} M_1(t) (\bm{\hat{\mathfrak B}} -\bm{\mathfrak{B}_0})t|^2 \varphi(t)dt \\
&\leq& 2 n \|\hat{\theta}-\theta_0\|^2 \|\bm{\hat{\mathfrak B}} -\bm{\mathfrak{B}_0} \|^2 \int_{\mathbb{R}^{2d}} \| \frac{1}{n} \sum_{j=1}^{n}  g'(\theta_0, X_j)X_j^{\top} \sin(t^{\top} \bm{{\mathfrak B}_0}^{\top}X_j)- M_1(t)\|^2 \times \\
&& \|t\|^2 \varphi(t)dt + C n \|\hat{\theta}-\theta_0\|^2 \|\bm{\hat{\mathfrak B}} -\bm{\mathfrak{B}_0} \|^2.
\end{eqnarray*}
Note that
\begin{eqnarray*}
&& E \left( \int_{\mathbb{R}^{2d}} \| \frac{1}{n} \sum_{j=1}^{n} g'(\theta_0, X_j)X_j^{\top} \sin(t^{\top}\bm{{\mathfrak B}_0}^{\top}X_j)- M_1(t)\|^2 \|t\|^2 \varphi(t)dt \right) \\
&=& \sum_{k=1}^{m}\sum_{l=1}^{p} \int_{\mathbb{R}^{2d}} E | \frac{1}{n} \sum_{j=1}^{n}  g'_k(\theta_0, X_j)X_{jl} \sin(t^{\top} \bm{{\mathfrak B}_0}^{\top}X_j)- M_{1kl}(t) |^2 \|t\|^2 \varphi(t)dt \\
&\leq& \frac{1}{n} \sum_{k=1}^{m}\sum_{l=1}^{p}  E|g'_k(\theta_0, X_j)X_{jl}|^2 \leq C \frac{p^2}{n}.
\end{eqnarray*}
Thus we obtain $ \int_{\mathbb{R}^{2d}} |V_{n411}(t)|^2 \varphi(t)dt = (p^2/n) O_p(1)=o_p(1)$ and then $\int_{\mathbb{R}^{2d}} |V_{n41}(t)|^2 \varphi(t)dt = o_p(1)$. The same arguments show that $\int_{\mathbb{R}^{2d}} |V_{n42}(t)|^2 \varphi(t)dt = o_p(1)$. Hence we have
$$\int_{\mathbb{R}^{2d}} |V_{n4}(t)|^2 \varphi(t)dt = o_p(1).  $$

For the term $V_{n3}(t)$, recall that
$$ V_{n3}(t)=\frac{1}{\sqrt{n}} \sum_{j=1}^{n} [g(\hat{\theta}, X_j)-g(\theta_0, X_j)] [\cos(t^{\top} \bm{\mathfrak{B}_0}^{\top}X_j)+\sin(t^{\top} \bm{\mathfrak{B}_0}^{\top}X_j)].$$
Decompose $V_{n3}(t)$ as follows,
\begin{eqnarray*}
V_{n3}(t)&=& \frac{1}{\sqrt{n}} \sum_{j=1}^{n} (\hat{\theta}-\theta_0)^{\top} g'(\theta_0, X_j)[\cos(t^{\top} \bm{\mathfrak{B}_0}^{\top}X_j)+\sin(t^{\top} \bm{\mathfrak{B}_0}^{\top}X_j)] \\
&& + \frac{1}{2\sqrt{n}} \sum_{j=1}^{n} (\hat{\theta}-\theta_0)^{\top} g''(\theta_0, X_j)(\hat{\theta}-\theta_0)[\cos(t^{\top} \bm{\mathfrak{B}_0}^{\top}X_j)+\sin(t^{\top} \bm{\mathfrak{B}_0}^{\top}X_j)] \\
&& + \frac{1}{6\sqrt{n}} \sum_{j=1}^{n} (\hat{\theta}-\theta_0)^{\top}[(\hat{\theta}-\theta_0)^{\top} g'''(\theta_2, X_j) (\hat{\theta}-\theta_0)] [\cos(t^{\top} \bm{\mathfrak{B}_0}^{\top}X_j)+ \sin(t^{\top} \bm{\mathfrak{B}_0}^{\top}X_j)]\\
&=:& V_{n31}(t)+\frac{1}{2}V_{n32}(t)+\frac{1}{6} V_{n33}(t),
\end{eqnarray*}
where $\theta_2$ lies between $\hat{\theta}$ and $\theta_0$. Similar to the arguments as that for $V_{n213}$ and $V_{n411}$, we have
\begin{eqnarray*}
\int_{\mathbb{R}^{2d}}|V_{n33}(t)|^2 \varphi(t)dt &=& \frac{p^6}{n^2} O_p(1)=o_p(1); \\
\int_{\mathbb{R}^{2d}}|V_{n32}(t)|^2 \varphi(t)dt &=& \frac{p^2}{n} O_p(1)=o_p(1);
\end{eqnarray*}
For the first term $V_{n31}(t)$ in $V_{n3}(t)$, recall that
$$M(t)=E\{g'(\theta_0, X)[\cos(t^{\top} \bm{\mathfrak{B}_0}^{\top}X)+\sin(t^{\top} \bm{\mathfrak{B}_0}^{\top}X)]\}.$$
Then it follows that
\begin{eqnarray*}
V_{n31}(t) &=& \sqrt{n} (\hat{\theta}-\theta_0)^{\top}\{\frac{1}{n} \sum_{j=1}^{n} g'(\theta_0, X_j)[\cos(t^{\top} \bm{\mathfrak{B}_0}^{\top}X_j)+\sin(t^{\top} \bm{\mathfrak{B}_0}^{\top}X_j)]-M(t)\} \\
&& \sqrt{n} (\hat{\theta}-\theta_0)^{\top}M(t)\\
&=:& V_{n311}(t) + \sqrt{n} (\hat{\theta}-\theta_0)^{\top}M(t).
\end{eqnarray*}
Note that
\begin{eqnarray*}
&& \int_{\mathbb{R}^{2d}}|V_{n311}(t)|^2 \varphi(t)dt  \\
&\leq& n \|\hat{\theta}-\theta_0\|^2 \int_{\mathbb{R}^{2d}} \|\frac{1}{n} \sum_{j=1}^{n} g'(\theta_0, X_j)[\cos(t^{\top}\bm{\mathfrak{B}_0}^{\top}X_j) +\sin(t^{\top} \bm{\mathfrak{B}_0}^{\top}X_j)]-M(t)\|^2 \varphi(t)dt
\end{eqnarray*}
and
\begin{eqnarray*}
&&  E\left(\int_{\mathbb{R}^{2d}} \|\frac{1}{n} \sum_{j=1}^{n} g'(\theta_0, X_j) [\cos(t^{\top} \bm{\mathfrak{B}_0}^{\top}X_j)+\sin(t^{\top} \bm{\mathfrak{B}_0}^{\top}X_j)]-M(t)\|^2 \varphi(t)dt \right) \\
&=& \int_{\mathbb{R}^{2d}} E\{ \|\frac{1}{n} \sum_{j=1}^{n} g'(\theta_0, X_j) [\cos(t^{\top} \bm{\mathfrak{B}_0}^{\top}X_j)+\sin(t^{\top} \bm{\mathfrak{B}_0}^{\top}X_j)]-M(t)\|^2\} \varphi(t)dt \\
&=& \sum_{k=1}^m \int_{\mathbb{R}^{2d}} E\{ |\frac{1}{n} \sum_{j=1}^{n} g'_k(\theta_0, X_j)[\cos(t^{\top} \bm{\mathfrak{B}_0}^{\top}X_j)+\sin(t^{\top} \bm{\mathfrak{B}_0}^{\top}X_j)]-M_k(t) |^2\} \varphi(t)dt \\
&\leq& \frac{1}{n} \sum_{k=1}^m E|g'_k(\theta_0, X)|^2 \leq C \frac{p}{n}.
\end{eqnarray*}
Thus we obtain that
$$\int_{\mathbb{R}^{2d}}|V_{n311}(t)|^2 \varphi(t)dt = \frac{p^2}{n} O_p(1)=o_p(1).$$
By Theorem~\ref{theorem2.2}, we have
$$ \sqrt{n} M(t)^{\top} (\hat{\theta}-\theta_0)=M(t)^{\top} \Sigma^{-1} \frac{1}{\sqrt{n}}\sum_{j=1}^n \varepsilon_j g'(\theta_0, X_j) + o_p(1), $$
where $o_p(1)$ is uniformly in $t$. Altogether we obtain
\begin{eqnarray*}
\hat{V}_n^1(t)&=&\frac{1}{\sqrt{n}} \sum_{j=1}^{n} \varepsilon_j [\cos(t^{\top} \bm{\mathfrak{B}_0}^{\top}X_j)+\sin(t^{\top} \bm{\mathfrak{B}_0}^{\top}X_j)]\\
&& +M(t)^{\top} \Sigma^{-1} \frac{1}{\sqrt{n}}\sum_{j=1}^n \varepsilon_j g'(\theta_0, X_j) +R_n(t)\\
&=:& V_n^1(t)+R_n(t),
\end{eqnarray*}
where the remainder $R_n(t)$ satisfies
$ \int_{\mathbb{R}^{2d}}|R_n(t)|^2 \varphi(t)dt =o_p(1) $.

{\it Step 2}. In this step we will show that
$$ V_n^1(t) \longrightarrow V_{\infty}^1(t)   \quad {\rm in \ distribution \ in \ the \ space \ \mathbb{C}(\Gamma)}, $$
where $\Gamma$ is any compact subset in $\mathbb{R}^{2d}$, $\mathbb{C}(\Gamma)$ is the space of real-valued continuous functions on $\Gamma$, and $V_{\infty}^1(t)$ is a zero-mean Gaussian process with a covariance function $\psi(s, t)$. By the Continuous Mapping Theorem, we obtain
$$ \int_{\Gamma} |V_n^1(t)|^2 \varphi(t)dt \longrightarrow \int_{\Gamma} |V_{\infty}^1(t)|^2 \varphi(t)dt  \quad {\rm in \ distribution}.$$

First we show that $V_n^1(t)$ is asymptotically tight. Recall that
$$V_n^1(t)=\frac{1}{\sqrt{n}} \sum_{j=1}^{n} \varepsilon_j [\cos(t^{\top} \bm{\mathfrak{B}_0}^{\top}X_j)+\sin(t^{\top} \bm{\mathfrak{B}_0}^{\top}X_j)] - M(t)^{\top} \Sigma^{-1} \frac{1}{\sqrt{n}}\sum_{j=1}^n \varepsilon_j g'(\theta_0, X_j).
$$
For any $s, t \in \Gamma$, we have
\begin{eqnarray*}
&&  E|V_n^1(s)-V_n^1(t)|^2 \\
%&=& E|\frac{1}{\sqrt{n}} \sum_{j=1}^{n} \varepsilon_j [\cos(s^{\top} \bm{\mathfrak{B}_0}^{\top}X_j)-\cos(t^{\top} \bm{\mathfrak{B}_0}^{\top}X_j)]+ \\
%&&  \frac{1}{\sqrt{n}} \sum_{j=1}^{n} \varepsilon_j  [\sin(s^{\top} \bm{\mathfrak{B}_0}^{\top}X_j)-\sin(t^{\top} \bm{\mathfrak{B}_0}^{\top}X_j)]+ \\
%&&  [M(s)-M(t)]^{\top} \Sigma^{-1} \frac{1}{\sqrt{n}}\sum_{j=1}^n \varepsilon_j g'(\theta_0, X_j)X_j |^2 \\
&\leq& 3E|\frac{1}{\sqrt{n}} \sum_{j=1}^{n} \varepsilon_j [\cos(s^{\top} \bm{\mathfrak{B}_0}^{\top}X_j)-\cos(t^{\top} \bm{\mathfrak{B}_0}^{\top}X_j)] |^2 + \\
&& 3E|\frac{1}{\sqrt{n}} \sum_{j=1}^{n} \varepsilon_j  [\sin(s^{\top} \bm{\mathfrak{B}_0}^{\top}X_j)-\sin(t^{\top} \bm{\mathfrak{B}_0}^{\top}X_j)]|^2 +\\
&& 3E|[M(s)-M(t)]^{\top} \Sigma^{-1} \frac{1}{\sqrt{n}}\sum_{j=1}^n \varepsilon_j g'(\theta_0, X_j)|^2 \\
&=& 3E|\varepsilon [\cos(s^{\top} \bm{\mathfrak{B}_0}^{\top}X)-\cos(t^{\top} \bm{\mathfrak{B}_0}^{\top}X)]|^2 + \\
&&  3 E|\varepsilon [\sin(s^{\top} \bm{\mathfrak{B}_0}^{\top}X)-\sin(t^{\top} \bm{\mathfrak{B}_0}^{\top}X)]|^2 +\\
&& 3[M(s)-M(t)]^{\top} \Sigma^{-1}E[\varepsilon^2 g'(\theta_0, X) g'(\theta_0, X)^{\top}] \Sigma^{-1}[M(s)-M(t)] \\
&\leq& 6 E[\varepsilon^2 (s^{\top} \bm{\mathfrak{B}_0}^{\top}X-t^{\top} \bm{\mathfrak{B}_0}^{\top}X)^2] + C \|M(s)-M(t)\|^2 \\
&\leq& 6 \|s-t\|^2 E[\varepsilon^2 \| \bm{\mathfrak{B}_0}^{\top}X \|^2] + C\|\frac{\partial M}{\partial t}(\xi) (s-t)\|^2,
\end{eqnarray*}
where $\xi$ lies between $s$ and $t$. Note that
$$ \frac{\partial M}{\partial t}(\xi) = E\{ g'(\theta_0, X) [\cos(\xi^{\top} \bm{\mathfrak{B}_0}^{\top}X)-\sin(\xi^{\top} \bm{\mathfrak{B}_0}^{\top}X)] X^{\top} \bm{\mathfrak{B}_0} \}.$$
Thus we have
$$ \|\frac{\partial M}{\partial t}(\xi) (s-t)\|^2 \leq C  \|(s-t)\|^2. $$
By condition (B1), we have $E[\varepsilon^2 \| \bm{\mathfrak{B}_0}^{\top}X \|^2] \leq C$. Then it follows that
$$ E|V_n^1(s)-V_n^1(t)|^2 \leq C \|(s-t)\|^2. $$
By Theorem 12.3 of Billingsley (1968), we obtain that $V_n^1(t)$ is asymptotically tight.

Next we consider the convergence of finite-dimensional distributions of $V_n^1(t)$. Putting
\begin{eqnarray*}
Z_{nj}(t) = \frac{1}{\sqrt{n}} \varepsilon_j [\cos(t^{\top} \bm{\mathfrak{B}_0}^{\top}X_j)+\sin(t^{\top} \bm{\mathfrak{B}_0}^{\top}X_j)- M(t)^{\top} \Sigma^{-1} g'(\theta_0, X_j)].
\end{eqnarray*}
For any $t_1, \cdots, t_m \in \Gamma$, let $Z_{nj}=[Z_{nj}(t_1), \cdots, Z_{nj}(t_m)]^{\top}$. For any $\delta >0$, following the same arguments in Theorem 3.1 of Tan and Zhu (2019a), we obtain that
$$\sum_{j=1}^{n} E[\|Z_{nj} \|^2 I(\|Z_{nj}\| >\delta)]=O(\sqrt{p^3/n})=o(1).
$$
For the covariance matrix $\sum_{j=1}^{n} Cov(Z_{nj})$, we only need to calculate $\sum_{j=1}^{n} Cov\{Z_{nj}(s), Z_{nj}(t)\}$. It is easy to see that
\begin{eqnarray*}
&&  \sum_{j=1}^{n} Cov\{Z_{nj}(s), Z_{nj}(t)\} \\
&=& E\{ \varepsilon^2 [\cos((s-t)^{\top} \bm{\mathfrak{B}_0}^{\top}X)+\sin((s+t)^{\top} \bm{\mathfrak{B}_0}^{\top}X)]\}- \\
&& M(t)^{\top} \Sigma^{-1} E\{ \varepsilon^2 [\cos(s^{\top} \bm{\mathfrak{B}_0} ^{\top}X)+\sin(s^{\top} \bm{\mathfrak{B}_0}^{\top}X)] g'(\theta_0, X)\}- \\
&&  M(s)^{\top} \Sigma^{-1} E\{ \varepsilon^2 [\cos(t^{\top} \bm{\mathfrak{B}_0} ^{\top}X)+\sin(t^{\top} \bm{\mathfrak{B}_0}^{\top}X)] g'(\beta_0, X)\}+ \\
&& M(s)^{\top}\Sigma^{-1} E[\varepsilon^2 g'(\theta_0, X) g'(\theta_0, X)^{\top}] \Sigma^{-1} M(t).
\end{eqnarray*}
Note that $\sum_{j=1}^{n} Cov[Z_{nj}(s), Z_{nj}(t)] =K_n(s, t)$ and $K_n(s, t) \to K(s, t)$ point-wisely in $(s, t)$. Then it follows that $Z_{nj}$ satisfies conditions of Lindeberg-Feller Central limit theorem and then convergence of finite-dimensional distribution of $V_n^1(t)$ holds. Hence we obtain that
$$ V_n^1(t) \longrightarrow V_{\infty}^1(t)   \quad {\rm in \ distribution \ in \ the \ space \ \mathbb{C}(\Gamma)},$$
whence
$$ \int_{\Gamma} |V_n^1(t)|^2 \varphi(t)dt \longrightarrow \int_{\Gamma} |V_{\infty}^1(t)|^2 \varphi(t)dt  \quad {\rm in \ distribution}.$$

{\it Step 3}. In this step we will show that
$$ \int_{\mathbb{R}^{2d}} |V_n^1(t)|^2 \varphi(t)dt \longrightarrow \int_{\mathbb{R}^{2d}} |V_{\infty}^1(t)|^2 \varphi(t)dt  \quad {\rm in \ distribution}.$$
Fix $\varepsilon > 0$. As $\int_{\mathbb{R}^{2d}}K(t,t)\varphi(t)dt < \infty$, there exists a compact subset $\Gamma \subset \mathbb{R}^{2d}$ such that
$$\int_{\Gamma^c} K(t,t) \varphi(t)dt < \frac{\varepsilon^2}{2},  $$
where $\Gamma^c$ is the complementary set of $\Gamma$ in $\mathbb{R}^{2d}$.
Note that $K(t,t)=E|V_{\infty}^1(t)|^2$. Then it follows that
$$E\left(\int_{\Gamma^c} V_{\infty}^1(t)^2 \varphi(t)dt  \right)=\int_{\Gamma^c} E|V_{\infty}^1(t)|^2 \varphi(t)dt < \frac{\varepsilon^2}{2}.$$
Recall that
$$E|V_n^1(t)|^2=K_n(t,t) \quad {\rm and} \quad K_n(t,t) \longrightarrow K(t,t).  $$
Then it follows that
$$ E\left(\int_{\Gamma^c} V_n^1(t)^2 \varphi(t)dt \right) = \int_{\Gamma^c} K_n(t,t) \varphi(t)dt \longrightarrow \int_{\Gamma^c} K(t,t) \varphi(t)dt < \frac{\varepsilon^2}{2},$$
whence
$$ E\left(\int_{\Gamma^c} V_n^1(t)^2 \varphi(t)dt \right) < \varepsilon^2, \ {\rm for} \ n \ {\rm large \ enough}. $$
Putting
$$ V_{n1}^1 = \int_{\Gamma^c} V_n^1(t)^2 \varphi(t)dt, \quad V_{n2}^1 = \int_{\Gamma} V_n^1(t)^2 \varphi(t)dt,$$
$$ V_1=\int_{\Gamma^c} V_{\infty}^1(t)^2 \varphi(t)dt, \quad V_2=\int_{\Gamma} V_{\infty}^1(t)^2 \varphi(t)dt.$$
%Let $F_n$ and $F$ be the cumulative distribution function of $V_{n1}^1+V_{n2}^1$ and $V_1+V_2$ respectively.
Thus we have $E(V_1)< \varepsilon^2/2$, $E(V_{n1}^1) < \varepsilon^2$ for $n$ large enough, and $V_{n2}^1 \rightarrow V_2$ in distribution. Then it follows that
\begin{eqnarray*}
&& \mathbb{P}(V_1+V_2 \leq t-\varepsilon)-\varepsilon  \\
&\leq& \mathbb{P}(V_2 \leq t-\varepsilon)-\varepsilon \\
&=& \liminf_{n \to \infty} \mathbb{P}(V_{n2}^1 \leq t-\varepsilon)-\varepsilon \\
&=& \liminf_{n \to \infty} \{ \mathbb{P}(V_{n1}^1+V_{n2}^1 \leq t)+ \mathbb{P}(V_{n1}^1 \geq \varepsilon)\}-\varepsilon \\
&\leq&  \liminf_{n \to \infty}\mathbb{P}(V_{n1}^1+V_{n2}^1 \leq t) +\varepsilon-\varepsilon, \quad {\rm for} \ n \ {\rm large \ enough} \\
&\leq&  \limsup_{n \to \infty}\mathbb{P}(V_{n1}^1+V_{n2}^1 \leq t) \\
&\leq&  \limsup_{n \to \infty}\mathbb{P}(V_{n2}^1 \leq t) \\
&\leq&  \mathbb{P}(V_2 \leq t)  \\
&\leq&  \mathbb{P}(V_1+V_2 \leq t+\varepsilon)+\mathbb{P}(V_1 \geq \varepsilon) \\
&\leq&  \mathbb{P}(V_1+V_2 \leq t+\varepsilon)+ \frac{\varepsilon}{2}
\end{eqnarray*}
Let $\varepsilon$ tend to zero, we obtain that
$$ \lim_{n \to \infty}\mathbb{P}(V_{n1}^1+V_{n2}^1 \leq t)= \mathbb{P}(V_1+V_2 \leq t).$$
Consequently,
$$ \int_{\mathbb{R}^{2d}} |V_n^1(t)|^2 \varphi(t)dt \longrightarrow \int_{\mathbb{R}^{2d}} |V_{\infty}^1(t)|^2 \varphi(t)dt  \quad {\rm in \ distribution}.$$
In {\it Step} 1 we have shown
$$\hat{V}_n^1(t) = V_n^1(t)+R_n(t) \quad {\rm and} \quad \int_{\mathbb{R}^{2d}} |R_n(t)|^2 \varphi(t)dt=o_p(1).$$
Thus we have
$$\int_{\mathbb{R}^{2d}} |\hat{V}_n^1(t)|^2 \varphi(t)dt \longrightarrow \int_{\mathbb{R}^{2d}} |V_{\infty}^1(t)|^2 \varphi(t)dt  \quad {\rm in \ distribution}.$$
Hence we complete the proof.    \hfill$\Box$

\textbf{Proof of Theorem~\ref{theorem4.2}.} The proof follows the same line as that  in the Proposition 4 in Tan and Zhu (2019a). Thus we omit it here.  \hfill$\Box$

\textbf{Proof of Theorem~\ref{theorem4.3}.} We use the same notations as that in Theorem~\ref{theorem2.1}. Let $ F(\zeta)=\sum_{i=1}^{n}[Y_i-g(\theta_0+\zeta, X_i)]g'(\theta_0+\zeta, X_i)$ and $\psi_{\theta}(x,y)=[y-g(\theta, x)]g'(\theta, x)$. Then we have
\begin{eqnarray*}
\zeta^{\top}F(\zeta)=n\zeta^{\top}[P_n(\psi_{\theta_0+\zeta}-\psi_{\theta_0})-
P(\psi_{\theta_0+\zeta}-\psi_{\theta_0})]+ n \zeta^{\top} P_n\psi_{\theta_0} + n \zeta^{\top} P(\psi_{\theta_0+\zeta}-\psi_{\theta_0})
\end{eqnarray*}
Similar to the arguments in Theorem 2.1, we have
\begin{eqnarray*}
&&n \zeta^{\top} (P_n\psi_{\theta_0}-P\psi_{\theta_0})=\sqrt{np} \|\zeta\| O_p(1). \\
&&n \zeta^{\top} P(\psi_{\theta_0+\zeta}-\psi_{\theta_0})=n \zeta^{\top} \dot{P}\psi_{\theta_0} \zeta + \sqrt{p} n\|\zeta\|^3 O_p(1). \\
&&n\zeta^{\top}[P_n(\psi_{\theta_0+\zeta}-\psi_{\theta_0})-P(\psi_{\theta_0+\zeta}-
\psi_{\theta_0})]=\sqrt{np^2}\|\zeta\|^2 O_p(1).
\end{eqnarray*}
Note that $Y_n=g(\theta_0, X)+r_n G(\bm{B}^{\top}X)+\varepsilon$ and
$$ \dot{P}\psi_{\theta_0}=-E[g'(\theta_0, X)g'(\theta_0, X)^{\top}]+
E\{[Y_n-g(\theta_0, X)]g''(\theta_0, X)\} .$$
Then it follows that $\dot{P}\psi_{\theta_0}=-\Sigma + r_n E[G(\bm{B}^{\top}X)g''(\theta_0, X)] $.
Therefore,
\begin{eqnarray*}
n \zeta^{\top} P(\psi_{\theta_0+\zeta}-\psi_{\theta_0})&=&
-n \zeta^{\top} \Sigma \zeta + n r_n \zeta^{\top} E[G(\bm{B}^{\top}X)g''(\theta_0, X)] \zeta + \sqrt{p} n\|\zeta\|^3 O_p(1) \\
&=& -n \zeta^{\top} \Sigma \zeta + n r_n p \|\zeta\|^2 O_p(1)+ \sqrt{p} n\|\zeta\|^3 O_p(1).
\end{eqnarray*}
For $n \zeta^{\top} P\psi_{\theta_0}$, we have
\begin{eqnarray*}
n \zeta^{\top} P\psi_{\theta_0}&=&n \zeta^{\top} E\{[Y_n-g(\theta_0, X)]g'(\theta_0, X)\}\\
&=& n r_n \zeta^{\top} E[G(\bm{B}^{\top}X)g'(\theta_0, X)]=n r_n \sqrt{p} \|\zeta\| O_p(1)
\end{eqnarray*}
Let $\zeta=\sqrt{p} r_n U $. Then we obtain
\begin{eqnarray*}
\zeta^{\top} F(\zeta)&=&\sqrt{np^2}\|\zeta\|^2 O_p(1)+\sqrt{np} \|\zeta\| O_p(1) + n r_n \sqrt{p} \|\zeta\| O_p(1)\\
&& -n \zeta^{\top} \Sigma \zeta + n r_n p \|\zeta\|^2 O_p(1)+ \sqrt{p} n\|\zeta\|^3 O_p(1) \\
&\leq& \sqrt{n} p^2 r_n^2 \|U\|^2 O_p(1) + \sqrt{n} p r_n \|U\| O_p(1)+ n p r_n^2 \|U\| O_p(1) \\
&& -n p r_n^2 \lambda_{\min}(\Sigma) \|U\|^2 + n p^2 r_n^3 \| U\|^2 O_p(1) + np^2 r_n^3 \|U\|^3 O_p(1)\\
&=& n p r_n^2 \|U\| \{ \frac{p}{\sqrt{n}}\|U\| O_p(1)+ \frac{1}{\sqrt{n} r_n} O_p(1)+O_p(1)-\lambda_{\min}(\Sigma) \|U\| \\
&& + pr_n \|U\| O_p(1)+pr_n \|U\|^2 O_p(1) \}
\end{eqnarray*}
Note that $p^2/n \to 0$ and $pr_n \to 0$. Following the same arguments in Theorem 2.1, we have $\|\hat{\theta}-\theta_0\|=O_p(\sqrt{p}r_n)$.

Next we show the asymptotically linear representation of $\hat{\theta}-\theta_0$ under the alternatives $H_{1n}$. Following the same line of the proof in Theorem 2.2, if $n(p\log{n})^3r_n^4 \to 0$, we have
$$G_n[\alpha^{\top} (\psi_{\hat{\theta}}-\psi_{\theta_0})]=o_p(1), $$
uniformly in $\alpha \in \mathcal{S}^{q-1}$. Then it follows that
\begin{eqnarray*}
G_n(\alpha^{\top} \psi_{\hat{\theta}})
&=&\sqrt{n}(P_n \alpha^{\top} \psi_{\hat{\theta}}-P\alpha^{\top} \psi_{\hat{\theta}})
=-\sqrt{n}(P\alpha^{\top}\psi_{\hat{\theta}}-P\alpha^{\top} \psi_{\theta_0})- \sqrt{n} P\alpha^{\top} \psi_{\theta_0} \\
&=& G_n(\alpha^{\top} \psi_{\theta_0})+G_n[\alpha^{\top} (\psi_{\hat{\theta}}- \psi_{\theta_0})]=G_n(\alpha^{\top} \psi_{\theta_0})+o_p(1)
\end{eqnarray*}
Replacing $\alpha$ with $\Sigma^{-1} \alpha$, we have
\begin{eqnarray*}
\sqrt{n}P_n (\alpha^{\top}\Sigma^{-1}\psi_{\theta_0})=-\sqrt{n} \alpha^{\top} \Sigma^{-1}(P\psi_{\hat{\theta}}-P \psi_{\theta_0})+o_p(1).
\end{eqnarray*}
Applying the Taylor's expansion, we have
\begin{eqnarray*}
P\psi_{\hat{\theta}}-P\psi_{\theta_0}&=&\dot{P}\psi_{\theta}|_{\theta=\theta_0}
(\hat{\theta}-\theta_0)+(\hat{\theta}-\theta_0)^{\top} (\ddot{P}\psi_{\theta}|_{\theta=\tilde{\theta}_2})(\hat{\theta}-\theta_0)\\
&=& \{-\Sigma+ r_n E[G(\bm{B}^{\top}X)g''(\theta_0, X)] \} (\hat{\theta}-\theta_0)+(\hat{\theta}-\theta_0)^{\top} (\ddot{P}\psi_{\theta}|_{\theta=\tilde{\theta}_2})(\hat{\theta}-\theta_0) \\
&=& -\Sigma(\hat{\theta}-\theta_0) + r_n^2 p^{3/2} O(1).
\end{eqnarray*}
where $\Sigma=E[g'(\theta_0, X)g'(\theta_0, X)^{\top}]$ and $\tilde{\theta}_3$ lies between $\hat{\theta}$ and $\theta_0$. Therefore,
$$\sqrt{n}\alpha^{\top} \Sigma^{-1}(P\psi_{\hat{\theta}}-P\psi_{\theta_0})
= -\sqrt{n} \alpha^{\top} (\hat{\theta}-\theta_0) + \sqrt{np^3} r_n^2 O(1)
$$
Since
\begin{eqnarray*}
P_n (\psi_{\theta_0})&=&\frac{1}{n} \sum_{j=1}^{n}[Y_{nj}-g(\theta_0, X_j)]g'(\theta_0, X_j)\\
&=& \frac{1}{n} \sum_{j=1}^{n}\varepsilon_j g'(\theta_0, X_j)+ \frac{1}{n} \sum_{j=1}^{n} r_n G(\bm{B}^{\top}X_j) g'(\theta_0, X_j) \\
&=& \frac{1}{n} \sum_{j=1}^{n}\varepsilon_j g'(\theta_0, X_j)+ r_n E[G(\bm{B}^{\top}X) g'(\theta_0, X)] + \sqrt{p}r_n O_p(1),
\end{eqnarray*}
it follows that
$$\sqrt{n} \alpha^{\top} (\hat{\theta}-\theta_0)=\alpha^{\top} \Sigma^{-1} \frac{1}{\sqrt n} \sum_{j=1}^{n}\varepsilon_j g'(\theta_0, X_j)+ \sqrt{n}r_n \alpha^{\top} \Sigma^{-1} E[G(\bm{B}^{\top}X) g'(\theta_0, X)] +o_p(1). $$
Then we complete the proof.          \hfill$\Box$

\textbf{Proof of Theorem~\ref{theorem4.4}.} (1) Under the global alternative hypothesis $H_1$, we have $\mathbb{P}(\hat{q}=q) \to 1$. Thus we only work on the event $\{ \hat{q}=q \}$. Recall that
$$
\hat{V}_n^1(t) = \frac{1}{\sqrt{n}} \sum_{j=1}^{n} [Y_j-g(\hat{\theta}, X_j)] [\cos(t^{\top} \hat{\bm{\mathfrak B}}^{\top}X_j)+\sin(t^{\top} \hat{\bm{\mathfrak B}}^{\top}X_j)].
$$
Since $Y_j=G(\bm{B}^{\top}X_j)+\varepsilon_j$ under $H_1$, it follows that
\begin{eqnarray*}
\frac{1}{\sqrt{n}} \hat{V}_n^1(t) &=& \frac{1}{n} \sum_{j=1}^{n} \varepsilon_j [\cos(t^{\top} \bm{\hat{\mathfrak B}}^{\top}X_j)+\sin(t^{\top} \bm{\hat{\mathfrak B}}^{\top}X_j)] \\
&& -\frac{1}{n} \sum_{j=1}^{n} [g(\hat{\theta}, X_j)-g(\tilde{\theta}_0, X_j)] [\cos(t^{\top} \bm{\hat{\mathfrak B}}^{\top}X_j)+\sin(t^{\top} \bm{\bm{\hat{\mathfrak B}}}^{\top}X_j)] \\
&& +\frac{1}{n} \sum_{j=1}^{n} [G(\bm{B}^{\top}X_j)-g(\tilde{\theta}_0, X_j)] [\cos(t^{\top} \bm{\hat{\mathfrak B}}^{\top}X_j)+\sin(t^{\top} \bm{\hat{\mathfrak B}}^{\top}X_j)] \\
%&=& \frac{1}{\sqrt{n}} \sum_{j=1}^{n} \varepsilon_j [\cos(t^{\top} \bm{\mathfrak{B}_0}^{\top}X_j)+\sin(t^{\top} \bm{\mathfrak{B}_0}^{\top}X_j)] \\
%&&  +\frac{1}{\sqrt{n}} \sum_{j=1}^{n} \varepsilon_j [\cos(t^{\top} \bm{\hat{\mathfrak B}}^{\top}X_j)-\cos(t^{\top} \bm{\mathfrak{B}_0}^{\top}X_j) + \sin(t^{\top} \bm{\hat{\mathfrak B}}^{\top}X_j)-\sin(t^{\top} \bm{\mathfrak{B}_0}^{\top}X_j)]  \\
%&& -\frac{1}{\sqrt{n}} \sum_{j=1}^{n} [g(\hat{\theta}, X_j)-g(\theta_0, X_j)] [\cos(t^{\top} \bm{\mathfrak{B}_0}^{\top}X_j)+\sin(t^{\top} \bm{\mathfrak{B}_0}^{\top}X_j)]  \\
%&& -\frac{1}{\sqrt{n}} \sum_{j=1}^{n} [g(\hat{\theta}, X_j)-g(\theta_0, X_j)] [\cos(t^{\top} \bm{\hat{\mathfrak B}}^{\top}X_j)-\cos(t^{\top} \bm{\mathfrak{B}_0}^{\top}X_j) + \sin(t^{\top} \bm{\hat{\mathfrak B}}^{\top}X_j)-\sin(t^{\top} \bm{\mathfrak{B}_0}^{\top}X_j)] \\
&=:& T_{n1}(t)-T_{n2}(t)+T_{n3}(t),
\end{eqnarray*}
Following the same line of the arguments in the Theorem 3.1, we have
$$\int_{\mathbb{R}^{d+q}} |T_{n1}(t)|^2 \varphi(t)dt=o_p(1) \quad {\rm and} \quad \int_{\mathbb{R}^{d+q}} |T_{n2}(t)|^2 \varphi(t)dt=o_p(1). $$
For the third term $T_{n3}(t)$ in $(1/\sqrt{n})\hat{V}_n^1(t)$, we have
\begin{eqnarray*}
T_{n3}(t)&=&\frac{1}{n} \sum_{j=1}^{n} [G(\bm{B}^{\top}X_j)-g(\tilde{\theta}_0, X_j)] [\cos(t^{\top} \bm{\hat{\mathfrak B}}^{\top}X_j)+\sin(t^{\top} \bm{\hat{\mathfrak B}}^{\top}X_j)] \\
&=& \frac{1}{n} \sum_{j=1}^{n} [G(\bm{B}^{\top}X_j)-g(\tilde{\theta}_0, X_j)] [\cos(t^{\top} \bm{\tilde{\mathfrak B}}^{\top}X_j)+\sin(t^{\top} \bm{\tilde{\mathfrak B}}^{\top}X_j)] \\
&&-\frac{1}{n} \sum_{j=1}^{n} [G(\bm{B}^{\top}X_j)-g(\tilde{\theta}_0, X_j)] [\cos(t^{\top} \bm{\hat{\mathfrak B}}^{\top}X_j)-\cos(t^{\top} \bm{\tilde{\mathfrak B}}^{\top}X_j)]\\
&& -\frac{1}{n} \sum_{j=1}^{n} [G(\bm{B}^{\top}X_j)-g(\tilde{\theta}_0, X_j)] [\sin(t^{\top} \bm{\hat{\mathfrak B}}^{\top}X_j)-\sin(t^{\top} \bm{\tilde{\mathfrak B}}^{\top}X_j)] \\
&=:& T_{n31}(t)-T_{n32}(t)-T_{n33}(t)
\end{eqnarray*}
where $\bm{\tilde{\mathfrak B}}=(\bm{\tilde{\beta}_0}, \bm{B})$.
Similar to the argument for the term $V_{n21}(t)$ in the proof of Theorem 4.1, we have
$$ \int_{\mathbb{R}^{d+q}} |T_{n3k}(t)|^2 \varphi(t)dt = O_p(\frac{p^2}{n})=o_p(1), \quad k=2,3. $$
Let $L_n(t)=E[T_{n31}(t)]$. Then it follows that
\begin{eqnarray*}
&& E\left( \int_{\mathbb{R}^{d+q}} |T_{n31}(t)-L_n(t)|^2 \varphi(t)dt \right) \\
&=& \frac{1}{n} \int_{\mathbb{R}^{d+q}}  Var\{[G(\bm{B}^{\top}X_j)-g(\tilde{\theta}_0, X_j)] [\cos(t^{\top} \bm{\hat{\mathfrak B}}^{\top}X_j)+\sin(t^{\top} \bm{\hat{\mathfrak B}}^{\top}X_j)] \} \varphi(t)dt \\
&\leq& \frac{1}{n} E[G(\bm{B}^{\top}X_j)-g(\tilde{\theta}_0, X_j)]^2 \leq C \frac{1}{n}.
\end{eqnarray*}
Then we have $ \int_{\mathbb{R}^{d+q}} |T_{n31}(t)-L_n(t)|^2 \varphi(t)dt =o_p(1)$.
Altogether we obtain that
$$ \frac{1}{\sqrt{n}} \hat{V}_n^1(t)= L_n(t) + R_n^1(t), $$
where $R_n^1(t)$ is a remainder satisfying $\int_{\mathbb{R}^{d+q}} |R_n^1(t)|^2 \varphi(t)dt=o_p(1) $.
Therefore, we have in probability
$$\int_{\mathbb{R}^{d+q}} |\frac{1}{\sqrt{n}} \hat{V}_n^1(t)|^2 \varphi(t)dt=
\int_{\mathbb{R}^{d+q}} |L_n(t)|^2 \varphi(t)dt +o_p(1) \longrightarrow C_1 >0.  $$

(2) We only work on the event $\{\hat{q}=d\}$ as $\mathbb{P}(\hat{q}=d) \to 1$ under the alternatives $H_{n1}$ with $p r_n \to 0$.
Since $Y_{nj}=g(\hat{\theta}, X_j)+r_n G(\bm{B}^{\top}X_j) +\varepsilon_j$, we decompose $\hat{V}_n^1(t)$ as follows,
\begin{eqnarray*}
\hat{V}_n^1(t) &=& \frac{1}{\sqrt{n}} \sum_{j=1}^{n} [Y_j-g(\hat{\theta}, X_j)] [\cos(t^{\top} \bm{\hat{\mathfrak B}}^{\top}X_j)+\sin(t^{\top} \bm{\hat{\mathfrak B}}^{\top}X_j)]\\
&=& \frac{1}{\sqrt{n}} \sum_{j=1}^{n} \varepsilon_j [\cos(t^{\top} \bm{\hat{\mathfrak B}}^{\top}X_j)+\sin(t^{\top} \bm{\hat{\mathfrak B}}^{\top}X_j)] \\
&& -\frac{1}{\sqrt{n}} \sum_{j=1}^{n} [g(\hat{\theta}, X_j)-g(\theta_0, X_j)] [\cos(t^{\top} \bm{\hat{\mathfrak B}}^{\top}X_j)+\sin(t^{\top} \bm{\bm{\hat{\mathfrak B}}}^{\top}X_j)] \\
&& +\sqrt{n}r_n \frac{1}{n} \sum_{j=1}^{n} G(\bm{B}^{\top}X_j) [\cos(t^{\top} \bm{\hat{\mathfrak B}}^{\top}X_j)+\sin(t^{\top} \bm{\hat{\mathfrak B}}^{\top}X_j)] \\
&=:& W_{n1}(t)-W_{n2}(t)+W_{n3}(t)
\end{eqnarray*}
For the first term $W_{n1}(t)$ in $\hat{V}_n^1(t)$, we have
\begin{eqnarray*}
W_{n1}(t)&=& \frac{1}{\sqrt{n}} \sum_{j=1}^{n} \varepsilon_j [\cos(t^{\top} \bm{{\mathfrak B}_0}^{\top}X_j)+\sin(t^{\top} \bm{{\mathfrak B}_0}^{\top}X_j)] \\
&&+ \frac{1}{\sqrt{n}} \sum_{j=1}^{n} \varepsilon_j [\cos(t^{\top} \bm{\hat{\mathfrak B}}^{\top}X_j)-\cos(t^{\top} \bm{{\mathfrak B}_0}^{\top}X_j) +\sin(t^{\top} \bm{\hat{\mathfrak B}}^{\top}X_j)-\sin(t^{\top} \bm{{\mathfrak B}_0}^{\top}X_j)] \\
&=:& W_{n11}(t) + W_{n12}(t)
\end{eqnarray*}
Then it is easy to see that
$$ E \left( \int_{\mathbb{R}^{2d}} |W_{n11}(t)|^2 \varphi(t)dt \right) \leq 2 E (\varepsilon^2) \int_{\mathbb{R}^{2d}} \varphi(t)dt. $$
Thus we obtain that $ \int_{\mathbb{R}^{2d}} |W_{n11}(t)|^2 \varphi(t)dt =O_p(1) $. By the same arguments for the term $V_{n21}(t)$ in Theorem 4.1, we have
$$ \int_{\mathbb{R}^{2d}} |W_{n12}(t)|^2 \varphi(t)dt = O_p(pr_n)^2+ O_p(pr_n)^4 + O_p(n(pr_n)^6). $$
Note that $n(p\log{n})^3 r_n^4 \to 0$. Then it follows that
$$\int_{\mathbb{R}^{2d}} |\frac{1}{\sqrt{n}r_n} W_{n1}(t)|^2 \varphi(t)dt =o_p(1).  $$

For the second term $W_{n2}(t)$ in $\hat{V}_n^1(t)$, we have
\begin{eqnarray*}
W_{n2}(t)&=& \frac{1}{\sqrt{n}} \sum_{j=1}^{n} [g(\hat{\theta}, X_j)-g(\theta_0, X_j)] [\cos(t^{\top} \bm{{\mathfrak B}_0}^{\top}X_j)+\sin(t^{\top} \bm{{\mathfrak B}_0}^{\top}X_j)] \\
&& + \frac{1}{\sqrt{n}} \sum_{j=1}^{n} [g(\hat{\theta}, X_j)-g(\theta_0, X_j)] [\cos(t^{\top} \bm{\hat{\mathfrak B}}^{\top}X_j)-\cos(t^{\top} \bm{{\mathfrak B}_0}^{\top}X_j)] \\
&& + \frac{1}{\sqrt{n}} \sum_{j=1}^{n} [g(\hat{\theta}, X_j)-g(\theta_0, X_j)] [\sin(t^{\top} \bm{\hat{\mathfrak B}}^{\top}X_j)-\sin(t^{\top} \bm{{\mathfrak B}_0}^{\top}X_j)]\\
&=:& W_{n21}(t) + W_{n22}(t)+ W_{n23}(t)
\end{eqnarray*}
For the term $W_{n21}(t)$, we have
\begin{eqnarray*}
W_{n21}(t)&=& \frac{1}{\sqrt{n}} \sum_{j=1}^{n} (\hat{\theta}-\theta_0)^{\top} g'(\theta_0, X_j)[\cos(t^{\top} \bm{\mathfrak{B}_0}^{\top}X_j)+\sin(t^{\top} \bm{\mathfrak{B}_0}^{\top}X_j)] \\
&& + \frac{1}{2\sqrt{n}} \sum_{j=1}^{n} (\hat{\theta}-\theta_0)^{\top} g''(\theta_0, X_j)(\hat{\theta}-\theta_0)[\cos(t^{\top} \bm{\mathfrak{B}_0}^{\top}X_j)+\sin(t^{\top} \bm{\mathfrak{B}_0}^{\top}X_j)] \\
&& + \frac{1}{6\sqrt{n}} \sum_{j=1}^{n} (\hat{\theta}-\theta_0)^{\top}[(\hat{\theta}-\theta_0)^{\top} g'''(\theta_3, X_j) (\hat{\theta}-\theta_0)] [\cos(t^{\top} \bm{\mathfrak{B}_0}^{\top}X_j)+ \sin(t^{\top} \bm{\mathfrak{B}_0}^{\top}X_j)]\\
&=:& W_{n211}(t)+\frac{1}{2}W_{n212}(t)+\frac{1}{6} W_{n213}(t),
\end{eqnarray*}
where $\theta_3$ lies between $\hat{\theta}$ and $\theta_0$. Similar to the arguments for the term $V_{n3}(t)$ in the proof of Theorem 4.1, we obtain
\begin{eqnarray*}
\int_{\mathbb{R}^{2d}} |W_{n212}(t)|^2 \varphi(t)dt &=& O_p(np^2 r_n^4) + O_p(pr_n)^4 \\
\int_{\mathbb{R}^{2d}} |W_{n213}(t)|^2 \varphi(t)dt &=& O_p(np^6 r_n^6)
\end{eqnarray*}
For the term $W_{n211}(t)$, recall that $M(t)=\{g'(\theta_0, X)[\cos(t^{\top} \bm{\mathfrak{B}_0}^{\top}X)+\sin(t^{\top} \bm{\mathfrak{B}_0}^{\top}X)]\}$. By the same arguments in $V_{n31}(t)$ in the proof of Theorem 4.1, we have
$$W_{n211}(t)=\sqrt{n} (\hat{\theta}-\theta_0)^{\top}M(t) + R_{n211}(t) \quad {\rm and} \quad \int_{\mathbb{R}^{2d}} |R_{n211}(t)|^2 \varphi(t)dt=O_p(pr_n)^2.   $$
By Theorem 4.3, it follows that
\begin{eqnarray*}
\sqrt{n}M(t)^{\top} (\hat{\theta}-\theta_0)
&=& \frac{1}{\sqrt{n}} \sum_{j=1}^{n} \varepsilon_jM(t)^{\top}\Sigma^{-1} g'(\theta_0, X_j) \\
&& +\sqrt{n} r_n M(t)^{\top}\Sigma^{-1} E[g'(\theta_0, X)G(\bm{B}^{\top}X)] + o_p(1),
\end{eqnarray*}
where $o_p(1)$ is uniformly in $t$.  Note that $\int_{\mathbb{R}^{2d}} E[\varepsilon M(t)^{\top}\Sigma^{-1} g'(\theta_0, X)]^2 \varphi(t)dt =O(1)$.  Then it follows that
$$ \int_{\mathbb{R}^{2d}} |\frac{1}{\sqrt{n} r_n} W_{n211}(t)|^2 \varphi(t)dt = \int_{\mathbb{R}^{2d}} |M(t)^{\top}\Sigma^{-1} E[g'(\theta_0, X)G(\bm{B}^{\top}X)]|^2 \varphi(t)dt +o_p(1). $$
%By the same arguments for $W_{n212}(t)$ and $W_{n213}$, we have
%$$ \int_{\mathbb{R}^{2d}} |W_{n2k}(t)|^2 \varphi(t)dt = O_p(np^2 r_n^4) + O_p(pr_n)^4 +O_p(n p^6 r_n^6), \ {\rm for } \ k=2,3.  $$
Therefore we obtain
$$ \frac{1}{\sqrt{n} r_n} W_{n2}(t) = L_{n1}(t) +R_{n2}(t) \quad {\rm and } \quad \int_{\mathbb{R}^{2d}} |R_{n2}(t)|^2 \varphi(t)dt =o_p(1). $$
where $L_{n1}(t)=M(t)^{\top}\Sigma^{-1} E[g'(\theta_0, X)G(\bm{B}^{\top}X)]$.

For the third term $W_{n3}(t)$ in $\hat{V}_n^1(t)$, we decompose it as follows,
\begin{eqnarray*}
\frac{1}{\sqrt{n} r_n} W_{n3}(t) &=& \frac{1}{n} \sum_{j=1}^{n} G(\bm{B}^{\top}X_j) [\cos(t^{\top} \bm{\mathfrak{B}_0}^{\top}X_j)+\sin(t^{\top} \bm{\mathfrak{B}_0}^{\top}X_j)] \\
&& + \frac{1}{n} \sum_{j=1}^{n} G(\bm{B}^{\top}X_j) [\cos(t^{\top} \bm{\hat{\mathfrak B}}^{\top}X_j)-\cos(t^{\top} \bm{\mathfrak{B}_0}^{\top}X_j)]\\
&& + \frac{1}{n} \sum_{j=1}^{n} G(\bm{B}^{\top}X_j) [\sin(t^{\top} \bm{\hat{\mathfrak B}}^{\top}X_j)-\sin(t^{\top} \bm{\mathfrak{B}_0}^{\top}X_j)] \\
&=:& \frac{1}{\sqrt{n} r_n} W_{n31}(t)+\frac{1}{\sqrt{n} r_n} W_{n32}(t)+\frac{1}{\sqrt{n} r_n} W_{n33}(t).
\end{eqnarray*}
Let $L_{n2}(t)=E\{G(\bm{B}^{\top}X) [\cos(t^{\top} \bm{\mathfrak{B}_0}^{\top}X)+ \sin(t^{\top} \bm{\mathfrak{B}_0}^{\top}X)]\}$. It is easy to see that
$$ \int_{\mathbb{R}^{2d}} |\frac{1}{\sqrt{n} r_n} W_{n31}(t) - L_{n2}(t)|^2 \varphi(t)dt = O_p(\frac{1}{n}).  $$
For the terms $W_{n32}(t)$ and $W_{n33}(t)$, by the Taylor expansion, we obtain
$$ \int_{\mathbb{R}^{2d}} |\frac{1}{\sqrt{n} r_n} W_{n3k}(t)|^2 \varphi(t)dt =O_p(pr_n)^2+O_p(pr_n)^4, \ {\rm for} \ k=2,3. $$
Then it follows that
$$ \frac{1}{\sqrt{n} r_n} W_{n3}(t) = L_{n2}(t) +R_{n3}(t) \quad {\rm and } \quad \int_{\mathbb{R}^{2d}} |R_{n3}(t)|^2 \varphi(t)dt =o_p(1). $$
Altogether we obtain that
$$\int_{\mathbb{R}^{2d}} |\frac{1}{\sqrt{n} r_n} \hat{V}_n^1(t)|^2 \varphi(t)dt = \int_{\mathbb{R}^{2d}} |L_{n1}(t)+L_{n2}(t)|^2 \varphi(t)dt +o_p(1) \longrightarrow C_2 >0.   $$

(3) Following the same line as that in (2) with $r_n=1/\sqrt{n}$, we have
\begin{eqnarray*}
\hat{V}_n^1(t) &=&  \frac{1}{\sqrt{n}} \sum_{j=1}^{n} \varepsilon_j [\cos(t^{\top} \bm{{\mathfrak B}_0}^{\top}X_j)+\sin(t^{\top} \bm{{\mathfrak B}_0}^{\top} X_j)-M(t)^{\top}\Sigma^{-1} g'(\theta_0, X_j)]\\
&&-L_{n1}(t)+L_{n2}(t) + R_n^1(t),
\end{eqnarray*}
where $R_n^1(t)$ is a remainder with $\int_{\mathbb{R}^{2d}} |R_n^1(t)|^2 \varphi(t)dt = o_p(1)$ and
\begin{eqnarray*}
L_{n1}(t)&=& M(t)^{\top}\Sigma^{-1} E[g'(\theta_0, X)G(\bm{B}^{\top}X)]  \\
L_{n2}(t)&=& E\{G(\bm{B}^{\top}X) [\cos(t^{\top} \bm{\mathfrak{B}_0}^{\top}X)+ \sin(t^{\top} \bm{\mathfrak{B}_0}^{\top}X)] \}.
\end{eqnarray*}
The rest of proof is the same as that in the proof of step 2 and 3 of Theorem 4.1. Therefore we obtain
$$ \int_{\mathbb{R}^{2d}} |\hat{V}_n^1(t)|^2 \varphi(t)dt \rightarrow \int_{\mathbb{R}^{2d}} |V_{\infty}^1(t)+ L_1(t)-L_2(t)|^2 \varphi(t)dt,  $$
where $V_{\infty}^1(t)$ is a zero-mean Gaussian process given in Theorem 4.1 and $L_1(t)$ and $L_2(t)$ are the uniformly limits of $L_{n1}(t)$ and $L_{n2}(t)$ respectively.      \hfill$\Box$

\textbf{Proof of Theorem~\ref{theorem5.1}.} We only prove the theorem under the null hypothesis. The arguments under the alternatives are similar and thus we omit it here. Let $ P^*$ be the probability measure induced by the wild bootstrap resampling conditional on the original sample $\{(X_i, Y_i): i=1, \cdots, n\}$ and let  $E^*$ be the expectation under $P^*$. Recall that
$$\hat{V}_n^{1*}(t) = \frac{1}{\sqrt{n}} \sum_{j=1}^{n} [Y_j^*-g(\hat{\theta}^*, X_j)] [\cos(t^{\top} \hat{\bm{\mathfrak B}}^{\top}X_j)+\sin(t^{\top} \hat{\bm{\mathfrak B}}^{\top}X_j)].  $$
Decompose $\hat{V}_n^{1*}(t)$ as follows,
\begin{eqnarray*}
\hat{V}_n^{1*}(t) &=& \frac{1}{\sqrt{n}} \sum_{j=1}^{n} [\varepsilon_j^*+ g(\hat{\theta}, X_j) -g(\hat{\theta}^*, X_j)] [\cos(t^{\top} \hat{\bm{\mathfrak B}}^{\top}X_j)+\sin(t^{\top} \hat{\bm{\mathfrak B}}^{\top}X_j)] \\
&=& \frac{1}{\sqrt{n}} \sum_{j=1}^{n} \varepsilon_j^* [\cos(t^{\top} \hat{\bm{\mathfrak B}}^{\top}X_j)+\sin(t^{\top} \hat{\bm{\mathfrak B}}^{\top}X_j)] \\
&& -\frac{1}{\sqrt{n}} \sum_{j=1}^{n} [g(\hat{\theta}^*, X_j)-g(\hat{\theta}, X_j)] [\cos(t^{\top} \hat{\bm{\mathfrak B}}^{\top}X_j)+\sin(t^{\top} \hat{\bm{\mathfrak B}}^{\top}X_j)] \\
&=:& \hat{V}_{n1}^{1*}(t)-\hat{V}_{n2}^{1*}(t)
\end{eqnarray*}
For the term $\hat{V}_{n1}^{1*}(t)$, we can decompose it as follows
$$  \hat{V}_{n1}^{1*}(t)=\frac{1}{\sqrt{n}} \sum_{j=1}^{n} \varepsilon_j^* [\cos(t^{\top} \bm{\mathfrak B}_0^{\top}X_j)+\sin(t^{\top} \bm{\mathfrak B}_0^{\top}X_j)] + R_{n1}^{1*}(t). $$
Note that $E^* (\varepsilon_j^*)=0$ and $E^*(\varepsilon_j^*)^2=\hat{\varepsilon}_j^2$. Then it follows that
\begin{eqnarray*}
&& E^*\left(\int_{\mathbb{R}^{2d}}|R_{n1}^{1*}(t)|^2 \varphi(t)dt\right) \\
&=& \frac{1}{n} \sum_{j=1}^{n} \hat{\varepsilon}_j^2 [\cos(t^{\top} \hat{\bm{\mathfrak B}}^{\top}X_j)-\cos(t^{\top} \bm{\mathfrak B}_0^{\top}X_j)+ \sin(t^{\top} \hat{\bm{\mathfrak B}}^{\top}X_j) -\sin(t^{\top} \bm{\mathfrak B}_0^{\top}X_j)]^2 \\
&=& o_p(1).
\end{eqnarray*}
Thus we obtain that $\int_{\mathbb{R}^{2d}}|R_{n1}^{1*}(t)|^2\varphi(t)dt=o_{p^*}(1)$.

For the term $\hat{V}_{n2}^{1*}(t)$, applying the Taylor expansion, we have
\begin{eqnarray*}
\hat{V}_{n2}^{1*}(t) &=& \frac{1}{\sqrt{n}} \sum_{j=1}^{n} (\hat{\theta}^*-\hat{\theta})^{\top} g'(\hat{\theta}, X_j) [\cos(t^{\top} \hat{\bm{\mathfrak B}}^{\top}X_j)+\sin(t^{\top} \hat{\bm{\mathfrak B}}^{\top}X_j)] \\
&& + \frac{1}{2\sqrt{n}} \sum_{j=1}^{n} (\hat{\theta}^*-\hat{\theta})^{\top} g''(\hat{\theta}, X_j)(\hat{\theta}^*-\hat{\theta}) [\cos(t^{\top} \hat{\bm{\mathfrak B}}^{\top}X_j)+\sin(t^{\top} \hat{\bm{\mathfrak B}}^{\top}X_j)] \\
&& + \frac{1}{6\sqrt{n}} \sum_{j=1}^{n} (\hat{\theta}^*-\hat{\theta})^{\top} [(\hat{\theta}^*-\hat{\theta})^{\top} g'''(\hat{\theta}_2, X_j) (\hat{\theta}^*-\hat{\theta})] [\cos(t^{\top} \hat{\bm{\mathfrak B}}^{\top}X_j)+\sin(t^{\top} \hat{\bm{\mathfrak B}}^{\top}X_j)]\\
&=:& \hat{V}_{n21}^{1*}(t)+\hat{V}_{n22}^{1*}(t)+\hat{V}_{n23}^{1*}(t)
\end{eqnarray*}
where $\hat{\theta}_2$ lies between $\hat{\theta}^*$ and $\hat{\theta}$. Similar to the arguments for $V_{n32}(t)$ and $V_{n33}(t)$, we have in probability
$$ \int_{\mathbb{R}^{2d}} |\hat{V}_{n2k}^{1*}(t)|^2 \varphi(t)dt =O_{p^*} (\frac{p^3}{n}) \quad {\rm for} \ k=2,3. $$
For the term $\hat{V}_{n21}^{1*}(t)$, we obtain
\begin{eqnarray*}
\hat{V}_{n21}^{1*}(t)&=&\sqrt{n}(\hat{\theta}^*-\hat{\theta})^{\top} \frac{1}{n} \sum_{j=1}^{n} g'(\theta_0, X_j) [\cos(t^{\top} \bm{\mathfrak B}_0^{\top}X_j) +\sin(t^{\top} \bm{\mathfrak B}_0^{\top}X_j)] + R_{n211}^{1*}(t) \\
&=& \sqrt{n}(\hat{\theta}^*-\hat{\theta})^{\top} M(t)+R_{n211}^{1*}(t) +R_{n212}^{1*}(t),
\end{eqnarray*}
where $M(t)=E\{g'(\theta_0, X_j) [\cos(t^{\top} \bm{\mathfrak B}_0^{\top}X_j) +\sin(t^{\top} \bm{\mathfrak B}_0^{\top}X_j)]\}$ and the remainders $R_{n211}^{1*}(t)$ and $R_{n212}^{1*}(t)$ satisfy
$$\int_{\mathbb{R}^{2d}}|R_{n211}^{1*}(t)|^2 \varphi(t)dt=O_{p^*}(p^3/n) \quad {\rm and} \quad \int_{\mathbb{R}^{2d}}|R_{n212}^{1*}(t)|^2 \varphi(t)dt=O_{p^*}(p^2/n) . $$
Now we need the asymptotically linear expansion of $(\hat{\theta}^*-\hat{\theta})$. Following the arguments in Theorem 2.1 and 2.2, we obtain that
$$\sqrt{n} \alpha^{\top}(\hat{\theta}^*-\hat{\theta})= \frac{1}{\sqrt{n}} \sum_{j=1}^{n} \varepsilon_j^* \alpha^{\top} \Sigma_n(\hat{\theta})^{-1} g'(\hat{\theta}, X_j) + o_{p^*}(1), $$
where $\Sigma_n(\hat{\theta})=(1/n) \sum_{j=1}^{n} g'(\hat{\theta}, X_j)g'(\hat{\theta}, X_j)^{\top} $.
Then it follows that
\begin{eqnarray*}
\sqrt{n}M(t)^{\top}(\hat{\theta}^*-\hat{\theta})= \frac{1}{\sqrt{n}} \sum_{j=1}^{n} \varepsilon_j^* M(t)^{\top} \Sigma_n(\hat{\theta})^{-1} g'(\hat{\theta}, X_j) + o_{p^*}(1).
\end{eqnarray*}
Consequently we obtain
$$ \hat{V}_{n2}^{1*}(t)=\frac{1}{\sqrt{n}} \sum_{j=1}^{n} \varepsilon_j^* M(t)^{\top} \Sigma^{-1} g'(\theta_0, X_j) + R_{n2}^{1*}(t), $$
where $\Sigma=E[g'(\theta_0, X)g'(\theta_0, X)^{\top}]$ and the remainder $R_{n2}^{1*}(t)$ satisfies
$$\int_{\mathbb{R}^{2d}}|R_{n2}^{1*}(t)|^2 \varphi(t)dt=o_{p^*}(1). $$
Altogether we obtain that
\begin{eqnarray*}
\hat{V}_{n}^{1*}(t)&=&\frac{1}{\sqrt{n}} \sum_{j=1}^{n} \varepsilon_j^* [\cos(t^{\top} \bm{\mathfrak B}_0^{\top}X_j)+\sin(t^{\top} \bm{\mathfrak B}_0^{\top}X_j)-M(t)^{\top} \Sigma^{-1} g'(\theta_0, X_j)] + R_{n}^{1*}(t)\\
&=:& V_{n}^{1*}(t) + R_{n}^{1*}(t).
\end{eqnarray*}
Here $R_{n}^{1*}(t)=R_{n1}^{1*}(t)+R_{n2}^{1*}(t)$ and then $\int_{\mathbb{R}^{2d}}|R_{n}^{1*}(t)|^2 \varphi(t)dt=o_{p^*}(1)$.

It remains to derive the limiting distribution of $V_{n}^{1*}(t)$ under the probability measure $P^*$. Thus we need to deal with the covariance structure, convergence of the finite dimensional distributions of $V_{n}^{1*}(t)$ and the tightness under the induced bootstrap probability measure $P^*$.
Similar to the arguments in Lemma A.1-A.3 in Stute, Gonz$\acute{a}$lez Manteiga and Presedo Quindimil (1998), we have in probability,
$$ \hat{V}_n^{1*}(t) \longrightarrow V_{\infty}^{1*}(t) \quad {\rm in \ distribution \ in \ \mathbb{C}(\Gamma)}, $$
where $\Gamma$ is a compact subset of $\mathbb{R}^{2d}$ and $V_{\infty}^{1*}$ have the same distribution as the Gaussian process  $V_{\infty}^1$ given in Theorem~{\ref{theorem4.1}}. The rest of proof follows the same line as that in the arguments of Theorem~{\ref{theorem4.1}}, with the expectation $E$ replaced by the bootstrap expectation $E^*$. Hence we complete the proof.  \hfill$\Box$

\newpage
\begin{table}[ht!]\caption{Empirical sizes and powers of the tests for $H_0$ vs. $H_{11}$.}
\centering
{\small\scriptsize\hspace{12.5cm}
\renewcommand{\arraystretch}{1}\tabcolsep 0.5cm
\begin{tabular}{*{20}{c}}
\hline
&\multicolumn{1}{c}{a}&\multicolumn{1}{c}{n=100}&\multicolumn{1}{c}{n=200}&\multicolumn{1}{c}{n=400}&\multicolumn{1}{c}{n=600}\\
&&\multicolumn{1}{c}{p=8}& \multicolumn{1}{c}{p=12}  & \multicolumn{1}{c}{p=17}&\multicolumn{1}{c}{p=20}\\
\hline
$AICM_n$   &0.00               &0.0580  &0.0540  &0.0410  &0.0560\\
              &0.25               &0.3670  &0.6260  &0.8710  &0.9770\\
\hline
$ACM_n$       &0.00               &0.0440  &0.0495  &0.0570  &0.0510\\
              &0.25               &0.2720  &0.4960  &0.8125  &0.9335\\
\hline
$PCVM_n$      &0.00               &0.0660  &0.0530  &0.0540  &0.0650\\
              &0.25               &0.3080  &0.5140  &0.8330  &0.9350\\
\hline
$T_n^{SZ}$    &0.00               &0.0475  &0.0450  &0.0450  &0.0450\\
              &0.25               &0.2665  &0.4945  &0.7970  &0.9315\\
\hline
$ICM_n$       &0.00               &0.0020  &0.0000  &0.0000  &0.0000\\
              &0.25               &0.0200  &0.0000  &0.0000  &0.0000\\
\hline
$T_n^{GWZ}$   &0.00               &0.0345  &0.0580  &0.0490  &0.0605\\
              &0.25               &0.2550  &0.4435  &0.7570  &0.9100\\
\hline
$T_n^{ZH}$    &0.00               &0.0305  &0.0210  &0.0250  &0.0285\\
              &0.25               &0.0490  &0.0370  &0.0225  &0.0255\\
\hline
\end{tabular}}
\end{table}

\begin{table}[ht!]\caption{Empirical sizes and powers of  the tests for $H_0$ vs. $H_{12}$.}
\centering
{\small\scriptsize\hspace{12.5cm}
\renewcommand{\arraystretch}{1}\tabcolsep 0.5cm
\begin{tabular}{*{20}{c}}
\hline
&\multicolumn{1}{c}{a}&\multicolumn{1}{c}{n=100}&\multicolumn{1}{c}{n=200}&\multicolumn{1}{c}{n=400}&\multicolumn{1}{c}{n=600}\\
&&\multicolumn{1}{c}{p=8}& \multicolumn{1}{c}{p=12}  & \multicolumn{1}{c}{p=17}&\multicolumn{1}{c}{p=20}\\
\hline
$AICM_n$   &0.00         &0.0530 &0.0630 &0.0550 &0.0510\\
              &0.50         &0.5830 &0.8970 &0.9990 &1.0000\\
\hline
$ACM_n$       &0.00         &0.0435 &0.0595 &0.0495 &0.0480\\
              &0.50         &0.1465 &0.3045 &0.7255 &0.9580\\
\hline
$PCvM_n$      &0.00         &0.0790 &0.0470 &0.0500 &0.0790\\
              &0.50         &0.1650 &0.2330 &0.3860 &0.5060\\
\hline
$T_n^{SZ}$    &0.00         &0.0575 &0.0475 &0.0470 &0.0545\\
              &0.50         &0.1455 &0.3390 &0.7415 &0.9585\\
\hline
$ICM_n$       &0.00         &0.0010 &0.0000 &0.0000 &0.0000\\
              &0.50         &0.0250 &0.0000 &0.0000 &0.0000\\
\hline
$T_n^{GWZ}$   &0.00         &0.0410 &0.0440 &0.0455 &0.0400\\
              &0.50         &0.5645 &0.9090 &0.9965 &1.0000\\
\hline
$T_n^{ZH}$    &0.00         &0.0385 &0.0325 &0.0200 &0.0195\\
              &0.50         &0.0715 &0.0570 &0.0420 &0.0330\\
\hline
\end{tabular}}
\end{table}

\begin{table}[ht!]\caption{Empirical sizes and powers of the tests for $H_0$ vs. $H_{13}$.}
\centering
{\small\scriptsize\hspace{12.5cm}
\renewcommand{\arraystretch}{1}\tabcolsep 0.5cm
\begin{tabular}{*{20}{c}}
\hline
&\multicolumn{1}{c}{a}&\multicolumn{1}{c}{n=100}&\multicolumn{1}{c}{n=200}&\multicolumn{1}{c}{n=400}&\multicolumn{1}{c}{n=600}\\
&&\multicolumn{1}{c}{p=8}& \multicolumn{1}{c}{p=12}  &\multicolumn{1}{c}{p=17}&\multicolumn{1}{c}{p=20}\\
\hline
$AICM_n $   &0.00            &0.0510  &0.0500  &0.0500  &0.0650\\
               &0.25            &0.5620  &0.8530  &0.9920  &1.0000\\
\hline
$ACM_n$        &0.00            &0.0480  &0.0555  &0.0440  &0.0495\\
               &0.25            &0.5830  &0.8965  &0.9955  &1.0000\\
\hline
$PCvM_n$       &0.00            &0.0800  &0.0570  &0.0650  &0.0600\\
               &0.25            &0.6090  &0.9020  &0.9960  &1.0000\\
\hline
$T_n^{SZ}$     &0.00            &0.0490  &0.0525  &0.0590  &0.0495\\
               &0.25            &0.5820  &0.8885  &0.9960  &1.0000\\
\hline
$ICM_n$        &0.00            &0.0020  &0.0000  &0.0000  &0.0000\\
               &0.25            &0.0120  &0.0000  &0.0000  &0.0000\\
\hline
$T_n^{GWZ}$    &0.00            &0.0545  &0.0550  &0.0510  &0.0550\\
               &0.25            &0.3735  &0.6900  &0.9500  &0.9960\\
\hline
$T_n^{ZH}$     &0.00            &0.0245  &0.0290  &0.0235  &0.0180\\
               &0.25            &0.0525  &0.0405  &0.0270  &0.0255\\
\hline
\end{tabular}}
\end{table}

\newpage
\begin{table}[ht!]\caption{Empirical sizes and powers of the tests for $H_0$ vs. $H_{14}$.}
\centering
{\small\scriptsize\hspace{12.5cm}
\renewcommand{\arraystretch}{1}\tabcolsep 0.5cm
\begin{tabular}{*{20}{c}}
\hline
&\multicolumn{1}{c}{a}&\multicolumn{1}{c}{n=100}&\multicolumn{1}{c}{n=200}&\multicolumn{1}{c}{n=400}&\multicolumn{1}{c}{n=600}\\
&&\multicolumn{1}{c}{p=8}& \multicolumn{1}{c}{p=12}  & \multicolumn{1}{c}{p=17}&\multicolumn{1}{c}{p=20}\\
\hline
$AICM_n$      &0.0          &0.0640  &0.0520  &0.0540  &0.0500\\
                 &0.1          &0.3360  &0.5920  &0.8630  &0.9670\\
\hline
$ACM_n$          &0.0          &0.0505  &0.0470  &0.0475  &0.0575\\
                 &0.1          &0.3490  &0.5900  &0.8975  &0.9685\\
\hline
$PCvM_n$         &0.0          &0.0470  &0.0590  &0.0560  &0.0530\\
                 &0.1          &0.3600  &0.6510  &0.9010  &0.9780\\
\hline
$T_n^{SZ}$       &0.0          &0.0515  &0.0530  &0.0510  &0.0530\\
                 &0.1          &0.3325  &0.5830  &0.8860  &0.9765\\
\hline
$ICM_n$          &0.0          &0.0000  &0.0000  &0.0000  &0.0000\\
                 &0.1          &0.0170  &0.0000  &0.0000  &0.0000\\
\hline
$T_n^{GWZ}$      &0.0          &0.0465  &0.0505  &0.0560  &0.0515\\
                 &0.1          &0.1955  &0.3720  &0.6610  &0.8395\\
\hline
$T_n^{ZH}$       &0.0          &0.0280  &0.0225  &0.0170  &0.0200\\
                 &0.1          &0.0335  &0.0285  &0.0290  &0.0175\\
\hline
\end{tabular}}
\end{table}

\begin{table}[ht!]\caption{Empirical sizes and powers of the tests $AICM_n,  PCvM_n,  ICM_n,$ and $T_n^{ZH}$  for $H_{21}$.}
\centering
{\small\scriptsize\hspace{12.5cm}
\renewcommand{\arraystretch}{1}\tabcolsep 0.5cm
\begin{tabular}{*{20}{c}}
\hline
&\multicolumn{1}{c}{a}&\multicolumn{1}{c}{n=100}&\multicolumn{1}{c}{n=200}&\multicolumn{1}{c}{n=400}&\multicolumn{1}{c}{n=600}\\
&&\multicolumn{1}{c}{p=8}& \multicolumn{1}{c}{p=12}  & \multicolumn{1}{c}{p=17}&\multicolumn{1}{c}{p=20}\\
\hline
$AICM_n$       &0.0      &0.0650 &0.0640 &0.0580 &0.0650\\
                  &0.5      &0.2850 &0.6390 &0.9780 &1.0000\\
\hline
$PCvM_n$          &0.0      &0.0680 &0.0750 &0.0680 &0.0760\\
                  &0.5      &0.2390 &0.5030 &0.8620 &0.9820\\
\hline
$ICM_n$           &0.0      &0.0170 &0.0000 &0.0000 &0.0000\\
                  &0.5      &0.0460 &0.0000 &0.0000 &0.0000\\
\hline
$T_n^{ZH}$        &0.0      &0.0215 &0.0220 &0.0270 &0.0170\\
                  &0.5      &0.0405 &0.0335 &0.0325 &0.0345\\
\hline
\end{tabular}}
\end{table}

\begin{table}[ht!]\caption{Empirical sizes and powers of the tests $AICM_n,  PCvM_n,  ICM_n,$ and $T_n^{ZH}$ for $H_{22}$.}
\centering
{\small\scriptsize\hspace{12.5cm}
\renewcommand{\arraystretch}{1}\tabcolsep 0.5cm
\begin{tabular}{*{20}{c}}
\hline
&\multicolumn{1}{c}{a}&\multicolumn{1}{c}{n=100}&\multicolumn{1}{c}{n=200}&\multicolumn{1}{c}{n=400}&\multicolumn{1}{c}{n=600}\\
&&\multicolumn{1}{c}{p=8}& \multicolumn{1}{c}{p=12}  & \multicolumn{1}{c}{p=17}&\multicolumn{1}{c}{p=20}\\
\hline
$AICM_n$       &0.0          &0.0700   &0.0790   &0.0560   &0.0720\\
                  &0.5          &0.4030   &0.7780   &0.9880   &1.0000\\
\hline
$PCvM_n$          &0.0          &0.0880   &0.0770   &0.0610   &0.0500\\
                  &0.5          &0.3080   &0.5420   &0.8100   &0.9110\\
\hline
$ICM_n$           &0.0          &0.0170   &0.0000   &0.0000   &0.0000\\
                  &0.5          &0.0880   &0.0000   &0.0000   &0.0000\\
\hline
$T_n^{ZH}$        &0.0          &0.0230   &0.0230   &0.0240   &0.0220\\
                  &0.5          &0.0565   &0.0500   &0.0525   &0.0335\\
\hline
\end{tabular}}
\end{table}

\begin{table}[ht!]\caption{Empirical sizes and powers of the tests $AICM_n,  PCvM_n,  ICM_n,$ and $T_n^{ZH}$ for $H_{23}$.}
\centering
{\small\scriptsize\hspace{12.5cm}
\renewcommand{\arraystretch}{1}\tabcolsep 0.5cm
\begin{tabular}{*{20}{c}}
\hline
&\multicolumn{1}{c}{a}&\multicolumn{1}{c}{n=100}&\multicolumn{1}{c}{n=200}&\multicolumn{1}{c}{n=400}&\multicolumn{1}{c}{n=600}\\
&&\multicolumn{1}{c}{p=8}& \multicolumn{1}{c}{p=12}  & \multicolumn{1}{c}{p=17}&\multicolumn{1}{c}{p=20}\\
\hline
$AICM_n$       &0.00     &0.0790  &0.0740  &0.0700  &0.0610\\
               &0.75     &0.3060  &0.5250  &0.8390  &0.9680\\
\hline
$PCvM_n$       &0.00     &0.0690  &0.0730  &0.0600  &0.0630\\
               &0.75     &0.2170  &0.2820  &0.3320  &0.3730\\
\hline
$ICM_n$        &0.00     &0.0090  &0.0000  &0.0000  &0.0000\\
               &0.75     &0.0230  &0.0000  &0.0000  &0.0000\\
\hline
$T_n^{ZH}$     &0.00     &0.0270  &0.0270  &0.0225  &0.0225\\
               &0.75     &0.0380  &0.0365  &0.0375  &0.0280\\
\hline
\end{tabular}}
\end{table}

\newpage

\begin{figure}
  \centering
  \includegraphics[width=10cm,height=8cm]{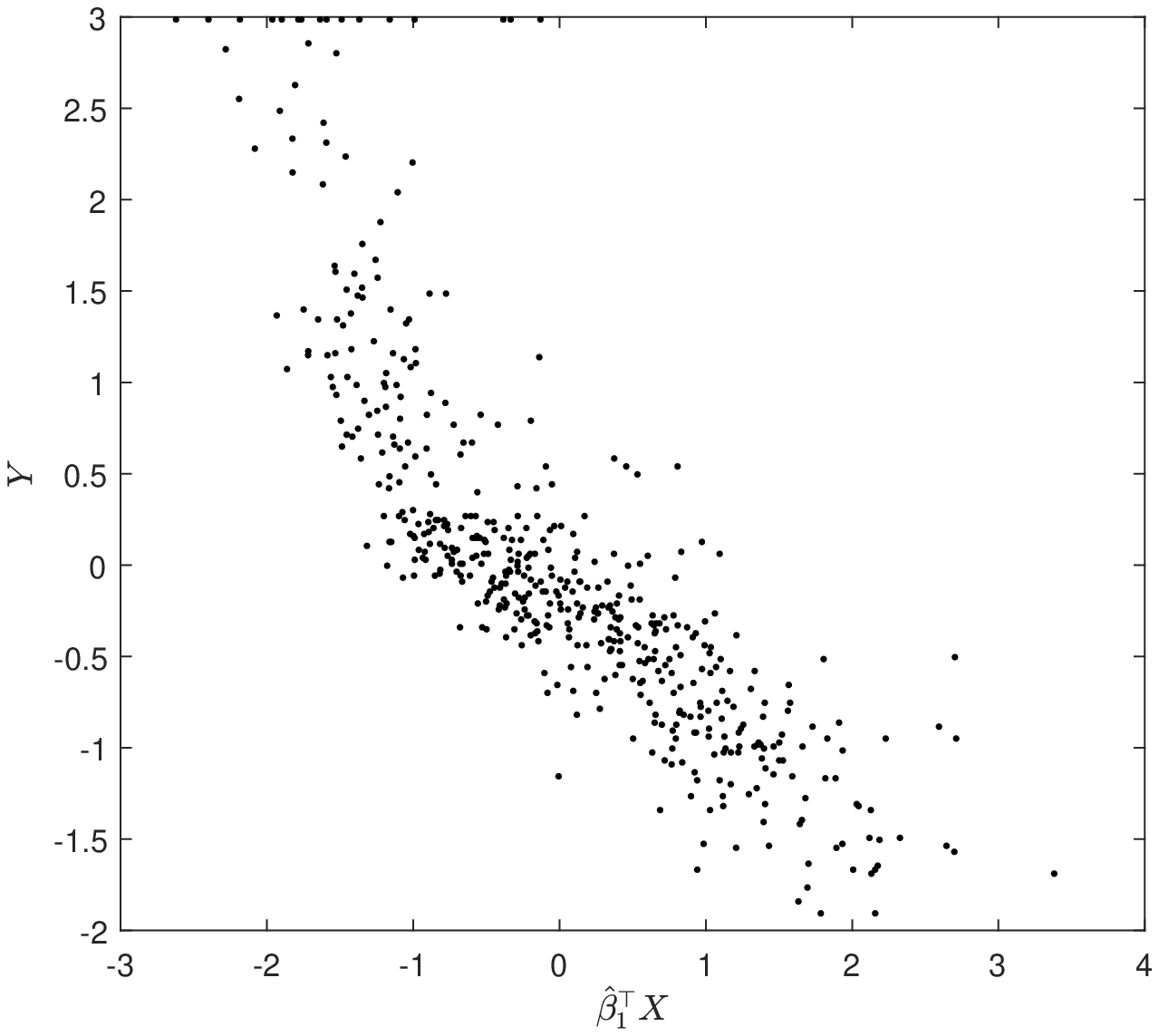}
  \caption{Scatter plot of $Y$ versus the projected covariate $\hat{\beta}_1^{\top}X$ in which the direction $\hat{\beta}_1$ is obtained by CSE.}\label{Figure2}
\end{figure}

\begin{figure}
  \centering
  \includegraphics[width=10cm,height=8cm]{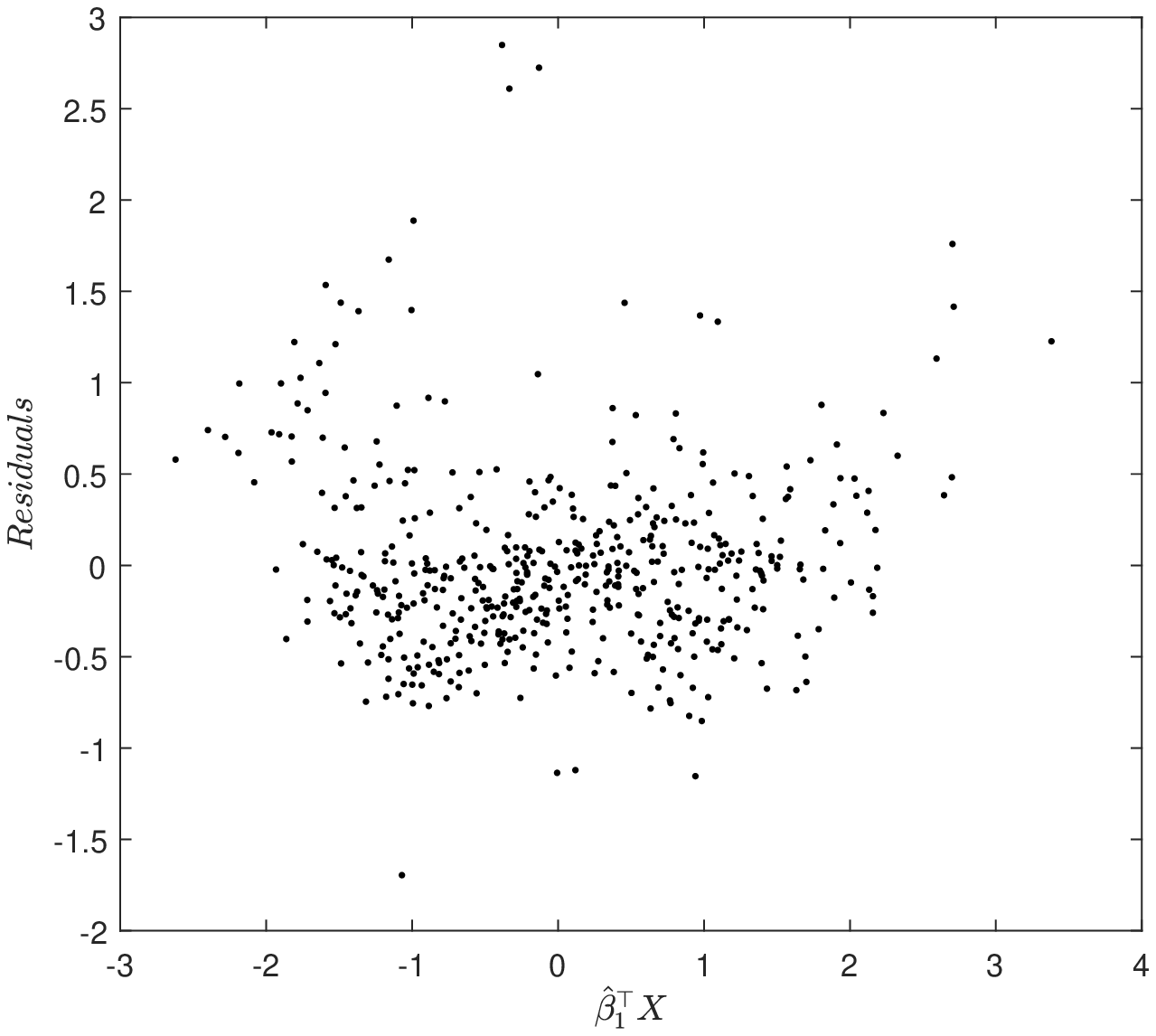}
  \caption{Scatter plot of Residuals from the linear model for $(Y,X)$ versus the projected covariate $\hat{\beta}_1$ in which the direction $\hat{\beta}_1^{\top}X$ is obtained by CSE.}\label{Figure2}
\end{figure}

\begin{figure}
  \centering
  \includegraphics[width=10cm,height=10cm]{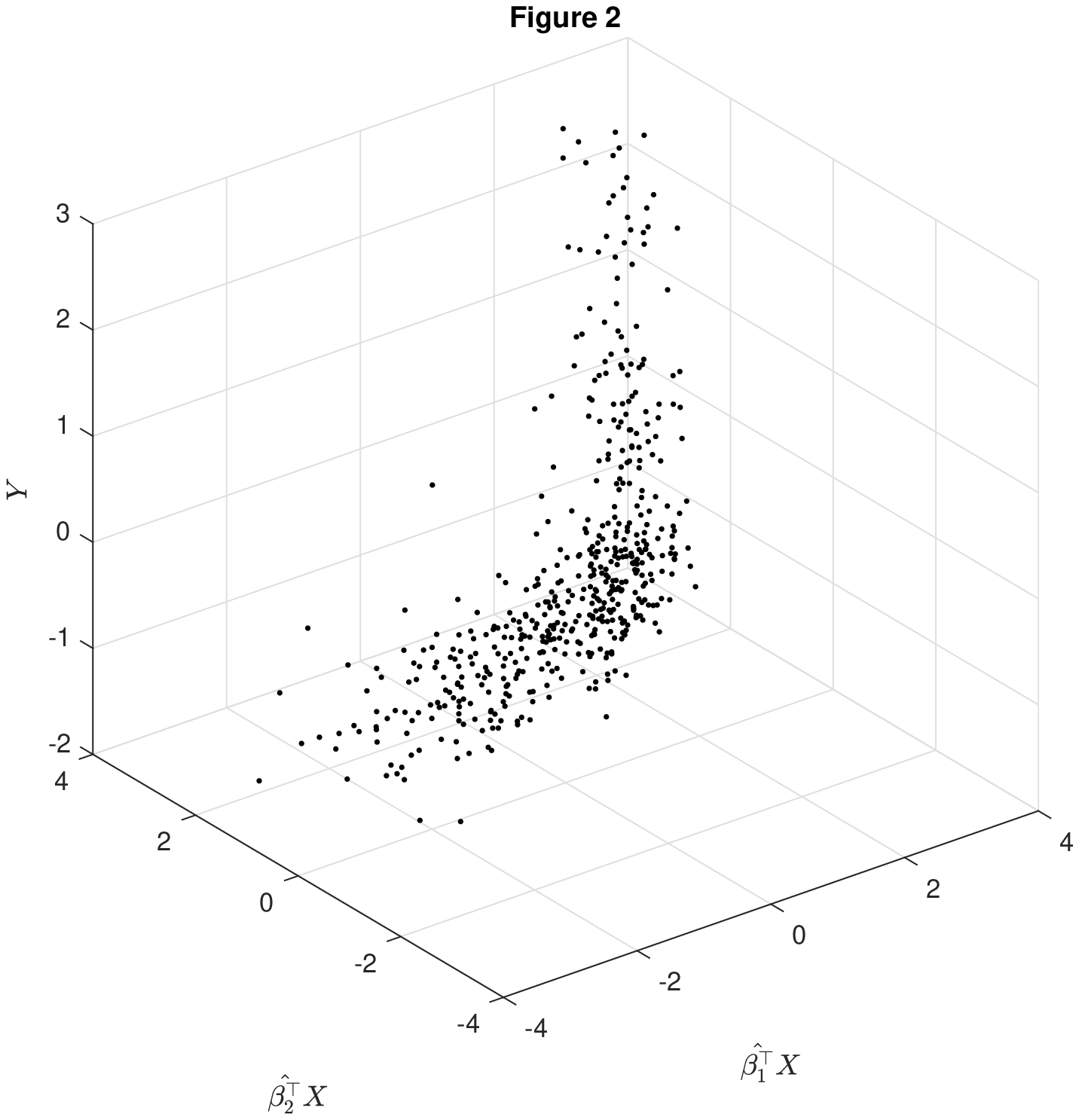}
  \caption{Scatter plot of $Y$ versus the projected covariates $(\hat{\beta}_1^{\top}X, \hat{\beta}_2^{\top}X)$ in which the directions $(\hat{\beta}_1, \hat{\beta}_2)$ are obtained by CSE.}\label{Figure3}
\end{figure}

\begin{figure}
  \centering
  \includegraphics[width=10cm,height=8cm]{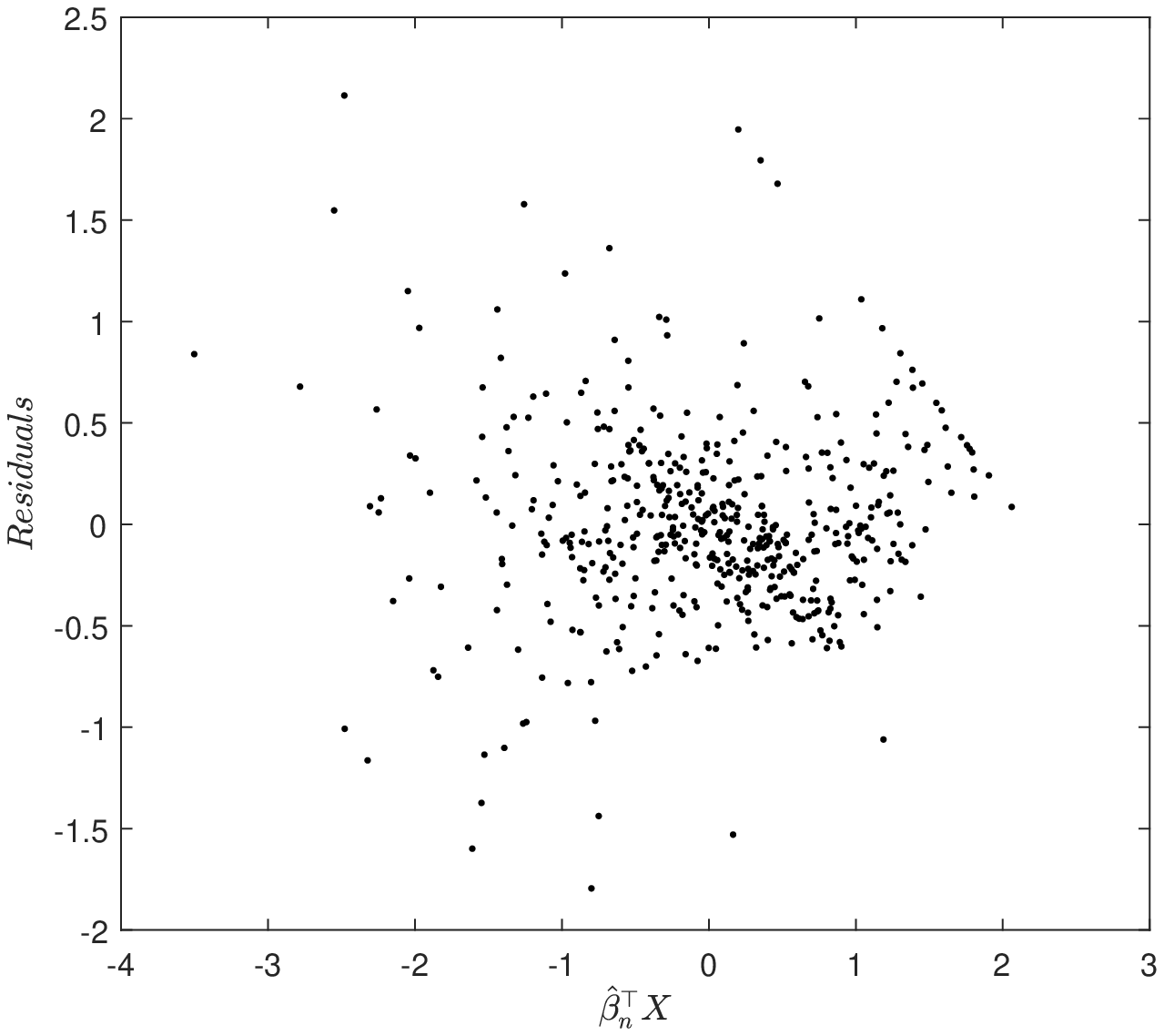}
  \caption{Scatter plot of Residuals from model (\ref{6.1}) versus the fitted values $\hat{\beta}_n^{\top}X$.}\label{Figure4}
\end{figure}

\end{document}